\documentclass{amsart}
\def\BookNumber{1001.4852}
\def\Defined{}
\ifx\UseRussian\Defined
	\usepackage[T2A,T2B]{fontenc}
	\usepackage[cp1251]{inputenc}
	\usepackage[english,russian]{babel}
\fi
\usepackage{footmisc}
\usepackage[all]{xy}
\usepackage{color}
\definecolor{UrlColor}{rgb}{.9,0,.3}
\definecolor{SymbColor}{rgb}{.4,0,.9}
\definecolor{IndexColor}{rgb}{1,.3,.6}
\definecolor{eml1}{rgb}{.8,.1,.1}
\definecolor{eml2}{rgb}{.1,.6,.6}

\usepackage{xr-hyper}
\usepackage[unicode]{hyperref}
\hypersetup{pdfdisplaydoctitle=true}
\hypersetup{colorlinks}
\hypersetup{citecolor=UrlColor}
\hypersetup{urlcolor=UrlColor}
\hypersetup{pdffitwindow=true}
\hypersetup{pdfnewwindow=true}
\hypersetup{pdfstartview={FitH}}

\def\hyph{\penalty0\hskip0pt\relax-\penalty0\hskip0pt\relax}
\def\Hyph{-\penalty0\hskip0pt\relax}

\newcommand{\Basis}[1]{\overline{\overline{#1}}{}}
\newcommand{\Vector}[1]{\overline{#1}{}}
\newcommand{\gi}[1]{\boldsymbol{\textcolor{IndexColor}{#1}}}

\makeatletter
\newcommand{\NameDef}[1]{%
	\expandafter\gdef\csname #1\endcsname%
}%
\newcommand{\ShowSymbol}[1]{%
	\@nameuse{ViewSymbol#1}%
}%
\newcommand{\symb}[3]{%
	\@ifundefined{ViewSymbol#3}{%
		\NameDef{ViewSymbol#3}{\textcolor{SymbColor}{#1}}%
		\NameDef{RefSymbol#3}{\pageref{symbol: #3}}%
		\@namedef{LabelSymbol#3}{\label{symbol: #3}}%
	}{%
		\message {extra symb #3}
	}%
	\ifcase#2%0
	\or%1
		$\@nameuse{ViewSymbol#3}$%
	\or%2
		\[\@nameuse{ViewSymbol#3}\]%
	\else%
	\fi%
	\@nameuse{LabelSymbol#3}%
}%
\newcommand{\DefEq}[2]{%
	\@ifundefined{ViewEq#2}{%
		\NameDef{ViewEq#2}{#1}%
	}{%
	}%
}%
\newcommand{\ShowEq}[1]{%
	\@ifundefined{ViewEq#1}{%
		\message {missed ShowEq #1}
  }{%
	\@nameuse{ViewEq#1}%
	}%
}%
\makeatother

\newcommand{\subs}{${}_*$\Hyph}
\newcommand{\sups}{${}^*$\Hyph}

\newcommand{\CRstar}{{}_*{}^*}
\newcommand{\RCstar}{{}^*{}_*}
\newcommand{\RCcirc}{{}_{\circ}{}^{\circ}}
\newcommand{\CRcirc}{{}^{\circ}{}_{\circ}}

\newcommand{\RC}{$\RCstar$\Hyph}
\newcommand{\CR}{$\CRstar$\Hyph}
\newcommand{\drc}{$D\RCstar$\Hyph}
\newcommand{\Drc}{$\mathcal D\RCstar$\Hyph}
\newcommand{\dcr}{$D\CRstar$\hyph}
\newcommand{\rcd}{$\RCstar D$\Hyph}
\newcommand{\crd}{$\CRstar D$\Hyph}

\newcommand\sT{$\star T$\Hyph}%
\newcommand\Ts{$T\star$\Hyph}%
\newcommand\sD{$\star D$\Hyph}%
\newcommand\Ds{$D\star$\Hyph}%

\newcommand\VirtVar{\vphantom{\overset{\rightarrow}{1}^1}}
\newcommand\pC[2]{{}_{(#1)#2}}%
\newcommand\DrcPartial[1]%
{%
	\def\tempa{}%
	\def\tempb{#1}%
	\ifx\tempa\tempc%
		(\partial\RCstar)%
	\else%
		({}_{#1}\partial\RCstar)%
	\fi%
}%
\newcommand\crDPartial[1]%
{%
	\def\tempa{}%
	\def\tempb{#1}%
	\ifx\tempa\tempc%
		(\CRstar\partial)%
	\else%
		(\CRstar\partial_{#1})%
	\fi%
}%
\newcommand\StandPartial[3]%
{%
	\frac{\partial^{\gi{#3}} #1}{\partial #2}%
}%

\renewcommand{\uppercasenonmath}[1]{}

\makeatletter
\newcommand\@dotsep{4.5}
\def\@tocline#1#2#3#4#5#6#7
{\relax
		\par \addpenalty\@secpenalty\addvspace{#2}%
		\begingroup %\hyphenpenalty\@M
		\@ifempty{#4}{%
			\@tempdima\csname r@tocindent\number#1\endcsname\relax
		}{%
			\@tempdima#4\relax
		}%
		\parindent\z@ \leftskip#3\relax \advance\leftskip\@tempdima\relax
		\rightskip\@pnumwidth plus1em \parfillskip-\@pnumwidth
		#5\leavevmode\hskip-\@tempdima #6\relax
		\leaders\hbox{$\m@th
		\mkern \@dotsep mu\hbox{.}\mkern \@dotsep mu$}\hfill
		\hbox to\@pnumwidth{\@tocpagenum{#7}}\par
		\nobreak
		\endgroup
}
\makeatother 

\ifx\PrintBook\undefined
	\makeatletter
	\renewcommand{\@indextitlestyle}{%
		\twocolumn[\section{\indexname}]%
		\def\IndexSpace{off}%
	}
	\makeatother 
	\thanks{\href{mailto:Aleks\_Kleyn@MailAPS.org}{Aleks\_Kleyn@MailAPS.org}}
	\thanks{\ \ \ \url{http://sites.google.com/site/AleksKleyn/}}
	\thanks{\ \ \ \url{http://arxiv.org/a/kleyn\_a\_1}}
	\thanks{\ \ \ \url{http://AleksKleyn.blogspot.com/}}
\else

\makeatletter
\def\@maketitle{%
  \cleardoublepage \thispagestyle{empty}%
  \begingroup \topskip\z@skip
  \null\vfil
  \begingroup
  \LARGE\bfseries \centering
  \openup\medskipamount
  \@title\par\vspace{24pt}%
  \def\and{\par\medskip}\centering
  \mdseries\authors\par\bigskip
  \endgroup
  \vfil\vspace{24pt}
  \ifx\@empty\addresses \else \@setaddresses \fi
  \vfil
  \ifx\@empty\@dedicatory
  \else \begingroup
    \centering{\footnotesize\itshape\@dedicatory\@@par}%
    \endgroup
  \fi
  \vfill
  \newpage\thispagestyle{empty}
  \@setabstract
  \begin{center}
    \ifx\@empty\@subjclass\else\@setsubjclass\fi
    \ifx\@empty\@keywords\else\@setkeywords\fi
    \ifx\@empty\@translators\else\vfil\@settranslators\fi
    \ifx\@empty\thankses\else\vfil\@setthanks\fi
  \end{center}
  \vfil
  \endgroup}
\makeatother 

	\makeatletter
	\renewcommand{\@indextitlestyle}{%
		\twocolumn[\chapter{\indexname}]%
		\def\IndexSpace{off}%
		\let\@secnumber\@empty
		\chaptermark{\indexname}%
	}
	\makeatother 
	\email{\href{mailto:Aleks\_Kleyn@MailAPS.org}{Aleks\_Kleyn@MailAPS.org}}
	\urladdr{\url{http://sites.google.com/site/alekskleyn/}}
	\urladdr{\url{http://arxiv.org/a/kleyn\_a\_1}}
	\urladdr{\url{http://AleksKleyn.blogspot.com/}}
\fi

\ifx\SelectlEnglish\undefined
	\ifx\UseRussian\undefined
		\def\SelectlEnglish{}
	\fi
\fi

\newcommand\wRefDef[2]
	{
	\def\Tempa{#1}
	\def\Tempb{0405.027}
	\ifx\Tempa\Tempb
	\def\wRef{gr-qc/pdf/0405/0405027v3.pdf}
	\fi
	\def\Tempb{0405.028}
	\ifx\Tempa\Tempb
	\def\wRef{gr-qc/pdf/0405/0405028v5.pdf}
	\fi
	\def\Tempb{0412.391}
	\ifx\Tempa\Tempb
	\def\wRef{math/pdf/0412/0412391v4.pdf}
	\fi
	\def\Tempb{0612.111}
	\ifx\Tempa\Tempb
	\def\wRef{math/pdf/0612/0612111v2.pdf}
	\fi
	\def\Tempb{0701.238}
	\ifx\Tempa\Tempb
	\def\wRef{math/pdf/0701/0701238v4.pdf}
	\fi
	\def\Tempb{0702.561}
	\ifx\Tempa\Tempb
	\def\wRef{math/pdf/0702/0702561v3.pdf}
	\fi
	\def\Tempb{0707.2246}
	\ifx\Tempa\Tempb
	\def\wRef{arxiv/pdf/0707/0707.2246v2.pdf}
	\fi
	\def\Tempb{0803.3276}
	\ifx\Tempa\Tempb
	\def\wRef{arxiv/pdf/0803/0803.3276v3.pdf}
	\fi
	\def\Tempb{0812.4763}
	\ifx\Tempa\Tempb
	\def\wRef{arxiv/pdf/0812/0812.4763v6.pdf}
	\fi
 	\def\Tempb{0906.0135}
	\ifx\Tempa\Tempb
	\def\wRef{arxiv/pdf/0906/0906.0135v3.pdf}
	\fi
 	\def\Tempb{0912.4061}
	\ifx\Tempa\Tempb
	\def\wRef{arxiv/pdf/0912/0912.4061v2.pdf}
	\fi
	\externaldocument[#1-#2-]{#1.#2}[http://arxiv.org/PS_cache/\wRef]
	}
\newcommand\LanguagePrefix{}%
\ifx\SelectlEnglish\undefined
	\newcommand\CurrentLanguage{Russian.}%
	\author{Александр Клейн}
	\newtheorem{theorem}{Теорема}[section]
	\newtheorem{corollary}[theorem]{Следствие}
	\theoremstyle{definition}
	\newtheorem{definition}[theorem]{Определение}
	\newtheorem{example}[theorem]{Пример}
	\newtheorem{xca}[theorem]{Exercise}
	\theoremstyle{remark}
	\newtheorem{remark}[theorem]{Замечание}
	
	\newcommand\Gbasis{$G$\Hyph базис}
	\newcommand\Gcoords{$G$\Hyph координат}
	\newcommand\Gspace{$G$\Hyph пространств}
	\newcommand\xRefDef[1]
		{
		\wRefDef{#1}{Russian}
		\NameDef{xRefDef#1}{}%
		}
	\makeatletter
	\newcommand\xRef[2]%
	{%
		\@ifundefined{xRefDef#1}{%
		\ref{#2}%
		}{% 
		\citeBib{#1}-\ref{#1-Russian-#2}%
		}%
	}%
	\newcommand\xEqRef[2]%
	{%
		\@ifundefined{xRefDef#1}{%
		\eqref{#2}%
		}{%
		\citeBib{#1}-\eqref{#1-Russian-#2}%
		}%
	}%
	\makeatother
	\ifx\PrintBook\undefined
		\newcommand{\BibTitle}{%
			\section{Список литературы}%
		}
	\else
		\newcommand{\BibTitle}{%
			\chapter{Список литературы}%
		}
	\fi
\else
	\newcommand\CurrentLanguage{English.}%
	\author{Aleks Kleyn}
	\newtheorem{theorem}{Theorem}[section]
	\newtheorem{corollary}[theorem]{Corollary}
	\theoremstyle{definition}
	\newtheorem{definition}[theorem]{Definition}

	\theoremstyle{remark}
	
	\newcommand\Gbasis{$G$\Hyph basis}
	\newcommand\Gcoords{$G$\Hyph coordinates}
	\newcommand\Gspace{$G$\Hyph space}
	\newcommand\xRefDef[1]
		{
		\wRefDef{#1}{English}
		\NameDef{xRefDef#1}{}%
		}
	\makeatletter
	\newcommand\xRef[2]%
	{%
		\@ifundefined{xRefDef#1}{%
		\ref{#2}%
		}{% 
		\citeBib{#1}-\ref{#1-English-#2}%
		}%
	}%
	\newcommand\xEqRef[2]%
	{%
		\@ifundefined{xRefDef#1}{%
		\eqref{#2}%
		}{%
		\citeBib{#1}-\eqref{#1-English-#2}%
		}%
	}%
	\makeatother
	\ifx\PrintBook\undefined
		\newcommand{\BibTitle}{%
			\section{References}%
		}
	\else
		\newcommand{\BibTitle}{%
			\chapter{References}%
		}
	\fi
\fi

\ifx\PrintBook\undefined
	\numberwithin{Hfootnote}{section}
\else
	\numberwithin{section}{chapter}
	\numberwithin{footnote}{chapter}
	\numberwithin{Hfootnote}{chapter}
\fi

\numberwithin{equation}{section}
\numberwithin{figure}{section}
\numberwithin{table}{section}
\numberwithin{Item}{section}

\makeatletter
\newcommand\org@maketitle{}
\let\org@maketitle\maketitle
\def\maketitle{%
	\hypersetup{pdftitle={\@title}}%
	\hypersetup{pdfauthor={\authors}}%
	\hypersetup{pdfsubject=\@keywords}%
	\org@maketitle
}
\def\make@stripped@name#1{%
	\begingroup
		\escapechar\m@ne
		\global\let\newname\@empty
		\protected@edef\Hy@tempa{\CurrentLanguage #1}%
		\edef\@tempb{%
			\noexpand\@tfor\noexpand\Hy@tempa:=%
			\expandafter\strip@prefix\meaning\Hy@tempa
		}%
		\@tempb\do{%
			\if\Hy@tempa\else
				\if\Hy@tempa\else
					\xdef\newname{\newname\Hy@tempa}%
				\fi
			\fi
		}%
	\endgroup
}%
\newenvironment{enumBib}{%
	\BibTitle
	\advance\@enumdepth \@ne
	\edef\@enumctr{enum\romannumeral\the\@enumdepth}\list
	{\csname biblabel\@enumctr\endcsname}{\usecounter
	{\@enumctr}\def\makelabel##1{\hss\llap{\upshape##1}}}
}{%
	\endlist
}
\def\Items#1{\ItemList#1,LastItem,}%
\def\LastItem{LastItem}%
\def\ItemList#1,{\def\ViewBook{#1}%
	\ifx\ViewBook\LastItem%
	\else%
		\ifx\ViewBook\BookNumber%
			\def\Semafor{on}%
		\fi%
		\expandafter\ItemList%
	\fi%
}%
\makeatletter

\newcommand{\ePrints}[1]
{%
	\def\Semafor{off}%
	\Items{#1}%
}%

\newcommand{\BiblioItem}[2]
{
	\def\Semafor{off}
	\@ifundefined{\LanguagePrefix ViewCite#1}{}{%
		\def\Semafor{on}%
	}%
	\ifx\Semafor\ValueOff
		\@ifundefined{xRefDef#1}{}{% 
		\def\Semafor{on}%
		}%
	\fi
	\ifx\Semafor\ValueOn
		\ifx\IndexState\ValueOff
			\begin{enumBib}
			\def\IndexState{on}
		\fi
		\item \label{\LanguagePrefix bibitem: #1}#2%
	\fi
}
\makeatother
\newcommand{\OpenBiblio}
{
	\def\IndexState{off}
}
\newcommand{\CloseBiblio}
{
	\ifx\IndexState\ValueOn
		\end{enumBib}
		\def\IndexState{off}
	\fi
}

\makeatletter
\def\StartCite{[}%
\def\citeBib#1{\protect\showCiteBib#1,endCite,}%
\def\endCite{endCite}%
\def\showCiteBib#1,{\def\temp{#1}%
\ifx\temp\endCite
]%
\def\StartCite{[}%
\else
	\StartCite\LanguagePrefix \ref{\LanguagePrefix bibitem: #1}%
	\@ifundefined{\LanguagePrefix ViewCite#1}{%
		\NameDef{\LanguagePrefix ViewCite#1}{}%
	}{%
	}%
	\def\StartCite{, }%
\expandafter\showCiteBib%
\fi}%
\makeatother

\newcommand{\arp}{\ar @{-->}}

\newcommand{\bundle}[4]%
{%
	\def\tempa{}%
	\def\tempb{#3}%
	\def\tempc{#1}%
	\ifx\tempa\tempb%
		\ifx\tempa\tempc%
			#2%
		\else%
			\xymatrix{#2:#1\arp[r]&#4}%
		\fi%
	\else%
		\ifx\tempa\tempc%
			#2[#3]%
		\else%
			\xymatrix{#2[#3]:#1\arp[r]&#4}%
		\fi%
	\fi%
}%
\makeatletter
\newcommand{\AddIndex}[2]%
{%
	\@ifundefined{RefIndex#2}{%
		\NameDef{RefIndex#2}{#2}%
	}{%
		\message {extra AddIndex #2}
	}%
	{\bf #1}%
	\label{index: #2}%
}%
\newcommand{\Index}[2]%
{%
	\def\Semafor{off}%
	\@ifundefined{RefIndex#2}{%
	}{%
		\def\Semafor{on}
	}%
	\ifx\Semafor\ValueOn%
		\def\tempa{}%
		\def\tempb{#2}%
		\ifx\IndexState\ValueOff%
			\begin{theindex}%
			\def\IndexState{on}%
		\fi%
		\ifx\IndexSpace\ValueOn%
			\indexspace%
			\def\IndexSpace{off}%
		\fi%
		\item #1%
		\ifx\tempa\tempb%
		\else%
			\ \pageref{index: #2}%
		\fi%
	\fi%
}%
\newcommand{\Symb}[2]
{
	\def\Semafor{off}
	\@ifundefined{ViewSymbol#2}{%
	}{%
		\def\Semafor{on}
	}%
	\ifx\Semafor\ValueOn
		\ifx\IndexState\ValueOff
			\begin{theindex}
			\def\IndexState{on}
		\fi
		\ifx\IndexSpace\ValueOn
			\indexspace
			\def\IndexSpace{off}
		\fi
		\item $\displaystyle\@nameuse{ViewSymbol#2}$\ \ #1
		\@nameuse{RefSymbol#2}%
	\fi
}

\makeatother
\newcommand{\SetIndexSpace}%
{%
	\def\IndexSpace{on}%
}%
\def\ValueOff{off}
\def\ValueOn{on}

\newcommand{\OpenIndex}
{
	\def\IndexState{off}
}
\newcommand{\CloseIndex}
{
	\ifx\IndexState\ValueOn
		\end{theindex}
		\def\IndexState{off}
	\fi
}

\def\LastMemo{LastMemo}%
\def\MemoList#1//{\def\temp{#1}%
	\ifx\temp\LastMemo
	\else%
		\par
		\textcolor{blue}{#1}%
		\expandafter\MemoList%
	\fi%
}%

\xRefDef{0701.238}
\xRefDef{0812.4763}
\xRefDef{0906.0135}
\xRefDef{0912.4061}
\begin{document}
\pdfbookmark[1]{The Matrix of Linear Mappings}{TitleEnglish}
\title{The Matrix of Linear Mappings}

\begin{abstract}
On the set of mappings of the given set,
we define the product of mappings.
If $A$ is associative algebra,
then we consider the set of matrices, whose elements
are linear mappings of algebra $A$.
In algebra of matrices of linear mappings we
define the operation of $\RCcirc$\Hyph product.
The operation is based on the product of mappings.

If the matrix $a$ of linear mappings
has an inverse matrix, then the quasideterminant of the matrix $a$
and the inverse matrix are matrices of linear mappings.
In the paper, I consider conditions when a matrix of linear mappings
has inverse matrix, as well methods of solving a
system of linear equations in an associative algebra.
\end{abstract}
\maketitle

\tableofcontents

%auto-ignore
%auto-ignore

\DefEq
{
\begin{equation}
\begin{pmatrix}
\displaystyle\frac{\partial f^1(\Vector x)}{\partial x^1}(h^1)
&...&
\displaystyle\frac{\partial f^m(\Vector x)}{\partial x^1}(h^1)
\\...&...&...\\
\displaystyle\frac{\partial f^1(\Vector x)}{\partial x^n}(h^n)
&...&
\displaystyle\frac{\partial f^m(\Vector x)}{\partial x^n}(h^n)
\end{pmatrix}
=
\begin{pmatrix}
\displaystyle\frac{\partial f^1(\Vector x)}{\partial x^1}
&...&
\displaystyle\frac{\partial f^m(\Vector x)}{\partial x^1}
\\...&...&...\\
\displaystyle\frac{\partial f^1(\Vector x)}{\partial x^n}
&...&
\displaystyle\frac{\partial f^m(\Vector x)}{\partial x^n}
\end{pmatrix}
\begin{pmatrix}
h^1\\...\\h^n
\end{pmatrix}
\label{eq: Jacobi Gateaux matrix of map, vector space}
\end{equation}
}
{Jacobi Gateaux matrix of map, vector space}

\DefEq
{
\begin{equation}
\label{eq: differential of composite map, D vector space}
\partial \Vector g\Vector f(\Vector x)(\Vector a)
=\partial \Vector g(\Vector f(\Vector x))(\partial \Vector f(\Vector x)(\Vector a))
\end{equation}
}
{differential of composite map, D vector space}

\DefEq
{
\begin{equation}
\Vector x'_i(\Vector x')(v'^i)
=\Vector x_j(\Vector x)
\left(
\frac{\partial x^j}{\partial x'^i}(v'^i)
\right)
\label{eq: map and coordinates, 5, D affine space}
\end{equation}
}
{map and coordinates, D affine space}

\DefEq
{
\[
f(x)=f\circ x
\]
}
{new notation for mapping}

\DefEq
{
\[
(a^{\gi{ij}}\Vector e_{\gi i}\otimes \Vector e_{\gi j})\circ x=
a^{\gi{ij}}\Vector e_{\gi i}x\Vector e_{\gi j}
\]
}
{linear mapping of division ring}

\DefEq
{
\begin{equation}
\partial f(x)\circ a=
\left(\StandPartial{f(x)}{x}{ij}
\Vector e_{\gi i}\otimes\Vector e_{\gi j}\right)\circ a
\label{eq: derivative of mapping f}
\end{equation}
}
{derivative of mapping f}

\section{Preface}

I gave the definition of tensor products of division rings
(section
\xRef{0701.238}{section: Tensor Product of Division Rings})
almost at the same time when I defined the linear mapping
of division rings (section
\xRef{0701.238}{section: Additive Map of Division Ring}).
It was evident for me that components of a linear mapping
are tensor of order $2$. However this statement was for
me so unexpected that until now I have written nothing
about this in my papers.

There was a contradiction in this statement.
I used a linear mapping to define a tensor
of order $1$. However to build a mapping I was need
a tensor of order $2$; in other words, I was need a bilinear mapping.

However, I do not see any contradiction now.
For any field $F$, tensor product $F\otimes F$
is isomorphic to $F$.
The consequence of this isomorphism is the possibility of replacing
tensor $a\otimes b$ by ordinary product $ab$.
Therefore, the set of linear mappings of the field $F$ is isomorphic
to the field $F$ or, to be more exact,
to the set of left shifts of the multiplicative
group of the field $F$.
In division ring, the product is noncomutative
and there is no such isomorphism.

Another problem that I met arose when I
wanted to represent a derivative of a mapping of vector space
as a matrix of partial derivatives (the equation
\xEqRef{0812.4763}{eq: Jacobi Gateaux matrix of map, D vector space}).
One may be tempted to write the Jacobi\Hyph G\^ateaux matrix of mapping
as a product of matrices
\ShowEq{Jacobi Gateaux matrix of map, vector space}
However, it is evident that the equation
\eqref{eq: Jacobi Gateaux matrix of map, vector space}
is not true.
However the matrix written in the left side of the equation
\eqref{eq: Jacobi Gateaux matrix of map, vector space}
has very strong deficiency.
I have impression that rank of this matrix depends on increment of
argument.

We observe similar picture in the theorem
\xRef{0812.4763}{theorem: composite map, differential, D vector space},
where the derivative has form
\ShowEq{differential of composite map, D vector space}

The last reason to change model under study
was the rule of transformation of vector
in curvilinear coordinates of affine space
(equations
\xEqRef{0906.0135}{eq: map and coordinates, 5, D affine space},
\xEqRef{0906.0135}{eq: map and coordinates, 6, D affine space})
\ShowEq{map and coordinates, D affine space}
This finally convinced me that I had met new mathematical
object. The name of this object is functional matrix, in other words,
the matrix whose elements are mappings.

%The product of functional matrices is different from the product
%of ordinary matrices тем, что вместо произведения элемента одной матрицы
%на элемент другой I consider product of соответствующих
%mappings.
%To подчеркнуть this difference, I offered
%new notation for the product of functional matrices.

All that remains is the final step.
If the element of the matrix is the mapping,
then I break the connection between the mapping and the argument,
or more exactly, I write the mapping as an operator
\ShowEq{new notation for mapping}

Using new notation, I can write the linear mapping of division ring
in the following form
\ShowEq{linear mapping of division ring}
Correspondingly, we can write derivative of mapping $f$ in the form
\ShowEq{derivative of mapping f}
In the equation \eqref{eq: derivative of mapping f},
I succeeded in separation of derivative and increment.

%auto-ignore
%auto-ignore

\DefEq
{
\[
\Vector r(\Vector a)=(r_1(a_1),...,r_n(a_n))
\]
}
{vector notation in tower of representations}

\DefEq
{
\[
a=a^{\gi i}\Vector e_{\gi i}
\]
}
{Expansion relative basis in algebra}

\section{Conventions}

\begin{enumerate}

\ePrints{0812.4763,0906.0135,0908.3307,0909.0855,0912.3315}
\ifx\Semafor\ValueOn
\item Function and map are synonyms. However according to
tradition, correspondence between either rings or vector
spaces is called map and map of
either real field or quaternion algebra is called function.
I also follow this tradition.
%However we can consider free
%algebra as vector space over field.
%For the same correspondence it is undesirable называть то функцией, то
%отображением, я буду пользоваться термином линейная функция алгебры.
\fi

\ePrints{0701.238,0812.4763,0908.3307,0912.4061,Report.10,1001.4852}
\ifx\Semafor\ValueOn
\item
In any expression where we use index I assume
that this index may have internal structure.
For instance, considering the algebra $A$ we enumerate coordinates of
$a\in A$ relative to basis $\Basis e$ by an index $i$.
This means that $a$ is a vector. However, if $a$
is matrix, then we need two indexes, one enumerates
rows, another enumerates columns. In the case, when index has
structure, we begin the index from symbol $\cdot$ in
the corresponding position. 
For instance, if I consider the matrix $a^i_j$ as an element of a vector
space, then I can write the element of matrix as $a^{\cdot}{}^i_j$.
\fi

\ePrints{0701.238,0812.4763,0908.3307,0912.4061}
\ifx\Semafor\ValueOn
\item
I assume sum over index $s$ in expression like
\[
\pC s0a\ x\ \pC s1a
\]
\fi

\ePrints{0701.238,0812.4763,0906.0135,0908.3307,0909.0855}
\ifx\Semafor\ValueOn
\item
We can consider division ring $D$ as $D$\Hyph vector space
of dimension $1$. According to this statement, we can explore not only
homomorphisms of division ring $D_1$ into division ring $D_2$,
but also linear maps of division rings.
This means that map is multiplicative over
maximum possible field. In particular, linear map
of division ring $D$ is multiplicative over center $Z(D)$. This statement
does not contradict with
definition of linear map of field because for field $F$ is true
$Z(F)=F$.
When field $F$ is different from
maximum possible, I explicit tell about this in text.
\fi

\ePrints{0912.4061,Report.10}
\ifx\Semafor\ValueOn
\item
For given field $F$, unless otherwise stated,
we consider finite dimensional $F$\Hyph algebra.
%without zero divisors.
\fi

\ePrints{0701.238,0812.4763,0906.0135,0908.3307}
\ifx\Semafor\ValueOn
\item
In spite of noncommutativity of product a lot of statements
remain to be true if we substitute, for instance, right representation by
left representation or right vector space by left
vector space.
To keep this symmetry in statements of theorems
I use symmetric notation.
For instance, I consider \Ds vector space
and \sD vector space.
We can read notation \Ds vector space
as either D\Hyph star\Hyph vector space or
left vector space.
We can read notation \Ds linear dependent vectors
as either D\Hyph star\Hyph linear dependent vectors or
vectors that are linearly dependent from left.
\fi

\ePrints{0701.238,0812.4763,0906.0135,0908.3307,0909.0855,0912.4061,1001.4852}
\ifx\Semafor\ValueOn
\item We consider algebra $A$ which is finite
dimensional vector space over center.
Considering expansion of element of algebra $A$ relative basis $\Basis e$
we use the same root letter to denote this element and its coordinates.
However we do not use vector notation in algebra.
In expression $a^2$, it is not clear whether this is component
of expansion of element
$a$ relative basis, or this is operation $a^2=aa$.
To make text more clear we use separate color for index of element
of algebra. For instance,
\ShowEq{Expansion relative basis in algebra}

\item When we consider finite dimensional algebra we identify
the vector of basis $\Vector e_{\gi 0}$ with unit of algebra.
\fi

\ePrints{0906.0135,0912.3315}
\ifx\Semafor\ValueOn
\item
Since the number of $\mathfrak{F}$\Hyph algebras
in the tower of representations is varying,
then we use vector notation for a tower of
representations. We denote the set
$(A_1,...,A_n)$ of $\mathfrak{F}_i$\Hyph algebras $A_i$, $i=1$, ..., $n$
as $\Vector A$. We denote the set of representations
$(f_{1,2},...,f_{n-1,n})$ of these algebras as $\Vector f$.
Since different algebras have different type, we also
talk about the set of $\Vector{\mathfrak{F}}$\Hyph algebras.
In relation to the set $\Vector A$,
we also use matrix notations 
that we discussed
in section \xRef{0701.238}{section: Concept of Generalized Index}.
For instance, we use the symbol $\Vector A_{[1]}$ to denote the
set of $\Vector{\mathfrak{F}}$\Hyph algebras $(A_2,...,A_n)$.
In the corresponding notation $(\Vector A_{[1]},\Vector f)$ of tower
of representation, we assume that $\Vector f=(f_{2,3},...,f_{n-1,n})$.

\item Since we use vector notation for elements of the
tower of representations, we need convention about notation of operation.
We assume that we get result of operation componentwise. For instance,
\ShowEq{vector notation in tower of representations}
\fi

\ePrints{0912.3315}
\ifx\Semafor\ValueOn
\item
I believe that diagrams of maps are an important tool.
However, sometimes I want
to see the diagram as three dimensional figure
and I expect that this would increase its expressive
power. Who knows what surprises the future holds.
In 1992, at a conference in Kazan, I have described to my colleagues
what advantages the computer preparation of papers has.
8 years later I learned from the letter from Kazan that now we can
prepare paper using LaTeX.
\fi

\item
Without a doubt, the reader of my articles may have questions,
comments, objections. I will appreciate any response.

\end{enumerate}

%auto-ignore
%auto-ignore

\DefEq
{
\[
f:A\rightarrow A
\]
}
{set of mappings A}

\DefEq
{
\begin{equation}
f\circ g=f(g)
\label{eq: set of mappings A, product}
\end{equation}
}
{set of mappings A, product}

\DefEq
{
\[
f\circ g=g\circ f
\]
}
{f circ g = g circ f}

\DefEq
{
\[
\xymatrix
{
A\ar[rr]^f\ar[d]^g&&A\ar[d]^g
\\
A\ar[rr]^f&&A
}
\]
}
{f circ g = g circ f, diagram}

\DefEq
{
\begin{equation}
f_a(x)=a
\label{eq: mapping f a=a}
\end{equation}
}
{mapping f a=a}

\DefEq
{
\begin{equation}
f\circ a=f(a)
\label{eq: set of mappings A, product over scalar}
\end{equation}
}
{set of mappings A, product over scalar}

\DefEq
{
\[
\begin{matrix}
\Vector a_1=
\begin{pmatrix}
a_1^1\\...\\a_1^n
\end{pmatrix}
&...&
\Vector a_p=
\begin{pmatrix}
a_1^p\\...\\a_p^n
\end{pmatrix}
\end{matrix}
\]
}
{Cartesian power of algebra, 1}

\DefEq
{
\[
\omega(\Vector a_1, ..., \Vector a_p)=
\begin{pmatrix}
\omega(a_1^1, ..., a_p^1)
\\...\\
\omega(a_1^n, ..., a_p^n)
\end{pmatrix}
\]
}
{Cartesian power of algebra, 2}

\DefEq
{
$\Vector a\in A^n$
\[
f\RCcirc\Vector a=
\begin{pmatrix}
f^1_1&...&f^1_n
\\...\\
f^n_1&...&f^n_n
\end{pmatrix}
\RCcirc
\begin{pmatrix}
a^1\\...\\a^n
\end{pmatrix}
=
\begin{pmatrix}
f^1_i\circ a^i\\...\\f^n_i\circ a^i
\end{pmatrix}
\]
}
{matrix of endomorphism, 1}

\DefEq
{
\begin{align*}
f\RCcirc\omega(\Vector a_1,...,\Vector a_p)
&=
\begin{pmatrix}
f^1_1&...&f^1_n
\\...\\
f^n_1&...&f^n_n
\end{pmatrix}
\RCcirc
\begin{pmatrix}
\omega(a_1^1, ..., a_p^1)
\\...\\
\omega(a_1^n, ..., a_p^n)
\end{pmatrix}
\\
&=
\begin{pmatrix}
f^1_i\circ \omega(a_1^i, ..., a_p^i)
\\...\\
f^n_i\circ \omega(a_1^i, ..., a_p^i)
\end{pmatrix}
\\
&=
\begin{pmatrix}
\omega(f^1_i\circ a_1^i, ..., f^1_i\circ a_p^i)
\\...\\
\omega(f^n_i\circ a_1^i, ..., f^n_i\circ a_p^i)
\end{pmatrix}
\\
&=
\omega(f\RCcirc\Vector a_1, ..., f\RCcirc\Vector a_p)
\end{align*}
}
{matrix of endomorphism, 2}

\DefEq
{
\begin{align*}
(f\RCcirc g)\RCcirc\omega(\Vector a_1,...,\Vector a_p)
&=
f\RCcirc (g\RCcirc\omega(\Vector a_1,...,\Vector a_p))
=
f\RCcirc \omega(g\RCcirc\Vector a_1,...,g\RCcirc\Vector a_p)
\\
&=
\omega(f\RCcirc (g\RCcirc\Vector a_1),...,f\RCcirc (g\RCcirc\Vector a_p))
\\
&=
\omega((f\RCcirc g)\RCcirc\Vector a_1,...,(f\RCcirc g)\RCcirc\Vector a_p)
\end{align*}
}
{matrix of endomorphism, 3}

\DefEq
{
\begin{align*}
(f\RCcirc g)\RCcirc h
&=
\left(
\left(
f\RCcirc g
\right)^i_j\circ h^j_k
\right)
=
\left(
(f^i_m\circ g^m_j)\circ h^j_k
\right)
\\
&=
\left(
f^i_m\circ (g^m_j\circ h^j_k)
\right)
=
\left(
f^i_m\circ \left(g\RCcirc h\right)^m_k
\right)
\\
&=
f\RCcirc (g\RCcirc h)
\end{align*}
}
{product of matrices of endomorphisms is associative}

\DefEq
{
\symb{b\CRcirc c}0{cr product of functional matrices}
\begin{equation}
\label{eq: cr product of functional matrices}
\left\{ \begin{array}{rcl}
\ShowSymbol{cr product of functional matrices}&
=&(b_b^c\circ c_c^a)\\
(\ShowSymbol{cr product of functional matrices})^a_b&
=&b_b^c\circ c_c^a
\end{array}\right.
\end{equation}
}
{cr product of functional matrices}

\DefEq
{
\symb{b\RCcirc c}0{rc product of functional matrices}
\begin{equation}
\label{eq: rc product of functional matrices}
\left\{ \begin{array}{rcl}
\ShowSymbol{rc product of functional matrices}&
=&(b^a_c\circ c^c_b)\\
(\ShowSymbol{rc product of functional matrices})^a_b&
=&b^a_c\circ c^c_b
\end{array}\right.
\end{equation}
}
{rc product of functional matrices}

\DefEq
{
\begin{equation}
\Vector e_{\gi i}\Vector e_{\gi j}=
B^{\gi k}_{\gi{ij}}\Vector e_{\gi k}
\label{eq: product in associative algebra}
\end{equation}
}
{product in associative algebra}

\DefEq
{
\[
(\Vector e_{\gi i}\Vector e_{\gi j})\Vector e_{\gi k}
=
\Vector e_{\gi i}(\Vector e_{\gi j}\Vector e_{\gi k})
\]
}
{associative product in algebra}

\DefEq
{
\begin{equation}
B^{\gi p}_{\gi{ij}}B^{\gi q}_{\gi{pk}}
=
B^{\gi q}_{\gi{ip}}B^{\gi p}_{\gi{jk}}
\label{eq: associative product in algebra, 1}
\end{equation}
}
{associative product in algebra, 1}

\DefEq
{
\[
(a\otimes b)\circ x=axb
\]
}
{linear mapping of algebra, a b}

\DefEq
{
\begin{equation}
(a\otimes b)\circ(c\otimes d)=(ac)\otimes(db)
\label{eq: linear mapping of algebra, product}
\end{equation}
}
{linear mapping of algebra, product}

\DefEq
{
\begin{equation}
(a^{\gi{ij}}\Vector e_{\gi i}\otimes\Vector e_{\gi j})
\circ
(b^{\gi{kl}}\Vector e_{\gi k}\otimes\Vector e_{\gi l})
=
a^{\gi{ij}}b^{\gi{kl}}B^{\gi p}_{\gi{ik}}B^{\gi q}_{\gi{lj}}
\Vector e_{\gi p}\otimes\Vector e_{\gi q}
\label{eq: linear mapping of algebra, standard representation, product}
\end{equation}
}
{linear mapping of algebra, standard representation, product}

\DefEq
{
\begin{align}
(a^{\gi{ij}}\Vector e_{\gi i}\otimes\Vector e_{\gi j})\circ
((b^{\gi{kl}}\Vector e_{\gi k}\otimes\Vector e_{\gi l})\circ x)
&=
(a^{\gi{ij}}\Vector e_{\gi i}\otimes\Vector e_{\gi j})\circ
(b^{\gi{kl}}\Vector e_{\gi k}x\Vector e_{\gi l})
\nonumber
\\
&=
a^{\gi{ij}}\Vector e_{\gi i}b^{\gi{kl}}\Vector e_{\gi k}x\Vector e_{\gi l}
\Vector e_{\gi j}
\label{eq: linear mapping of algebra, standard representation, product 1}
\\
&=
a^{\gi{ij}}b^{\gi{kl}}B^{\gi p}_{\gi{ik}}B^{\gi q}_{\gi{lj}}
\Vector e_{\gi p}x\Vector e_{\gi q}
\nonumber
\end{align}
}
{linear mapping of algebra, standard representation, product 1}

\DefEq
{
\begin{align*}
a&=a^{\gi{ij}}\Vector e_{\gi i}\otimes\Vector e_{\gi j}
\\
b&=b^{\gi{ij}}\Vector e_{\gi i}\otimes\Vector e_{\gi j}
\\
c&=c^{\gi{ij}}\Vector e_{\gi i}\otimes\Vector e_{\gi j}
\end{align*}
}
{linear mapping of algebra, product associative, 1}

\DefEq
{
\begin{equation}
\begin{array}{r@{}l}
(a\circ b)\circ c
&=
(a^{\gi{ij}}b^{\gi{kl}}B^{\gi p}_{\gi{ik}}B^{\gi q}_{\gi{lj}}
\Vector e_{\gi p}\otimes\Vector e_{\gi q})
\circ
(c^{\gi{ab}}\Vector e_{\gi a}\otimes\Vector e_{\gi b})
\\
\VirtVar&=
a^{\gi{ij}}b^{\gi{kl}}B^{\gi p}_{\gi{ik}}B^{\gi q}_{\gi{lj}}
c^{\gi{ab}}B^{\gi s}_{\gi{pa}}B^{\gi t}_{\gi{bq}}
\Vector e_{\gi s}\otimes\Vector e_{\gi t}
\end{array}
\label{eq: linear mapping of algebra, product associative, 2 1}
\end{equation}
\begin{equation}
\begin{array}{r@{}l}
a\circ( b\circ c)
&=
(a^{\gi{ij}}\Vector e_{\gi i}\otimes\Vector e_{\gi j})
\circ
(b^{\gi{kl}}c^{\gi{ab}}B^{\gi p}_{\gi{ka}}B^{\gi q}_{\gi{bl}}
\Vector e_{\gi p}\otimes\Vector e_{\gi q})
\\
\VirtVar&=
a^{\gi{ij}}b^{\gi{kl}}c^{\gi{ab}}B^{\gi p}_{\gi{ka}}B^{\gi q}_{\gi{bl}}
B^{\gi s}_{\gi{ip}}B^{\gi t}_{\gi{qj}}
\Vector e_{\gi s}\otimes\Vector e_{\gi t}
\end{array}
\label{eq: linear mapping of algebra, product associative, 2 2}
\end{equation}
}
{linear mapping of algebra, product associative, 2}

\DefEq
{
\begin{equation}
a\RCcirc a^{-1\RCcirc} = \delta
\label{eq: rc inverce matrix of endomorphisms, definition}
\end{equation}
}
{rc inverce matrix of endomorphisms, definition}

\DefEq
{
\begin{align}
\label{eq: rc inverse minor, matrix of endomorphisms}
\left(
a^{-1\RCcirc}{}_\cdot{}_-^J{}^\cdot{}_I^-
\right)^{-1\RCcirc}
&=a^J_I
-a^J_{[I]}\RCcirc
\left(
a_\cdot{}_-^{[I]}{}^\cdot{}_{[J]}^-
\right)^{-1\RCcirc}\RCcirc
a^{[J]}_I
\end{align}
}
{rc inverse minor, matrix of endomorphisms}

\DefEq
{
\begin{align}
\label{eq: rc inverse minor, matrix of endomorphisms, 1, 1}
a^{[J]}_{[I]}\RCcirc
a^{-1\RCcirc}{}^{[I]}_J
+ a^{[J]}_I\RCcirc
a^{-1\RCcirc}{}^I_J&=0\\
\label{eq: rc inverse minor, matrix of endomorphisms, 1, 2}
a^J_{[I]}\RCcirc
a^{-1\RCcirc}{}^{[I]}_J
+a^J_I\RCcirc a^{-1\RCcirc}{}^I_J&=\delta
\end{align}
}
{rc inverse minor, matrix of endomorphisms, 1}

\DefEq
{
$\left(a_\cdot{}_-^{[I]}{}^\cdot{}_{[J]}^-\right)^{-1\RCcirc}$
\begin{align}
\label{eq: rc inverse minor, matrix of endomorphisms,3}
a^{-1\RCcirc}{}^{[I]}_J
+\left(a_\cdot{}_-^{[I]}{}^\cdot{}_{[J]}^-\right)^{-1\RCcirc}\RCcirc
a^{[J]}_I\RCcirc
a^{-1\RCcirc}{}^I_J&=0
\end{align}
}
{rc inverse minor, matrix of endomorphisms, 3}

\DefEq
{
\begin{align}
\label{eq: rc inverse minor, matrix of endomorphisms,4}
-a^J_{[I]}\RCcirc
\left(a_\cdot{}_-^{[I]}{}^\cdot{}_{[J]}^-\right)^{-1\RCcirc}\RCcirc
a^{[J]}_I\RCcirc
a^{-1\RCcirc}{}^I_J
+a^J_I\RCcirc a^{-1\RCcirc}{}^I_J&=\delta
\end{align}
}
{rc inverse minor, matrix of endomorphisms, 4}

\DefEq
{
$(A_\cdot{}_-^i{}^\cdot{}_j^-)^{-1}$
}
{index of inverse element, 1}

\DefEq
{
\[
(A_\cdot{}_-^i{}^\cdot{}_j^-)^{-1}=
\frac 1{A_i^j}
\]
}
{index of inverse element, 2}

\DefEq
{
\begin{align}
\label{eq: rc inverse matrix of endomorphisms}
\left(\mathcal{H}a^{-1\RCcirc}\right)_i^j
&=a^j_i
-a^j_{[i]}\RCcirc
\left(a_\cdot{}_-^{[i]}{}^\cdot{}_{[j]}^-\right)^{-1\RCcirc}\RCcirc
a^{[j]}_i
\end{align}
}
{rc inverse matrix of endomorphisms}

\DefEq
{
\symb{\det\left(a,\RCcirc\right)^a_b}
0{a b RCcirc-quasideterminant definition}
\begin{equation}
\label{eq: a b RCcirc-quasideterminant definition}
\ShowSymbol{a b RCcirc-quasideterminant definition}
=\left(\mathcal{H}a^{-1\RCcirc}\right)^a_b
\end{equation}
}
{a b RCcirc-quasideterminant definition}

\DefEq
{
\begin{align}
\label{eq: RCcirc quasideterminant, expression, 1}
\det\left(a,\RCcirc\right)_a^b&=a_a^b
-a^b_{[a]}\RCcirc
\left(a_\cdot{}_-^{[a]}{}^\cdot{}_{[b]}^-\right)^{-1\RCcirc}\RCcirc
a^{[b]}_a\\
\label{eq: RCcirc quasideterminant, expression, 2}
\det\left(a,\RCcirc\right)^b_a&=a^b_a
-a^b_{[a]}\RCcirc
\mathcal{H}\det\left(a_\cdot{}_-^{[a]}{}^\cdot{}_{[b]}^-,\RCcirc\right)
\RCcirc
a^{[b]}_a
\end{align}
}
{RCcirc quasideterminant, expression}

\DefEq
{
$\left(a_\cdot{}_-^{[a]}{}^\cdot{}_{[b]}^-\right)^{-1\RCcirc}$
}
{minor of RCcirc quasideterminant}

\DefEq
{
\[
\det\left(a,\RCcirc\right)_1^1=a_1^1
\]
}
{quasideterminant of matrix of linear mappings, 1}

\DefEq
{
\symb{\det\left(a,\RCcirc\right)}1{RCcirc-quasideterminant definition}
}
{RCcirc-quasideterminant definition}

\DefEq
{
\begin{equation}
\label{eq: quasideterminant and inverse, RCcirc}
a^{-1\RCcirc}=
\mathcal{H}\det\left(a,\RCcirc\right)
\end{equation}
}
{quasideterminant and inverse, RCcirc}

\DefEq
{
\begin{align}
B^{\gi p}_{\gi{ik}}B^{\gi s}_{\gi{pa}}
&=
B^{\gi p}_{\gi{ka}}B^{\gi s}_{\gi{ip}}
\label{eq: linear mapping of algebra, product associative, 3 1}
\\
B^{\gi q}_{\gi{lj}}B^{\gi t}_{\gi{bq}}
&=
B^{\gi q}_{\gi{bl}}B^{\gi t}_{\gi{qj}}
\label{eq: linear mapping of algebra, product associative, 3 2}
\end{align}
}
{linear mapping of algebra, product associative, 3}

\DefEq
{
\begin{equation}
(a\otimes b)\circ((c\otimes d)\circ x)
=(a\otimes b)\circ(cxd)=a(cxd)b=(ac)x(db)
\label{eq: linear mapping of algebra, product 1}
\end{equation}
}
{linear mapping of algebra, product 1}

\DefEq
{
\begin{equation}
(a^{\gi{ij}}\Vector e_{\gi i}\otimes\Vector e_{\gi j})\circ x=
a^{\gi{ij}}\Vector e_{\gi i} x\Vector e_{\gi j}
\label{eq: linear mapping of algebra, standard representation}
\end{equation}
}
{linear mapping of algebra, standard representation}

\DefEq
{
\begin{equation}
f=
\begin{pmatrix}
f^1_1{}^{\cdot\gi{ij}}\Vector e_{\gi i}\otimes \Vector e_{\gi j}
&...&
f^1_n{}^{\cdot\gi{ij}}\Vector e_{\gi i}\otimes \Vector e_{\gi j}
\\...&...&...\\
f^m_1{}^{\cdot\gi{ij}}\Vector e_{\gi i}\otimes \Vector e_{\gi j}
&...&
f^m_n{}^{\cdot\gi{ij}}\Vector e_{\gi i}\otimes \Vector e_{\gi j}
\end{pmatrix}
\label{eq: linear matrix, standard representation}
\end{equation}
}
{linear matrix, standard representation}

\DefEq
{
\begin{equation}
\begin{pmatrix}
f^1_1{}^{\cdot\gi{ij}}\Vector e_{\gi i}\otimes \Vector e_{\gi j}
&...&
f^1_n{}^{\cdot\gi{ij}}\Vector e_{\gi i}\otimes \Vector e_{\gi j}
\\...&...&...\\
f^m_1{}^{\cdot\gi{ij}}\Vector e_{\gi i}\otimes \Vector e_{\gi j}
&...&
f^m_n{}^{\cdot\gi{ij}}\Vector e_{\gi i}\otimes \Vector e_{\gi j}
\end{pmatrix}
\RCcirc
\begin{pmatrix}
a^1\\...\\a^n
\end{pmatrix}
=
\begin{pmatrix}
f^1_k{}^{\cdot\gi{ij}}\Vector e_{\gi i}a^k \Vector e_{\gi j}
\\...\\
f^m_k{}^{\cdot\gi{ij}}\Vector e_{\gi i}a^k \Vector e_{\gi j}
\end{pmatrix}
\label{eq: linear matrix, standard representation, RCcirc product}
\end{equation}
}
{linear matrix, standard representation, RCcirc product}

\DefEq
{
\begin{equation}
\left\{
\begin{matrix}
a^1_k{}^{\cdot\gi{ij}}\Vector e_{\gi i}x^k \Vector e_{\gi j}=b^1
\\...\\
a^m_k{}^{\cdot\gi{ij}}\Vector e_{\gi i}x^k \Vector e_{\gi j}=b^m
\end{matrix}
\right.
\label{eq: system of linear equations in associative algebra}
\end{equation}
}
{system of linear equations in associative algebra}

\DefEq
{
\begin{align}
\label{eq: nonsingular system of linear equations in algebra, solution, matrix}
x&=a^{-1\RCcirc}\RCcirc b\\
\label{eq: nonsingular system of linear equations in algebra, solution, quasideterminant}
x&=\mathcal H\det\left(a,\RCcirc\right)\RCcirc b
\end{align}
}
{nonsingular system of linear equations in algebra}

\DefEq
{
\begin{align}
&\begin{pmatrix}
a^1_1{}^{\cdot\gi{ij}}\Vector e_{\gi i}\otimes \Vector e_{\gi j}
&...&
a^1_n{}^{\cdot\gi{ij}}\Vector e_{\gi i}\otimes \Vector e_{\gi j}
\\...&...&...\\
a^m_1{}^{\cdot\gi{ij}}\Vector e_{\gi i}\otimes \Vector e_{\gi j}
&...&
a^m_n{}^{\cdot\gi{ij}}\Vector e_{\gi i}\otimes \Vector e_{\gi j}
\end{pmatrix}
\RCcirc
\begin{pmatrix}
x^1\\...\\x^n
\end{pmatrix}
=
\begin{pmatrix}
b^1
\\...\\
b^m
\end{pmatrix}
\label{eq: system of linear equations in associative algebra, matrix, 1}
\\
&a\RCcirc x
=b
\label{eq: system of linear equations in associative algebra, matrix, 2}
\end{align}
}
{system of linear equations in associative algebra, matrix}

\DefEq
{
\begin{equation}
\left\{
\begin{matrix}
f^1_1{}^{\cdot\gi{ij}}\Vector e_{\gi i}a^1 \Vector e_{\gi j}
+
f^1_2{}^{\cdot\gi{ij}}\Vector e_{\gi i}a^2 \Vector e_{\gi j}
=b^1
\\
\VirtVar
f^2_1{}^{\cdot\gi{ij}}\Vector e_{\gi i}a^1 \Vector e_{\gi j}
+
f^2_2{}^{\cdot\gi{ij}}\Vector e_{\gi i}a^2 \Vector e_{\gi j}
=b^2
\end{matrix}
\right.
\label{eq: system of linear equations in associative algebra, 2}
\end{equation}
}
{system of linear equations in associative algebra, 2}

\DefEq
{
\begin{equation}
\begin{pmatrix}
f^1_1{}^{\cdot\gi{ij}}\Vector e_{\gi i}\otimes \Vector e_{\gi j}
&
f^1_2{}^{\cdot\gi{ij}}\Vector e_{\gi i}\otimes \Vector e_{\gi j}
\\
\VirtVar
f^2_1{}^{\cdot\gi{ij}}\Vector e_{\gi i}\otimes \Vector e_{\gi j}
&
f^2_2{}^{\cdot\gi{ij}}\Vector e_{\gi i}\otimes \Vector e_{\gi j}
\end{pmatrix}
\RCcirc
\begin{pmatrix}
a^1\\a^2
\end{pmatrix}
=
\begin{pmatrix}
b^1
\\
b^2
\end{pmatrix}
\label{eq: system of linear equations in associative algebra, matrix, 2}
\end{equation}
}
{system of linear equations in associative algebra, matrix, 2}

\DefEq
{
\[
a^k=a^{k\gi p}\Vector e_{\gi p}
\]
}
{ak expansion basis e}

\DefEq
{
\begin{equation}
f^{l\gi p}_{k\gi q}a^{k\gi q}=b^{l\gi p}
\label{eq: system of linear equations in associative algebra, basis}
\end{equation}
}
{system of linear equations in associative algebra, basis}

\DefEq
{
\[
\begin{matrix}
b=(b^i_j{}^{\cdot\gi{kl}}\Vector e_{\gi k}\otimes \Vector e_{\gi l})
&
c=(c^i_j{}^{\cdot\gi{kl}}\Vector e_{\gi k}\otimes \Vector e_{\gi l})
\end{matrix}
\]
}
{linear matrices B C}

\DefEq
{
\begin{equation}
b\RCcirc c=
(
b^i_p{}^{\cdot\gi{ab}}
c^p_j{}^{\cdot\gi{cd}}
B^{\gi k}_{\gi{ac}}
B^{\gi l}_{\gi{bd}}
\Vector e_{\gi k}\otimes \Vector e_{\gi l})
\label{eq: linear matrices B rc C}
\end{equation}
}
{linear matrices B rc C}

\DefEq
{
\begin{align}
(b\RCcirc c)^i_j
&=
(b^i_p{}^{\cdot\gi{ab}}\Vector e_{\gi a}\otimes \Vector e_{\gi b})
\circ
(c^p_j{}^{\cdot\gi{cd}}\Vector e_{\gi c}\otimes \Vector e_{\gi d})
\nonumber\\
&=
b^i_p{}^{\cdot\gi{ab}}c^p_j{}^{\cdot\gi{cd}}
(\Vector e_{\gi a}\Vector e_{\gi c})
\otimes
(\Vector e_{\gi b}\Vector e_{\gi d})
\label{eq: linear matrices B rc C, 1}
\\
&=
b^i_p{}^{\cdot\gi{ab}}c^p_j{}^{\cdot\gi{cd}}
(B^{\gi k}_{\gi{ac}}\Vector e_{\gi k})
\otimes
(B^{\gi l}_{\gi{bd}}\Vector e_{\gi l})
\nonumber
\end{align}
}
{linear matrices B rc C, 1}

\DefEq
{
\begin{equation}
e^i_j{}^{\cdot\gi{km}}=\delta^i_j
\delta^{\gi k}_{\gi 0}
\delta^{\gi m}_{\gi 0}
\label{eq: standard identity matrix}
\end{equation}
}
{standard identity matrix}

\DefEq
{
\begin{equation}
a^i_p{}^{\cdot\gi{ab}}
b^p_j{}^{\cdot\gi{cd}}
B^{\gi k}_{\gi{ac}}
B^{\gi m}_{\gi{bd}}
=\delta^i_j
\delta^{\gi k}_{\gi 0}
\delta^{\gi m}_{\gi 0}
\label{eq: inverce linear matrix}
\end{equation}
}
{inverce linear matrix}

\DefEq
{
\[
f^{l\gi p}_{k\gi q}=
f^1_k{}^{\cdot\gi{ij}}
B^{\gi a}_{\gi{iq}}B^{\gi p}_{\gi{aj}}
\]
}
{system of linear equations in associative algebra, basis, 1}

\section{Product of Mappings}

On the set of mappings
\ShowEq{set of mappings A}
we define product according to rule
\ShowEq{set of mappings A, product}
The equation
\ShowEq{f circ g = g circ f}
is true iff the diagram
\ShowEq{f circ g = g circ f, diagram}
is commutative.

For $a\in A$, there exists mapping
\ShowEq{mapping f a=a}
If we denote mapping $f_a$ by letter $a$,
then using equation
\eqref{eq: set of mappings A, product},
assume
\ShowEq{set of mappings A, product over scalar}

\section{Biring of Functional Matrices}

If $A$ is $\mathfrak{H}$\Hyph algebra
(\citeBib{Burris Sankappanavar,Cohn: Universal Algebra}),
where the operation of addition is defined,
then we consider the set of
\AddIndex{functional matrices}{functional matrix},
whose elements
are mappings
\ShowEq{set of mappings A}
According to definition
\xRef{0701.238}{definition: rc-product of matrices},
we define
\AddIndex{$\RCcirc$\Hyph product of functional matrices}
{rc product of functional matrices}
\ShowEq{rc product of functional matrices}
According to definition
\xRef{0701.238}{definition: cr-product of matrices},
we define
\AddIndex{$\CRcirc$\Hyph product of functional matrices}
{cr product of functional matrices}
\ShowEq{cr product of functional matrices}

\ifx\texFuture\Defined
\begin{definition}
Пусть $A$ - $\mathfrak{H}$\Hyph алгебра.
Тогда на множестве $A^n$ можно определить структуру
$\mathfrak{H}$\Hyph алгебры таким образом, что все операции
выполняются покомпонентно. А именно,
если $\omega$ - $p$\Hyph арная операция,
то для заданных кортежей
\ShowEq{Cartesian power of algebra, 1}
действие операции $\omega$ определено согласно правилу
\ShowEq{Cartesian power of algebra, 2}
$\mathfrak{H}$\Hyph алгебра $A^n$ называется
\AddIndex{Cartesian power $n$ of $\mathfrak{H}$\Hyph algebra}
{Cartesian power of algebra} $A$.
\qed
\end{definition}

\begin{definition}
Пусть $A$ - $\mathfrak{H}$\Hyph алгебра,
в которой определена операция сложения.
Функциональная матрица $a$ называется
\AddIndex{matrix of endomorphisms of $\mathfrak{H}$\Hyph algebra}
{matrix of endomorphisms of H algebra}
$A$, если $a_i^j$ is endomorphism $\mathfrak{H}$\Hyph алгебры $A$.
\qed
\end{definition}

\begin{theorem}
Пусть $A$ - $\mathfrak{H}$\Hyph алгебра,
в которой определена операция сложения.
$n\times n$ matrix of endomorphisms of $\mathfrak{H}$\Hyph algebra $A$
is endomorphism of $\mathfrak{H}$\Hyph algebra $A^n$.
\end{theorem}
\begin{proof}
Let $f$ be the matrix of endomorphisms of $\mathfrak{H}$\Hyph algebra $A$.
Тогда для
\ShowEq{matrix of endomorphism, 1}
Если $\omega$ - $p$\Hyph арная операция, то
\ShowEq{matrix of endomorphism, 2}
\end{proof}

\begin{theorem}
$\RCcirc$\Hyph произведение матриц of endomorphisms of $\mathfrak{H}$\Hyph algebra $A$
является матрицей of endomorphisms of $\mathfrak{H}$\Hyph algebra $A$.
\end{theorem}
\begin{proof}
Пусть $f$, $g$ - матрицы of endomorphisms $\mathfrak{H}$\Hyph алгебры $A$.
Тогда для операции $\omega$, following equation is true
\ShowEq{matrix of endomorphism, 3}
\end{proof}

\begin{theorem}
If the product of endomorphisms of $\mathfrak{H}$\Hyph algebra $A$
is associative,
then the product of matrices of endomorphisms of $\mathfrak{H}$\Hyph algebra $A$
is associative.
\end{theorem}
\begin{proof}
Утверждение теоремы следует from the chain of equations
\ShowEq{product of matrices of endomorphisms is associative}
\end{proof}

\section{Quasideterminant of Matrix of Endomorphisms}

\begin{theorem}
\label{theorem: rc inverse minor, matrix of endomorphisms}
Suppose $n\times n$ matrix $a$
of endomorphisms of $\mathfrak{H}$\Hyph algebra $A$ has
$\RCcirc$\Hyph inverse matrix\footnote{This statement and its proof are based
on statement 1.2.1 from \citeBib{math.QA-0208146}
(\href{http://arxiv.org/PS_cache/math/pdf/0208/0208146.pdf\#Page=8}{page 8})
for matrix over free division ring.}
\ShowEq{rc inverce matrix of endomorphisms, definition}
Then $k\times k$ minor of $\RCcirc$\Hyph inverse matrix satisfy
следующему равенству, при условии, что рассматриваемые
обратные матрицы существуют,\footnote{
The notation
\ShowEq{index of inverse element, 1}
means that we exchange rows and columns
in Hadamard inverse.
\label{footnote: index of inverse element}
We can
formally write this expression in following form
\ShowEq{index of inverse element, 2}
}
\ShowEq{rc inverse minor, matrix of endomorphisms}
\end{theorem}
\begin{proof}
Definition \eqref{eq: rc inverce matrix of endomorphisms, definition}
of $\RCcirc$\Hyph inverse matrix
leads to system of linear equations
\ShowEq{rc inverse minor, matrix of endomorphisms, 1}
We multiply \eqref{eq: rc inverse minor, matrix of endomorphisms, 1, 1} by
\ShowEq{rc inverse minor, matrix of endomorphisms, 3}
Now we can substitute
\eqref{eq: rc inverse minor, matrix of endomorphisms,3} into
\eqref{eq: rc inverse minor, matrix of endomorphisms, 1, 2}
\ShowEq{rc inverse minor, matrix of endomorphisms, 4}
\eqref{eq: rc inverse minor, matrix of endomorphisms}
follows from \eqref{eq: rc inverse minor, matrix of endomorphisms,4}.
\end{proof}

\begin{corollary}
\label{corollary: rc inverse matrix of endomorphisms}
Suppose $n\times n$ matrix $a$
of endomorphisms of $\mathfrak{H}$\Hyph algebra $A$ has
$\RCcirc$\Hyph inverse matrix.
Then elements of $\RCcirc$\Hyph inverse matrix satisfy
to the equation\footref{footnote: index of inverse element}
\ShowEq{rc inverse matrix of endomorphisms}
\qed
\end{corollary}

\begin{definition}
\label{definition: rc-quasideterminant, matrix of endomorphisms}
\AddIndex{$(^a_b)$\hyph $\RCcirc$\Hyph quasideterminant}
{a b RCcirc-quasideterminant}
of $n\times n$ matrix $a$
is formal expression\footref{footnote: index of inverse element}
\ShowEq{a b RCcirc-quasideterminant definition}
According to the remark
\xRef{0701.238}{remark: combine the notation of indexes}
we can get $(^b_a)$\hyph \RC quasideterminant
as an element of the matrix
\ShowEq{RCcirc-quasideterminant definition}
which we call
\AddIndex{$\RCcirc$\Hyph quasideterminant}
{RCcirc-quasideterminant definition}.
\qed
\end{definition}

\begin{theorem}
Expression for elements of $\RCcirc$\Hyph inverse matrix has form
\ShowEq{quasideterminant and inverse, RCcirc}
\end{theorem}
\begin{proof}
\eqref{eq: quasideterminant and inverse, RCcirc} follows
from \eqref{eq: a b RCcirc-quasideterminant definition}.
\end{proof}

\begin{theorem}
Expression for
$(^a_b)$\hyph $\RCcirc$\Hyph quasideterminant
can be evaluated by either form
\ShowEq{RCcirc quasideterminant, expression}
\end{theorem}
\begin{proof}
Statement follows
from \eqref{eq: rc inverse matrix of endomorphisms} and
\eqref{eq: a b RCcirc-quasideterminant definition}.
\end{proof}
\fi

\section{Algebra of Linear Mappings}

Let $A$ be associative algebra over field $F$.
Let $\Basis e$ be basis of algebra $A$ over field $F$.
Let $\Vector e_{\gi 0}$ be unit of algebra $A$.
The product in the algebra $A$ is defined according to rule
\ShowEq{product in associative algebra}

The equation
\ShowEq{associative product in algebra, 1}
follows from equation
\ShowEq{associative product in algebra}

For given $a$, $b\in A$ we define linear mapping
$a\otimes b$ according to rule
\ShowEq{linear mapping of algebra, a b}
The sum of linear mappings is also a linear mapping.
The set of linear mappings is algebra $A\otimes A$.
We define the product in algebra $A\otimes A$ according to
rule
\ShowEq{linear mapping of algebra, product}
The equation \eqref{eq: linear mapping of algebra, product}
follows from the equation
\ShowEq{linear mapping of algebra, product 1}

Hereinafter we will use the standard representation
\ShowEq{linear mapping of algebra, standard representation}

\begin{theorem}
The product of linear mappings
in the standard representation
has form
\ShowEq{linear mapping of algebra, standard representation, product}
\end{theorem}
\begin{proof}
The statement of the theorem follows from the equation
\ShowEq{linear mapping of algebra, standard representation, product 1}
\end{proof}

\begin{theorem}
\label{theorem: product of linear mappings is associative}
The product of linear mappings of algebra $A$
is associative.
\end{theorem}
\begin{proof}
Consider linear mapping
\ShowEq{linear mapping of algebra, product associative, 1}
According to equation
\eqref{eq: linear mapping of algebra, standard representation, product}
\ShowEq{linear mapping of algebra, product associative, 2}
According to equation
\eqref{eq: associative product in algebra, 1}
\ShowEq{linear mapping of algebra, product associative, 3}
The statement of theorem follows from equations
\eqref{eq: linear mapping of algebra, product associative, 2 1},
\eqref{eq: linear mapping of algebra, product associative, 2 2},
\eqref{eq: linear mapping of algebra, product associative, 3 1},
\eqref{eq: linear mapping of algebra, product associative, 3 2}.
\end{proof}

\section{Matrix of Linear Mappings}

Let $A$ be associative algebra over field $F$.

\begin{definition}
Functional matrix $a$ is called
\AddIndex{matrix of linear mappings}
{matrix of linear mappings},
if $a_i^j$ is linear mapping of algebra $A$.
\qed
\end{definition}

\begin{theorem}
The product of matrices of linear mappings of algebra $A$
is associative.
\end{theorem}
\begin{proof}
The statement of theorem follows from theorem
\ref{theorem: product of linear mappings is associative}
and chain of equations
\ShowEq{product of matrices of endomorphisms is associative}
\end{proof}

\ifx\texFuture\Defined
\begin{theorem}
$n\times n$ matrix of linear mappings of algebra $A$
is matrix of linear mapping
of vector space over algebra $A$.
\end{theorem}
\begin{proof}
Let $f$ be the matrix of endomorphisms of $\mathfrak{H}$\Hyph algebra $A$.
Than for
\ShowEq{matrix of endomorphism, 1}
If $\omega$ is $p$\Hyph арная операция, то
\ShowEq{matrix of endomorphism, 2}
\end{proof}
\fi

\begin{theorem}
$\RCcirc$\Hyph product of matrices of linear mappings is
a matrix of linear mappings.
\end{theorem}
\begin{proof}
The statement is corollary of equation
\eqref{eq: rc product of functional matrices}
and theorem
\xRef{0701.238}{theorem: product of additive map, D D D}.
\end{proof}

We will use
the standard representation for notation of linear mapping.
The standard representation of matrix of linear mappings has form
\ShowEq{linear matrix, standard representation}
According to the definition
\eqref{eq: set of mappings A, product over scalar},
\ShowEq{linear matrix, standard representation, RCcirc product}

\begin{theorem}
\label{theorem: linear matrices B rc C}
For given matrices of linear mappings
\ShowEq{linear matrices B C}
$\RCcirc$\Hyph product has form
\ShowEq{linear matrices B rc C}
\end{theorem}
\begin{proof}
From the equation
\eqref{eq: cr product of functional matrices},
it follows that
\ShowEq{linear matrices B rc C, 1}
Equation
\eqref{eq: linear matrices B rc C}
follows from equation
\eqref{eq: linear matrices B rc C, 1}.
\end{proof}

\section{Quasideterminant of Matrix of Linear Mappings}

\begin{theorem}
\label{theorem: rc inverse minor, matrix of linear mappings}
Suppose $n\times n$ matrix $a$
of linear mappings of algebra $A$ has
$\RCcirc$\Hyph inverse matrix\footnote{This statement and its proof are based
on statement 1.2.1 from \citeBib{math.QA-0208146}
(\href{http://arxiv.org/PS_cache/math/pdf/0208/0208146.pdf\#Page=8}{page 8})
for matrix over free division ring.}
\ShowEq{rc inverce matrix of endomorphisms, definition}
Then $k\times k$ minor of $\RCcirc$\Hyph inverse matrix satisfy
to the following equation, on conditions that considered
inverse matrices exist,\footnote{
The notation
\ShowEq{index of inverse element, 1}
means that we exchange rows and columns
in Hadamard inverse.
\label{footnote: index of inverse element}
We can
formally write this expression in following form
\ShowEq{index of inverse element, 2}
}
\ShowEq{rc inverse minor, matrix of endomorphisms}
\end{theorem}
\begin{proof}
Definition \eqref{eq: rc inverce matrix of endomorphisms, definition}
of $\RCcirc$\Hyph inverse matrix
leads to system of linear equations
\ShowEq{rc inverse minor, matrix of endomorphisms, 1}
We multiply \eqref{eq: rc inverse minor, matrix of endomorphisms, 1, 1} by
\ShowEq{rc inverse minor, matrix of endomorphisms, 3}
Now we can substitute
\eqref{eq: rc inverse minor, matrix of endomorphisms,3} into
\eqref{eq: rc inverse minor, matrix of endomorphisms, 1, 2}
\ShowEq{rc inverse minor, matrix of endomorphisms, 4}
\eqref{eq: rc inverse minor, matrix of endomorphisms}
follows from \eqref{eq: rc inverse minor, matrix of endomorphisms,4}.
\end{proof}

\begin{corollary}
\label{corollary: rc inverse matrix of linear mappings}
Suppose $n\times n$ matrix $a$
of linear mappings of algebra $A$ has
$\RCcirc$\Hyph inverse matrix.
Then elements of $\RCcirc$\Hyph inverse matrix satisfy
to the equation\footref{footnote: index of inverse element}
\ShowEq{rc inverse matrix of endomorphisms}
\qed
\end{corollary}

\begin{definition}
\label{definition: rc-quasideterminant, matrix of endomorphisms}
\AddIndex{$(^a_b)$\hyph $\RCcirc$\Hyph quasideterminant}
{a b RCcirc-quasideterminant}
of $n\times n$ matrix $a$
is formal expression\footref{footnote: index of inverse element}
\ShowEq{a b RCcirc-quasideterminant definition}
According to the remark
\xRef{0701.238}{remark: combine the notation of indexes}
we can get $(^b_a)$\hyph \RC quasideterminant
as an element of the matrix
\ShowEq{RCcirc-quasideterminant definition}
which we call
\AddIndex{$\RCcirc$\Hyph quasideterminant}
{RCcirc-quasideterminant definition}.
\qed
\end{definition}

\begin{theorem}
\label{theorem: quasideterminant and inverse, RCcirc}
Expression for elements of $\RCcirc$\Hyph inverse matrix has form
\ShowEq{quasideterminant and inverse, RCcirc}
\end{theorem}
\begin{proof}
\eqref{eq: quasideterminant and inverse, RCcirc} follows
from \eqref{eq: a b RCcirc-quasideterminant definition}.
\end{proof}

\begin{theorem}
Expression for
$(^a_b)$\hyph $\RCcirc$\Hyph quasideterminant
can be evaluated by either form
\ShowEq{RCcirc quasideterminant, expression}
\end{theorem}
\begin{proof}
Statement follows
from \eqref{eq: rc inverse matrix of endomorphisms} and
\eqref{eq: a b RCcirc-quasideterminant definition}.
\end{proof}

\begin{theorem}
\label{theorem: quasideterminant of matrix of linear mappings}
Let $a$ be matrix of linear mappings.
Then matrices
$\det\left(a,\RCcirc\right)$ and
$a^{-1\RCcirc}$ are
matrices of linear mappings.
\end{theorem}
\begin{proof}
We will prove the theorem by induction over order of matrix.

For $n=1$, from the equation
\eqref{eq: RCcirc quasideterminant, expression, 1}
it follows that
\ShowEq{quasideterminant of matrix of linear mappings, 1}
Therefore, quasideterminant is a matrix of
linear mappings.
From theorem
\xRef{0912.4061}{theorem: linear equation in associative algebra, root},
it follows that the matrix $a^{-1\RCcirc}$ is
a matrix of linear mappings.

Let the statement of the theorem be true for $n-1$.
Let $a$ be $n\times n$ matrix.
According to assumption of induction,
the matrix
\ShowEq{minor of RCcirc quasideterminant}
in the equation
\eqref{eq: RCcirc quasideterminant, expression, 1}
is a matrix of linear mappings.
Therefore,
$(^a_b)$\hyph $\RCcirc$\Hyph quasideterminant
is linear mapping.
From theorems
\xRef{0912.4061}{theorem: linear equation in associative algebra, root},
\ref{theorem: quasideterminant and inverse, RCcirc},
it follows that the matrix $a^{-1\RCcirc}$ is
a matrix of linear mappings.
\end{proof}

\section{System of Linear Equations in Associative Algebra}

\begin{theorem}
\label{theorem: standard identity matrix}
Identity matrix $e$ has standard representation
\ShowEq{standard identity matrix}
\end{theorem}
\begin{proof}
First of all, elements of identity matrix are different from $0$ only
on diagonal. Therefore, $i=j$.
Since every element on diagonal equal to identity of the field, then
element of the matrix has form
$\Vector e_{\gi 0}\otimes\Vector e_{\gi 0}$.
Therefore, $\gi k=\gi 0$, $\gi m=\gi 0$.
\end{proof}

\begin{theorem}
\label{theorem: inverce linear matrix}
Let $a$ be matrix of linear mappings of algebra $A$.
Let $b$ be matrix, $\RCcirc$\Hyph inverce matrix $a$.
Elements of matrices $a$ and $b$ satisfy to the equation
\ShowEq{inverce linear matrix}
\end{theorem}
\begin{proof}
Equation
\eqref{eq: inverce linear matrix}
follows from equations
\eqref{eq: linear matrices B rc C},
\eqref{eq: standard identity matrix}.
\end{proof}

\begin{definition}
If $n\times n$ matrix $a$
of linear mappings of algebra $A$ has
$\RCcirc$\Hyph inverse matrix
we call matrix $a$
\AddIndex{$\RCcirc$\Hyph nonsingular matrix of linear mappings}
{RCcirc nonsingular matrix of linear mappings}.
Otherwise, we call such matrix
\AddIndex{$\RCcirc$\Hyph singular matrix of linear mappings}
{RCcirc singular matrix of linear mappings}.
\qed
\end{definition}

The system of linear equations in associative algebra
has form
\ShowEq{system of linear equations in associative algebra}
We can write the system of linear equations
\eqref{eq: system of linear equations in associative algebra}
in matrix form
\ShowEq{system of linear equations in associative algebra, matrix}

\begin{definition}
Suppose $a$
is $\RCcirc$\Hyph nonsingular matrix. Appropriate system
of linear equations
\eqref{eq: system of linear equations in associative algebra, matrix, 2}
is called
\AddIndex{$\RCcirc$\Hyph nonsingular system of linear equations}
{RCcirc nonsingular system of linear equations}.
\qed
\end{definition}

\begin{theorem}
\label{theorem: nonsingular system of linear equations in algebra}
Solution of nonsingular system of linear equations
\eqref{eq: system of linear equations in associative algebra, matrix, 2}
is determined uniquely and can be presented
in either form
\ShowEq{nonsingular system of linear equations in algebra}
\end{theorem}
\begin{proof}
Multiplying both sides of equation
\eqref{eq: system of linear equations in associative algebra, matrix, 2}
from left by $a^{-1\RCcirc}$ we get
\eqref{eq: nonsingular system of linear equations in algebra, solution, matrix}.
Using definition
\eqref{eq: a b RCcirc-quasideterminant definition}
we get \eqref{eq: nonsingular system of linear equations in algebra, solution, quasideterminant}.
%Since theorem
%\xEqRef{0612.111}{theorem: two products equal}
%the solution is unique.
\end{proof}

We can also consider solving of
the system of linear equations
\eqref{eq: system of linear equations in associative algebra, matrix, 1}
the same way as is done in the theorem
\xRef{0912.4061}
{theorem: determinant of linear equation in associative algebra over field}
If $a^k$ has expansion
\ShowEq{ak expansion basis e}
then the system of linear equations
\eqref{eq: system of linear equations in associative algebra}
is equivalent to the system of linear equations
\ShowEq{system of linear equations in associative algebra, basis}
where
\ShowEq{system of linear equations in associative algebra, basis, 1}

\ifx\texFuture\Defined
Однако нас интересует решение of system of linear equations, оставаясь
в рамках алгебры $A$.

Consider the system of two linear equations с двумя неизвестными.
Эту систему можно записать в виде
\ShowEq{system of linear equations in associative algebra, 2}
Мы можем записать систему линейных уравнений
\eqref{eq: system of linear equations in associative algebra, 2}
в матричной форме
\ShowEq{system of linear equations in associative algebra, matrix, 2}

показать, что если система имеет единственное решение, то по крайней
мере одна из функций $f^1_1$, $f^2_1$ обратима.
\fi

%auto-ignore
\OpenBiblio

%22. A.G. Kurosh, Lectures on General Algebra, Chelsea, New York, 1963

\BiblioItem{Einstein: Electrodynamics of Moving Bodies}
{
Albert Einstein,
On the Electrodynamics of Moving Bodies, 1905,\\
Hendrik Antoon Lorentz, The Principle of Relativity, 37 - 65,
Translated by W Perrett, G B Jeffery
Courier Dover Publications, 1952
}%

\BiblioItem{Einstein: Foundations of general relativity}
{
Albert Einstein,
Die Grundlage der allgemeinen Relativit\"atstheorie,
Ann. Phys., 1916, {\bf 49}, 769 - 822,\\
Einstein's Annalen Papers: The Complete Collection 1901-1922,
edited by J\"urgen Renn, 517 - 571,\\
Wiley-VCH Verlag GmbH \& Co. KGaA, 2005
}%

\BiblioItem{Einstein: Geometry and Experience}
{
Albert Einstein, Geometry and Experience, (1921)\\
Albert Einstein, Sidelights on Relativity, 25 - 56,\\
Courier Dover Publications, 1983
}%

\BiblioItem{Einstein: Main problems of general relativity}
{
Albert Einstein,
Grundgedanken und Probleme der Relativit\"atstheorie, (1923),\\
Nobelstiftelsen, Les Prix Nobel en 1921 - 1922,
Imprimerie Royale, Stockholm, 1923
}%

\BiblioItem{Einstein: Noneuclidean Geometry and Physics}
{
Albert Einstein,
Nichtenklidische Geometrie in der Physik Neue Rundschan, (1925)
Berlin, S. 16 - 20
}%

\BiblioItem{Einstein: Isaak Newton}
{
Albert Einstein,
Isaak Newton, 1927,
Out of My Later Years, 
Citadel Press, 1995, 219 - 223
}%

\BiblioItem{Einstein: On Science}
{
Albert Einstein,
On Science, 
Cosmic Religion, with Other Opinions and Aphorisms,142 - 146,
New York, 1931, 97 - 103
}%

\BiblioItem{Einstein: Autobiographical Notes}
{
Albert Einstein,
Autobiographical Notes, 1949,\\
Paul A. Schilpp, editor, Albert Einstein: Philosopher-Scientist,
Evanston, 
Illinois, The Library of Living Philosophers, 1949, 1 - 95
}%

\BiblioItem{Cite: 104}
{
Cite 104, Source unknown
}%

\BiblioItem{Ghez}
{
Ghez et al.,
The First Measurement of Spectral Lines in a Short-Period Star Bound to the Galaxy's Central Black Hole: A Paradox of Youth,
\href{http://www.journals.uchicago.edu/ApJ/journal/issues/ApJL/v586n2/16990/brief/16990.abstract.html}{ApJL, 586, L127} (2003),
eprint \href{http://arxiv.org/abs/astro-ph/0302299}{arXiv:astro-ph/0302299} (2003)
}%

\BiblioItem{Schodel}
{
R. Sch\"odel et al.,
A star in a 15.2-year orbit around the supermassive black hole at the centre of the Milky Way,
\href{http://www.nature.com/cgi-taf/DynaPage.taf?file=/nature/journal/v419/n6908/abs/nature01121_fs.html}{Nature 419, 694} (2002)
}%

\BiblioItem{Mielke}
{
Eckehard W. Mielke, Affine generalization of the Komar complex of general relativity,
\href{http://prola.aps.org/searchabstract/PRD/v63/i4/e044018}{Phys. Rev. D 63, 044018} (2001)
}%

\BiblioItem{Obukhov}
{
Yu. N. Obukhov and J. G. Pereira, Metric\hyph affine approach to teleparallel gravity,
\href{http://scitation.aip.org/getabs/servlet/GetabsServlet?prog=normal&id=PRVDAQ000067000004044016000001&idtype=cvips&gifs=Yes}
{Phys. Rev. D 67, 044016} (2003),
eprint \href{http://arxiv.org/abs/gr-qc/0212080}{arXiv:gr-qc/0212080} (2002)
}%

\BiblioItem{Sardanashvily}
{
Giovanni Giachetta, Gennadi Sardanashvily, Dirac Equation in Gauge and Affine-Metric Gravitation Theories,
eprint \href{http://arxiv.org/abs/gr-qc/9511035}{arXiv:gr-qc/9511035} (1995)
}%

\BiblioItem{Gauge}
{
Frank Gronwald and Friedrich W. Hehl, On the Gauge Aspects of Gravity, eprint
\href{http://arxiv.org/abs/gr-qc/9602013}{arXiv:gr-qc/9602013} (1996)
}%

\BiblioItem{Neeman}
{
Yuval Neeman, Friedrich W. Hehl, Test Matter in a Spacetime with Nonmetricity, eprint
\href{http://arxiv.org/abs/gr-qc/9604047}{arXiv:gr-qc/9604047} (1996)
}%

\BiblioItem{torsion}
{
F. W. Hehl, P. von der Heyde, G. D. Kerlick, and J. M. Nester,
General relativity with spin and torsion: Foundations and prospects,\\
\href{http://prola.aps.org/abstract/RMP/v48/i3/p393_1}{Rev. Mod. Phys. 48, 393} (1976)
}%

\BiblioItem{Megged}
{
O. Megged, Post-Riemannian Merger of Yang-Mills Interactions with Gravity,
eprint \href{http://arxiv.org/abs/hep-th/0008135}{arXiv:hep-th/0008135} (2001)
}%

%\BiblioItem{Hehl}
%{
%Friedrich W. Hehl, Uwe Muench,
%eprint \href{http://arxiv.org/abs/gr-qc/9708007}{arXiv:gr-qc/9708007} (1997)
%}%

\BiblioItem{gr-qc-9604027}
{
Yu.N. Obukhov, E.J. Vlachynsky, W. Esser, R. Tresguerres and F.W. Hehl,
An exact solution of the metric\hyph affine gauge theory with dilation, shear, and spin charges,
eprint \href{http://arxiv.org/abs/gr-qc/9604027}{arXiv:gr-qc/9604027} (1996)
}%

\BiblioItem{Weinberg}
{
Steven Weinberg.
The Quantum Theory of Fields. Volume I. Foundations.
Cambridge university press, 1995
}%

\BiblioItem{Reinhardt}
{
Greiner Reinhardt. Field Quantization. Springer.
}%

\BiblioItem{Landau}
{
L. D. Landau, E. M. Lifshich, The classical theory of fields.
Oxford, New York, Pergamon Press
}%

\BiblioItem{Wheeler}
{
Ignazio Ciufolini, John Wheeler. Gravitation and Inertia.
Princeton university press.
}%

\BiblioItem{Anderson02}
{
J. D. Anderson, P. A. Laing, E. L. Lau, A. S. Liu, M. M. Nieto, and S. G. Turyshev,
Study of the anomalous acceleration of Pioneer 10 and 11,
\href{http://prola.aps.org/searchabstract/PRD/v65/i8/e082004}{Phys. Rev. D 65, 082004, 50 pp.}, (2002),
eprint \href{http://arxiv.org/abs/gr-qc/0104064}{arXiv:gr-qc/0104064} (2001)
}%

\BiblioItem{Anderson98}
{
J. D. Anderson, P. A. Laing, E. L. Lau, A. S. Liu, M. M. Nieto, and S. G. Turyshev,
Indication, from Pioneer 10/11, Galileo, and Ulysses Data, of an Apparent Anomalous, Weak, Long-Range Acceleration,
\href{http://prola.aps.org/abstract/PRL/v81/i14/p2858_1}{Phys. Rev. Lett. 81, 2858}, (1998),
eprint \href{http://arxiv.org/abs/gr-qc/9808081}{arXiv:gr-qc/9808081} (1998)
}%

%\BiblioItem{Havas} Peter Havas, The Classical Equations of Motion of Point Particles, I,
%{
%\href{http://prola.aps.org/abstract/PR/v87/i2/p309_1}{Phys. Rev. 87, 309} (1952)
%}%

\BiblioItem
{Serge Lang}
{
Serge Lang,
Algebra, Springer, 2002
}%

\BiblioItem{Burris Sankappanavar}
{
S. Burris, H.P. Sankappanavar,
A Course in Universal Algebra, Springer-Verlag (March, 1982),
\\eprint
\href{http://www.math.uwaterloo.ca/~snburris/htdocs/ualg.html}
{http://www.math.uwaterloo.ca/~snburris/htdocs/ualg.html}
\\(The Millennium Edition)
}%

\BiblioItem{Shilov}
{
G. E. Shilov,
Calculus, Multivariable Functions,
Moscow, Nauka, 1972
}%

\BiblioItem{Kolmogorov Fomin}
{
A. N. Kolmogorov and S. V. Fomin,
Elements of the Theory of Functions and Functional Analysis,
Courier Dover Publication, 1999
}%

\BiblioItem{Lebedev Vorovich}
{
L. P. Lebedev, I. I. Vorovich,
Functional Analysis in Mechanics,
Springer, 2002
}%

\BiblioItem
{Rashevsky}
{
P. K. Rashevsky, Riemann Geometry and Tensor Calculus,\\
Moscow, Nauka, 1967
}%

\BiblioItem
{Kurosh: High Algebra}
{
A. G. Kurosh, High Algebra,
Moscow, Nauka, 1968
}%

\BiblioItem
{Kurosh: General Algebra}
{
A. G. Kurosh, Lectures on General Algebra,
Chelsea Pub Co, 1965 
}%

\BiblioItem
{Sabinin: Smooth Quasigroups}
{
Lev V. Sabinin, Smooth Quasigroups and Loops,
Kluwer Academic Publisher, 1999 
}%

\BiblioItem{Dubrovin Fomenko Novikov part 1}
{
B. A. Dubrovin, A. T. Fomenko, S. P. Novikov,
Modern Geometry - Methods and Applications,\\
Part 1, The Geometry of Surfaces, Transformation Groups, and Fields,\\
Translated by Robert G. Burns,\\
Springer - New York, 1992
}%

\BiblioItem{Korn}
{
Granino A. Korn, Theresa M. Korn,
Mathematical Handbook for Scientists and Engineer,
McGraw-Hill Book Company, New York, San Francisco,
Toronto, London, Sydney, 1968
}%

\BiblioItem{Hocking Young Topology}
{
John G. Hocking, Gail S. Young,
Topology,\\
Courier Dover Publications, 1988
}%

%\BiblioItem{Kleyn} http://www.geocities.com/aleks\_kleyn/Derivative/Derivative.htm

\BiblioItem{Tartaglia}
{
Angelo Tartaglia and Matteo Luca Ruggiero,
Angular Momentum Effects in Michelson\Hyph Morley Type Experiments,
Gen.Rel.Grav. 34, 1371-1382 (2002),\\
eprint \href{http://arxiv.org/abs/gr-qc/0110015}{arXiv:gr-qc/0110015} (2001)
}%

\BiblioItem{Tomozawa}
{
Yukio Tomozawa, Speed of Light in Gravitational Fields, eprint
\href{http://arxiv.org/abs/astro-ph/0303047}{arXiv:astro-ph/0303047} (2004)
}%

\BiblioItem{Magueijo}
{
Joao Magueijo,
Covariant and locally Lorentz-invariant varying speed of light theories,
\href{http://prola.aps.org/abstract/PRD/v62/i10/e103521}{Phys. Rev. D 62, 103521} (2000),
eprint \href{http://arxiv.org/abs/gr-qc/0007036}{arXiv:gr-qc/0007036} (2000)
}%

\BiblioItem{Bassett}
{
Bruce A. Bassett, Stefano Liberati, Carmen Molina-Paris, and Matt Visser,
Geometrodynamics of variable-speed-of-light cosmologies,
\href{http://prola.aps.org/abstract/PRD/v62/i10/e103518}{Phys. Rev. D 62}, 103518 (2000),
eprint \href{http://arxiv.org/abs/astro-ph/0001441}{arXiv:astro-ph/0001441} (2000)
}%

\BiblioItem{C.A. Deavours The Quaternion Calculus}
{
C.A. Deavours, The Quaternion Calculus, 
American Mathematical Monthly, {\bf 80} (1973), pp. 995 - 1008
}%

\BiblioItem{Straumann}
{
Lochlainn O'Raifeartaigh and Norbert Straumann,
Gauge theory: Historical origins and some modern developments,
\href{http://prola.aps.org/abstract/RMP/v72/i1/p1_1}{Rev. Mod. Phys. 72, 1} (2000)
}%

\BiblioItem{Lammerzahl}
{
Claus L\"ammerzahl, Mark P. Haugan,
On the interpretation of Michelson\Hyph Morley experiments,
%\href{http://www.sciencedirect.com/science?\_ob=ArticleURL&\_udi=B6TVM-42WP7CR-1&\_user=10&\_handle=W-WA-A-A-AZ-MsSAYZW-UUW-AUDDYZYZAU-WZCBYCEDW-AZ-U&\_fmt=summary&\_coverDate=04%2F23%2F2001&\_rdoc=1&\_orig=browse&\_srch=%23toc%235538%232001%23997179995%23246657!&\_cdi=5538&view=c&\_acct=C000050221&\_version=1&\_urlVersion=0&\_userid=10&md5=385478cda8c5568dea1aeaf0c43669da}
{Phys. Lett. A282 223-229} (2001),\\
eprint \href{http://arxiv.org/abs/gr-qc/0103052}{arXiv:gr-qc/0103052} (2001)
}%

\BiblioItem{Muller}
{
Holger Muller et al.,
Modern Michelson-Morley Experiment using Cryogenic Optical Resonators,
\href{http://prola.aps.org/searchabstract/PRL/v91/i2/e020401}{Phys. Rev. Lett. 91, 020401} (2003),
eprint \href{http://arxiv.org/abs/physics/0305117}{arXiv:physics/0305117} (2000)
}%

\BiblioItem{Ranada}
{
Antonio F. Ranada,
Pioneer acceleration and variation of light speed: experimental situation,
eprint \href{http://arxiv.org/abs/gr-qc/0402120}{arXiv:gr-qc/0402120} (2004)
}%

\BiblioItem{math.QA-0208146}
{
I. Gelfand, S. Gelfand, V. Retakh, R. Wilson,
Quasideterminants,\\
eprint \href{http://arxiv.org/abs/math.QA/0208146}{arXiv:math.QA/0208146} (2002)
}%

\BiblioItem{q-alg-9705026}
{
I.Gelfand, V.Retakh,
Quasideterminants, I,\\
eprint \href{http://arxiv.org/abs/q-alg/9705026}{arXiv:q-alg/9705026} (1997)
}%

\BiblioItem{Gelfand Retakh 1991}
{
I. Gelfand and V. Retakh, Determinants of Matrices over Noncommutative Rings, Funct.
Anal. Appl. 25 (1991), no. 2, 91-102
}%

\BiblioItem{Gelfand Retakh 1992}
{
I. Gelfand and V. Retakh, A Theory of Noncommutative Determinants and Characteristic
Functions of Graphs, Funct. Anal. Appl. 26 (1992), no. 4, 1-20
}%

\BiblioItem{hep-th-9407124}
{
I. M. Gelfand, D. Krob, A. Lascoux, B. Leclerc, V.S. Retakh and J.-Y. Thibon,
Noncommutative symmetric functions,\\
eprint \href{http://arxiv.org/abs/hep-th/9407124}{arXiv:hep-th/9407124} (1994)
}%

\BiblioItem{Carl Faith 1}
{
Carl Faith, Algebra: Rings, Modules and Categories I,
Springer - Verlag, Berlin - Heidelberg - New York, 1973
}%

%\BiblioItem{Pareigis}
%{
%Bodo Pareigis, Categories and Functors,
%Academic Press - New York - London, 1970
%}%

%\BiblioItem{Beachy}
%{
%John A. Beachy, Introductory Lectures on Rings i Modules,
%Cambridge University Press, 1999
%}%

\BiblioItem{0412.391}
{
Aleks Kleyn,
Basis Manifold,
eprint \href{http://arxiv.org/abs/math.DG/0412391}{arXiv:math.DG/0412391} (2007)
}%

\BiblioItem{0405.027}
{
Aleks Kleyn,
Reference Frame in General Relativity,\\
eprint \href{http://arxiv.org/abs/gr-qc/0405027}{arXiv:gr-qc/0405027} (2008)
}%

\BiblioItem{0405.028}
{
Aleks Kleyn, Metric\hyph Affine Manifold,\\
eprint \href{http://arxiv.org/abs/gr-qc/0405028}{arXiv:gr-qc/0405028} (2008)
}%

\BiblioItem{0612.111}
{
Aleks Kleyn,
Biring of Matrices,\\
eprint \href{http://arxiv.org/abs/math.OA/0612111}{arXiv:math.OA/0612111} (2007)
}%

\BiblioItem{0701.238}
{
Aleks Kleyn,
Lectures on Linear Algebra over Division Ring,\\
eprint \href{http://arxiv.org/abs/math.GM/0701238}{arXiv:math.GM/0701238} (2010)
}%

\BiblioItem{0702.561}
{
Aleks Kleyn,
Fibered $\mathfrak{F}$\Hyph Algebra,\\
eprint \href{http://arxiv.org/abs/math.DG/0702561}{arXiv:math.DG/0702561} (2007)
}%

\BiblioItem{math.RA-0501237}
{
Aleks Kleyn,
Vector Space Over Division Ring,\\
eprint \href{http://arxiv.org/abs/math.RA/0412391}{arXiv:math.RA/0501237} (2007)
}%

\BiblioItem{math.RA-0501237v1}
{
Aleks Kleyn,
Module Over Division Ring, version 1,\\
eprint \href{http://arxiv.org/abs/math/0501237v1}{arXiv:math.RA/0501237v1} (2005)
}%

\BiblioItem{0707.2246}
{
Aleks Kleyn,
Fibered Correspondence,\\
eprint \href{http://arxiv.org/abs/0707.2246}{arXiv:0707.2246} (2007)
}%

\BiblioItem{0803.2620}
{
Aleks Kleyn,
Morphism of \Ts Representations,\\
eprint \href{http://arxiv.org/abs/0803.2620}{arXiv:0803.2620} (2008)
}%

\BiblioItem{0803.3276}
{
Aleks Kleyn,
Lorentz Transformation and General Covariance Principle,\\
eprint \href{http://arxiv.org/abs/0803.3276}{arXiv:0803.3276} (2009)
}%

\BiblioItem{0812.4763}
{
Aleks Kleyn,
Introduction into Calculus over Division Ring,\\
eprint \href{http://arxiv.org/abs/0812.4763}{arXiv:0812.4763} (2010)
}%

\BiblioItem{0906.0135}
{
Aleks Kleyn,
Introduction into Geometry over Division Ring,\\
eprint \href{http://arxiv.org/abs/0906.0135}{arXiv:0906.0135} (2009)
}%

\BiblioItem{0912.4061}
{
Aleks Kleyn,
Linear Equation in Finite Dimensional Algebra,\\
eprint \href{http://arxiv.org/abs/0912.4061}{arXiv:0912.4061} (2009)
}%

\BiblioItem{9705.009}%q-alg-9705009
{
John C. Baez,
An Introduction to n-Categories,\\
eprint \href{http://arxiv.org/abs/q-alg/9705009}{arXiv:q-alg/9705009} (1997)
}%

\BiblioItem{Tolstoi about Anna Karenina}
{
Tolstoi about Anna Karenina,
in book A Karenina Companion, by C. J. G. Turner,
published by Wilfrid Laurier University Press (August 1993)
}%

\BiblioItem
{Cohn: Universal Algebra}
{
Paul M. Cohn,
Universal Algebra,
Springer, 1981
}%

\BiblioItem
{Maunder: Algebraic Topology}
{
C. R. F. Maunder,
Algebraic Topology,
Dover Publications, Inc, Mineola, New York, 1996
}%

\BiblioItem{Pommaret: Partial Differential Equations}
{
J.-F. Pommaret,
Partial Differential Equations and Group Theory,
Springer, 1994
}%

\BiblioItem{Bourbaki: Set Theory}
{
N. Bourbaki,
Theory of sets,
Springer, 2004
}%

\BiblioItem
{Bourbaki: General Topology 1}
{
N. Bourbaki,
General Topology, Chapters 1 - 4,
Springer, 1989
}

\BiblioItem{Bourbaki: General Topology: Chapter 5 - 10}
{
N. Bourbaki,
General Topology, Chapters 5 - 10,
Springer, 1989
}

\BiblioItem{Bourbaki: Topological Vector Space}
{
N. Bourbaki,
Topological Vector Spaces, Chapters 1 - 5,
Transl. by H. G. Eggleston $\&$ S. Madan,
Springer, 2003
}

\BiblioItem{Bourbaki: Real Group Lie}
{
N. Bourbaki,
Lie Groups and Lie Algebras, Chapters 7 - 9,
Translator Andrew Pressley,
Springer, 2005
}

\BiblioItem{Pontryagin: Topological Group}
{
L. S. Pontryagin,
Selected Works, Volume Two, Topological Groups,
Gordon and Breach Science Publishers, 1986
}

\BiblioItem
{Eisenhart: Continuous Groups of Transformations}
{
Eisenhart,
Continuous Groups of Transformations,
Dover Publications, New York, 1961
}

\BiblioItem
{Condon Odabasi}
{
Edward Uhler Condon, Halis Odabasi,
Atomic Structure,
CUP Archive, 1980
}

\BiblioItem{Postnikov: Differential Geometry}
{
Postnikov M. M.,
Geometry IV: Differential geometry,
Moscow, Nauka, 1983
}

\BiblioItem{Fihtengolts: Calculus volume 1}
{
Fihtengolts G. M.,
Differential and Integral Calculus Course, volume 1,
Moscow, Nauka, 1969
}

\BiblioItem{Hatcher: Algebraic Topology}
{
Allen Hatcher,
Algebraic Topology,
Cambridge University Press, 2002
}

\BiblioItem{geometry of differential equations}
{
Vinogradov, A. M., Krasil'shchik, I. S., and Lychagin, V. V.,
Introduction to geometry of nonlinear differential equations,
Nauka, Moscow, 1986
}

\BiblioItem{cohomological analysis}
{
A. M. Vinogradov,
Cohomological Analysis of Partial Differential Equations
and Secondary Calculus,
American Mathematical Society, 2001
}

\BiblioItem{0801.1734}
{
Brandon S. DiNunno, Richard A. Matzner,
The Volume Inside a Black Hole,\\
eprint \href{http://arxiv.org/abs/0801.1734v1}{arXiv:0801.1734v1} (2008)
}

\BiblioItem{Izrail M. Gelfand: Quaternion Groups}
{
I. M. Gelfand, M. I. Graev,
Representation of Quaternion Groups over Localy Compact and
Functional Fields,\\
Funct. Anal. Appl. {\bf 2} (1968) 19 - 33;\\
Izrail Moiseevich Gelfand, Semen Grigorevich Gindikin,\\
Izrail M. Gelfand: Collected Papers, volume II, 435 - 449,\\
Springer, 1989
}

\BiblioItem{Bamberg Sternberg}
{
Paul Bamberg, Shlomo Sternberg,
A course in mathematics for students of physics,
Cambridge University Press, 1991
}

\BiblioItem{Conway Smith}
{
John Horton Conway, Derek Alan Smith,
On quaternions and octonions: their geometry, arithmetic, and symmetry,
A K Peters, Natick, Massachussets, 2003
}

\BiblioItem{Sudbery Quaternionic Analysis}
{
A. Sudbery,
Quaternionic Analysis,
Math. Proc. Camb. Phil. Soc. (1979), {\bf 85}, 199 - 225
}

\BiblioItem{0902.4771}
{
Fabrizio Colombo, Graziano Gentili, Irene Sabadini,
A Cauchy kernel for slice regular functions,\\
eprint \href{http://arxiv.org/abs/0902.4771v1}{arXiv:0902.4771v1} (2009)
}

\BiblioItem{Vadim Komkov}
{
Vadim Komkov,
Variational Principles of Continuum Mechanics with Engineering Applications: Critical Points Theory,\\
Springer, 1986
}

\BiblioItem{Alain Connes 1994}
{
Alain Connes,
Noncommutative Geometry,\\
Academic Press, 1994
}

\BiblioItem{Hamilton papers 3}
{
Sir William Rowan Hamilton,
The Mathematical Papers, Vol. III, Algebra,\\
Cambridge at the University Press, 1967
}

\BiblioItem{Hamilton Elements of Quaternions 1}
{
Sir William Rowan Hamilton,
Elements of Quaternions, Volume I,\\
Longmans, Green, and Co., London, New York, and Bombay, 1899
}

\BiblioItem{Cartan geometry in reper}
{
Elie Cartan, Vladislav V. Goldberg, Serge\u{i} Pavlovich Finikov,\\
Riemannian geometry in an orthogonal frame:
from lectures delivered by Elie Cartan at the Sorbonne in 1926-1927,\\
translated by Vladislav V. Goldberg,\\
World Scientific, 2001
}

\BiblioItem{Moore Yaqub}
{
Hal G. Moore, Adil Yaqub,
A first course in linear algebra with applications,
Edition 3, Academic Press, 1998 
}

\CloseBiblio

%auto-ignore
\OpenIndex
\SetIndexSpace%
\Index%1%1%1-drc form
   {$1$-\drc form}%
   {1-drc form, vector spaces}%
\SetIndexSpace%
\Index%2%2%2- ary fibered relation
   {$2$\Hyph ary fibered relation}%
   {2 ary fibered relation}%
\SetIndexSpace%
\Index%596%A%a b $RCcirc quasideterminant
   {$(^a_b)$\hyph $\RCcirc$\Hyph quasideterminant}%
   {a b RCcirc-quasideterminant}%
\Index%3%A%a b-CR quasideterminant
   {$(^a_b)$\hyph \CR quasideterminant}%
   {a b cr-quasideterminant}%
\Index%4%A%a b-RC quasideterminant
   {$(^a_b)$\hyph \RC quasideterminant}%
   {a b RC-quasideterminant}%
\Index%5%A%A- valued function
   {$A$\Hyph valued function}%
   {A valued function}%
\Index%284%A%absolute value on division ring
   {absolute value on division ring}%
   {absolute value on division ring}%
\Index%100%A%active representation
   {active representation}%
   {active representation}%
\Index%569%A%active representation of loop in basis manifold
   {active representation of loop $\mathfrak A(f)$ in basis manifold $\mathcal B(f)$}%
   {active representation in basis manifold}%
\Index%99%A%active starT- representation
   {active \sT representation}%
   {active representation, vector space}%
\Index%326%A%active transformation of basis manifold of representation
   {active transformation of basis manifold of representation}%
   {active transformation of basis, representation}%
\Index%101%A%active transformation on basis manifold
   {active transformation on basis manifold}%
   {active transformation}%
\Index%102%A%active transformation on the set of drc bases
   {active transformation on the set of \drc bases}%
   {active transformation, vector space}%
\Index%98%A%additive map of division ring generated by map G
   {additive map of division ring generated by map $G$}%
   {additive map generated by map, division ring}%
\Index%97%A%additive map of ring
   {additive map of ring}%
   {Additive map of Ring}%
\Index%96%A%additive mapping of D- vector spaces
   {additive mapping of $D$\Hyph vector spaces}%
   {additive map of D vector spaces}%
\Index%110%A%affine basis
   {affine basis}%
   {Affine Basis}%
\Index%141%A%affine transformation group
   {affine transformation group}%
   {drc affine transformation group}%
\Index%142%A%affine transformation group
   {affine transformation group}%
   {affine transformation group}%
\Index%109%A%affine transformation on basis manifold
   {affine transformation on basis manifold}%
   {affine transformation}%
\Index%105%A%alternative representation of matrix
   {alternative representation of matrix}%
   {Alternative representation}%
\Index%275%A%anholonomic coordinate
   {anholonomic coordinate}%
   {anholonomic coordinate}%
\Index%278%A%anholonomic coordinates of connection
   {anholonomic coordinates of connection}%
   {anholonomic coordinates of connection}%
\Index%276%A%anholonomic coordinates of vector
   {anholonomic coordinates of vector}%
   {vector anholonomic coordinates}%
\Index%277%A%anholonomic coordinates on manifold
   {anholonomic coordinates on manifold}%
   {anholonomic coordinates on manifold}%
\Index%293%A%anholonomity object
   {anholonomity object}%
   {anholonomity object}%
\Index%106%A%antihomomorphism of fibered groups
   {antihomomorphism of fibered groups}%
   {antihomomorphism of fibered groups}%
\Index%107%A%antisymmetric 2- ary fibered relation
   {antisymmetric $2$\Hyph ary fibered relation}%
   {antisymmetric 2 ary fibered relation}%
\Index%582%A%Arc basis for vector space
   {$A\RCstar$\Hyph basis for vector space}%
   {Arc basis, vector space}%
\Index%581%A%Arc linearly dependent vectors
   {$A\RCstar$\Hyph linearly dependent vectors}%
   {linearly dependent, A vector space}%
\Index%580%A%Arc linearly independent vectors
   {$A\RCstar$\Hyph linearly independent vectors}%
   {linearly independent, A vector space}%
\Index%108%A%arity of operation
   {arity of operation}%
   {arity of operation}%
\Index%585%A%associative law for Astar linear mappings of vector spaces
   {associative law for $A\star$\Hyph linear mappings of vector spaces}%
   {associative law for Astar linear mappings of vector spaces}%
\Index%576%A%associative law for Astar- vector space
   {associative law for $A\star$\Hyph vector space}%
   {associative law, Astar vector space}%
\Index%165%A%associative law for covariant starT- representation
   {associative law for covariant \sT representation}%
   {associative law for covariant starT representation}%
\Index%166%A%associative law for covariant Tstar- representation
   {associative law for covariant \Ts representation}%
   {associative law for covariant Tstar representation}%
\Index%162%A%associative law for Drc linear maps of vector bundles
   {associative law for \Drc linear maps of vector bundles}%
   {associative law for drc linear maps of vector bundles}%
\Index%161%A%associative law for drc linear maps of vector spaces
   {associative law for \drc linear maps of vector spaces}%
   {associative law for drc linear maps of vector spaces}%
\Index%164%A%associative law for Dstar- vector fields
   {associative law for $\mathcal D\star$\Hyph vector fields}%
   {associative law, Dstar vector fields}%
\Index%163%A%associative law for Dstar- vector space
   {associative law for $D\star$\Hyph vector space}%
   {associative law, Dstar vector space}%
\Index%167%A%associative law for twin representations
   {associative law for twin representations}%
   {associative law for twin representations}%
\Index%168%A%associative law of composition of fibered correspondences
   {associative law of composition of fibered correspondences}%
   {associative law, composition of fibered correspondences}%
\Index%583%A%Astar linear map of vector spaces
   {$A\star$\Hyph linear map of vector spaces}%
   {Astar linear map of vector spaces}%
\Index%575%A%Astar vector space
   {$A\star$\Hyph vector space}%
   {Astar vector space}%
\Index%579%A%Astar-product of vector over scalar
   {$A\star$\hyph product of vector over scalar}%
   {Astar product of vector over scalar, vector space}%
\Index%95%A%auto parallel line
   {auto parallel line}%
   {auto parallel line}%
\Index%94%A%automorphism of representation of F algebra
   {automorphism of representation of $\mathfrak{F}$\Hyph algebra}%
   {automorphism of representation}%
\Index%560%A%automorphism of tower of representations
   {automorphism of tower of representations}%
   {automorphism of tower of representations}%
\Index%282%A%norm of quaternion
   {norm of quaternion}%
   {norm of quaternion}%
\SetIndexSpace%
\Index%114%B%base of fibered correspondence
   {base of fibered correspondence}%
   {base of fibered correspondence}%
\Index%113%B%base of map
   {base of map}%
   {base of map}%
\Index%256%B%basis manifold of affine space
   {basis manifold of affine space}%
   {Basis Manifold, Affine Space}%
\Index%261%B%basis manifold of central affine space
   {basis manifold of central affine space}%
   {Basis Manifold, Central Affine Space, division ring}%
\Index%262%B%basis manifold of central affine space
   {basis manifold of central affine space}%
   {Basis Manifold, Central Affine Space}%
\Index%255%B%basis manifold of drc vector space
   {basis manifold of \drc vector space}%
   {basis manifold of drc vector space}%
\Index%259%B%basis manifold of Euclid space
   {basis manifold of Euclid space}%
   {Basis Manifold, Euclid Space}%
\Index%260%B%basis manifold of Euclid space
   {basis manifold of Euclid space}%
   {Basis Manifold, Euclid Space, division ring}%
\Index%536%B%basis manifold of representation
   {basis manifold of representation}%
   {basis manifold representation F algebra}%
\Index%258%B%basis manifold of vector space
   {basis manifold of vector space}%
   {basis manifold of vector space}%
\Index%329%B%basis of representation
   {basis of representation}%
   {basis of representation}%
\Index%567%B%basis of tower of representations
   {basis of tower of representations}%
   {basis of tower of representations}%
\Index%115%B%basis of vector space
   {basis of vector space}%
   {Basis}%
\Index%116%B%basis vector of starT- representation
   {basis vector of \sT representation}%
   {basis vector of starT representation}%
\Index%117%B%basis vector of Tstar- representation
   {basis vector of \Ts representation}%
   {basis vector of Tstar representation}%
\Index%122%B%biring
   {biring}%
   {biring}%
\Index%419%B%bundle of level 2
   {bundle of level $2$}%
   {bundle of level 2}%
\Index%418%B%bundle of level n
   {bundle of level $n$}%
   {bundle of level n}%
\SetIndexSpace%
\Index%6%C%c row of matrix
   {\subs row of matrix}%
   {c row}%
\Index%19%C%c rows dcr vector space
   {\subs rows \dcr vector space}%
   {subs rows dcr vector space}%
\Index%12%C%c-row of matrix
   {$c$\hyph row of matrix}%
   {c-row}%
\Index%443%C%Cartan connection
   {Cartan connection}%
   {Cartan connection}%
\Index%227%C%Cartan curvature
   {Cartan curvature}%
   {Cartan curvature}%
\Index%401%C%Cartan derivative
   {Cartan derivative}%
   {Cartan derivative}%
\Index%444%C%Cartan symbol
   {Cartan symbol}%
   {Cartan symbol}%
\Index%338%C%Cartan transport
   {Cartan transport}%
   {Cartan transport}%
\Index%149%C%Cartesian power \mathcal{A} of bundle \mathcal{B}
   {Cartesian power $\mathcal{A}$ of bundle $\mathcal{B}$}%
   {Cartesian power A of bundle B}%
\Index%148%C%Cartesian power A of set B
   {Cartesian power $A$ of set $B$}%
   {Cartesian power of set}%
\Index%150%C%Cartesian power n of bundle \mathcal{E}
   {Cartesian power $n$ of bundle $\mathcal{E}$}%
   {Cartesian power n of bundle E}%
\Index%594%C%Cartesian power n of H algebra
   {Cartesian power $n$ of $\mathfrak{H}$\Hyph algebra}%
   {Cartesian power of algebra}%
\Index%175%C%category of drc vector spaces
   {category of \drc vector spaces}%
   {category of drc vector spaces}%
\Index%179%C%category of fibered correspondences over diagonal
   {category of fibered correspondences over diagonal}%
   {category of fibered correspondences over diagonal}%
\Index%178%C%category of reduced fibered correspondences
   {category of reduced fibered correspondences}%
   {category of reduced fibered correspondences}%
\Index%176%C%category of Tstar- representations of \mathfrak{F}- algebra A
   {category of \Ts representations of $\mathfrak{F}$\Hyph algebra $A$}%
   {category of Tstar representations of F algebra}%
\Index%177%C%category of Tstar- representations of \mathfrak{F}- algebra from category A
   {category of \Ts representations of $\mathfrak{F}$\Hyph algebra from category $\mathcal A$}%
   {category of Tstar representations of F algebra from category}%
\Index%351%C%Cauchy sequence in valued division ring
   {Cauchy sequence in valued division ring}%
   {Cauchy sequence, valued division ring}%
\Index%507%C%center of ring D
   {center of ring $D$}%
   {center of ring}%
\Index%508%C%central affine basis
   {central affine basis}%
   {Central Affine Basis, division ring}%
\Index%509%C%central affine basis
   {central affine basis}%
   {Central Affine Basis}%
\Index%124%C%column vector
   {column vector}%
   {column vector}%
\Index%190%C%commutative diagram of correspondences
   {commutative diagram of correspondences}%
   {commutative diagram of correspondences}%
\Index%191%C%compact-open topology
   {compact\hyph open topology}%
   {compact open topology}%
\Index%349%C%complete division ring
   {complete division ring}%
   {complete division ring}%
\Index%348%C%complete system of linear partial differential equations
   {complete system of linear partial differential equations}%
   {Complete System of Linear Partial Differential Equations}%
\Index%128%C%completely integrable system
   {completely integrable system}%
   {completely integrable system}%
\Index%199%C%component of Gateaux derivative of map \Vector f(\Vector x)
   {component of the G\^ateaux derivative of map $\Vector f(\Vector x)$}%
   {component of Gateaux derivative of map, D vector space}%
\Index%198%C%component of Gateaux derivative of map f(x)
   {component of the G\^ateaux derivative of map $f(x)$}%
   {component of Gateaux derivative of map, division ring}%
\Index%197%C%component of Gateaux derivative of second order of map  of division ring
   {component of the G\^ateaux derivative of second order of map  of division ring}%
   {component of Gateaux derivative of Second Order, division ring}%
\Index%196%C%component of Gateaux derivative of second order of map \Vector f(\Vector x)
   {component of the G\^ateaux derivative of second order of map $\Vector f(\Vector x)$}%
   {component of Gateaux derivative of Second Order, D vector space}%
\Index%193%C%component of linear map f of division ring
   {component of linear map $f$ of division ring}%
   {component of linear map, division ring}%
\Index%192%C%component of linear map of D- vector space
   {component of linear map of $D$\Hyph vector space}%
   {component of linear map, D vector space}%
\Index%194%C%component of polyadditive map \Vector A
   {component of polyadditive map $\Vector A$}%
   {component of polyadditive map, D vector space}%
\Index%195%C%component of polylinear map of division ring
   {component of polylinear map of division ring}%
   {component of polylinear map, division ring}%
\Index%394%C%composition of fibered correspondences
   {composition of fibered correspondences}%
   {composition of fibered correspondences}%
\Index%393%C%composition of reduced fibered correspondences
   {composition of reduced fibered correspondences}%
   {composition of reduced fibered correspondences}%
\Index%495%C%condition of reducibility of products
   {condition of reducibility of products}%
   {condition of reducibility of products}%
\Index%538%C%connection coefficients in D affine space
   {connection coefficients in $D$\Hyph affine space}%
   {connection coefficients, D affine space}%
\Index%281%C%continuous correspondence
   {continuous correspondence}%
   {continuous correspondence}%
\Index%280%C%continuous function of division ring
   {continuous function of division ring}%
   {continuous function, division ring}%
\Index%200%C%contravariant starT- representation of fibered group
   {contravariant \sT representation of fibered group}%
   {contravariant starT representation of fibered group}%
\Index%202%C%contravariant starT- representation of group
   {contravariant \sT representation of group}%
   {contravariant starT representation of group}%
\Index%203%C%contravariant Tstar- representation of fibered group
   {contravariant \Ts representation of fibered group}%
   {contravariant Tstar representation of fibered group}%
\Index%201%C%contravariant Tstar- representation of group
   {contravariant \Ts representation of group}%
   {contravariant Tstar representation of group}%
\Index%214%C%coordinate drc isomorphism
   {coordinate \drc isomorphism}%
   {coordinate drc isomorphism}%
\Index%210%C%coordinate Drc vector bundle
   {coordinate \Drc vector bundle}%
   {coordinate drc vector bundle}%
\Index%209%C%coordinate drc vector space
   {coordinate \drc vector space}%
   {coordinate drc vector space}%
\Index%215%C%coordinate isomorphism
   {coordinate isomorphism}%
   {coordinate isomorphism}%
\Index%206%C%coordinate matrix of set of vectors in dcr rows vector space
   {coordinate matrix of set of vectors in \dcr rows vector space}%
   {coordinate matrix of set of vectors, dcr vector space}%
\Index%207%C%coordinate matrix of set of vectors in drc rows vector space
   {coordinate matrix of set of vectors in \drc rows vector space}%
   {coordinate matrix of set of vectors, drc vector space}%
\Index%205%C%coordinate matrix of vector field in Drc basis
   {coordinate matrix of vector field in \Drc basis}%
   {coordinate matrix of vector field in drc basis}%
\Index%204%C%coordinate matrix of vector in drc basis
   {coordinate matrix of vector in \drc basis}%
   {coordinate matrix of vector in drc basis}%
\Index%208%C%coordinate reference frame
   {coordinate reference frame}%
   {coordinate reference frame}%
\Index%212%C%coordinate representation in drc vector space
   {coordinate representation in \drc vector space}%
   {coordinate representation, drc vector space}%
\Index%213%C%coordinate representation of group in vector space
   {coordinate representation of group in vector space}%
   {coordinate representation, vector space}%
\Index%211%C%coordinate vector space
   {coordinate vector space}%
   {coordinate vector space}%
\Index%323%C%coordinates of automorphism of representation
   {coordinates of automorphism of representation}%
   {coordinates of automorphism of representation}%
\Index%325%C%coordinates of basis of representation
   {coordinates of basis of representation}%
   {coordinates of basis relative to basis, representation}%
\Index%558%C%coordinates of element relative to set
   {coordinates of element $m$ relative to set $X$}%
   {coordinates of element relative to set, representation}%
\Index%218%C%coordinates of geometrical object
   {coordinates of geometrical object}%
   {coordinates of geometrical object, vector space}%
\Index%220%C%coordinates of geometrical object in coordinate drc vector space
   {coordinates of geometrical object in coordinate \drc vector space}%
   {coordinates of geometrical object, coordinate drc vector space}%
\Index%221%C%coordinates of geometrical object in coordinate representation
   {coordinates of geometrical object in coordinate representation}%
   {coordinates of geometrical object, coordinate vector space}%
\Index%219%C%coordinates of geometrical object in drc vector space
   {coordinates of geometrical object in \drc vector space}%
   {coordinates of geometrical object, drc vector space}%
\Index%571%C%coordinates of passive map of basis manifold of representation
   {coordinates of passive map of basis manifold of representation}%
   {coordinates of passive map of basis manifold of representation}%
\Index%224%C%coordinates of representation
   {coordinates of representation}%
   {coordinates of representation, drc vector space}%
\Index%225%C%coordinates of representation
   {coordinates of representation}%
   {coordinates of representation}%
\Index%222%C%coordinates of set of vectors in dcr vector space
   {coordinates of set of vectors in \dcr vector space}%
   {coordinates of set of vectors, dcr vector space}%
\Index%223%C%coordinates of set of vectors in drc vector space
   {coordinates of set of vectors in \drc vector space}%
   {coordinates of set of vectors, drc vector space}%
\Index%217%C%coordinates of vector field in Drc basis
   {coordinates of vector field in \Drc basis}%
   {coordinates of vector field in drc basis}%
\Index%216%C%coordinates of vector in drc basis
   {coordinates of vector in \drc basis}%
   {coordinates of vector in drc basis}%
\Index%455%C%correspondence continuous on the set
   {correspondence continuous on the set}%
   {correspondence continuous on the set}%
\Index%454%C%correspondence of homomorphism
   {correspondence of homomorphism}%
   {correspondence of homomorphism}%
\Index%185%C%covariant starT- representation of fibered group
   {covariant \sT representation of fibered group}%
   {covariant starT representation of fibered group}%
\Index%184%C%covariant starT- representation of group
   {covariant \sT representation of group}%
   {covariant starT representation of group}%
\Index%187%C%covariant Tstar- representation of fibered group
   {covariant \Ts representation of fibered group}%
   {covariant Tstar representation of fibered group}%
\Index%186%C%covariant Tstar- representation of group
   {covariant \Ts representation of group}%
   {covariant Tstar representation of group}%
\Index%8%C%CR inverse element of biring
   {\CR inverse element of biring}%
   {cr-inverse element}%
\Index%7%C%CR matrix group
   {\CR matrix group}%
   {cr-matrix group}%
\Index%10%C%CR power
   {\CR power}%
   {cr power}%
\Index%589%C%cr product of functional matrices
   {$\CRcirc$\Hyph product of functional matrices}%
   {cr product of functional matrices}%
\Index%9%C%CR product of matrices
   {\CR product of matrices}%
   {cr-product of matrices}%
\Index%11%C%crd vector space
   {\crd vector space}%
   {crd vector space}%
\Index%537%C%curvilinear coordinates of point in affine space
   {curvilinear coordinates of point in affine space}%
   {curvilinear coordinates of point in affine space}%
\SetIndexSpace%
\Index%23%D%d affine space
   {$D$\Hyph affine space}%
   {d affine space}%
\Index%549%D%D- affine connection on manifold with affine connections
   {$D$\Hyph affine connection on manifold with affine connections}%
   {D affine connection, affine manifold}%
\Index%48%D%D- valued variable
   {$D$\Hyph valued variable}%
   {D valued variable}%
\Index%14%D%D- vector function
   {$D$\Hyph vector function}%
   {d vector function}%
\Index%550%D%D-affine connection coefficients on manifold
   {$D$\Hyph affine connection coefficients on manifold}%
   {D affine connection coefficients, manifold}%
\Index%15%D%D-vector space
   {$D$\hyph vector space}%
   {D vector space}%
\Index%16%D%dcr basis of c rows vector space
   {\dcr basis of \subs rows vector space}%
   {dcr basis, c rows vector space}%
\Index%17%D%dcr vector
   {\dcr vector}%
   {dcr vector}%
\Index%18%D%dcr vector space
   {\dcr vector space}%
   {dcr vector space}%
\Index%303%D%determinant of matrix
   {determinant of matrix}%
   {determinant}%
\Index%309%D%deviation of trajectories
   {deviation of trajectories}%
   {deviation of trajectories}%
\Index%153%D%diagonal in bundle
   {diagonal in bundle}%
   {diagonal in bundle}%
\Index%154%D%diagram of correspondences
   {diagram of correspondences}%
   {diagram of correspondences}%
\Index%417%D%dimension of drc vector space
   {dimension of \drc vector space}%
   {dimension of vector space}%
\Index%151%D%direct product of bundles
   {direct product of bundles}%
   {Cartesian product of bundles}%
\Index%410%D%direct product of D- vector spaces
   {direct product of $D$\Hyph vector spaces}%
   {direct product of D vector spaces}%
\Index%415%D%direct product of division rings
   {direct product of division rings}%
   {direct product of division rings}%
\Index%411%D%direct product of drc vector spaces
   {direct product of \drc vector spaces}%
   {direct product, drc vector space}%
\Index%414%D%direct product of representations of fibered group
   {direct product of representations of fibered group}%
   {direct product of representations of fibered group}%
\Index%413%D%direct product of representations of group
   {direct product of representations of group}%
   {direct product of representations of group}%
\Index%152%D%direct product of total spaces
   {direct product of total spaces}%
   {Cartesian product of total spaces}%
\Index%412%D%direct product of Tstar- representations of group
   {direct product of \Ts representations of group}%
   {direct product of Tstar representations of group}%
\Index%409%D%direct sum of representations
   {direct sum of representations}%
   {direct sum of representations}%
\Index%270%D%direction in D- vector space \Vector V over field P
   {direction in $D$\Hyph vector space $\Vector V$ over field $P$}%
   {direction over field, D vector space}%
\Index%271%D%direction in ring D over commutative ring P
   {direction in ring $D$ over commutative ring $P$}%
   {direction over commutative ring, ring}%
\Index%577%D%distributive law for Astar- vector space
   {distributive law for $A\star$\Hyph vector space}%
   {distributive law, Astar vector space}%
\Index%170%D%distributive law for Dstar- vector fields
   {distributive law for $\mathcal D\star$\Hyph vector fields}%
   {distributive law, Dstar vector fields}%
\Index%169%D%distributive law for Dstar- vector space
   {distributive law for $D\star$\Hyph vector space}%
   {distributive law, Dstar vector space}%
\Index%111%D%drc affine basis
   {\drc affine basis}%
   {drc affine basis, division ring}%
\Index%22%D%drc automorphism of vector space
   {\drc automorphism of vector space}%
   {automorphism of vector space}%
\Index%25%D%drc basis for r rows vector space
   {\drc basis for \sups rows vector space}%
   {drc basis, r rows vector space}%
\Index%26%D%Drc basis for vector  bundle
   {\Drc basis for vector  bundle}%
   {drc basis, vector bundle}%
\Index%24%D%drc basis for vector space
   {\drc basis for vector space}%
   {drc basis, vector space}%
\Index%257%D%drc basis manifold of affine space
   {\drc basis manifold of affine space}%
   {drc Basis Manifold, Affine Space, division ring}%
\Index%31%D%drc isomorphism of vector spaces
   {\drc isomorphism of vector spaces}%
   {isomorphism of vector spaces}%
\Index%36%D%Drc linear map of vector bundles
   {\Drc linear map of vector bundles}%
   {drc linear map of vector bundles}%
\Index%35%D%drc linear map of vector spaces
   {\drc linear map of vector spaces}%
   {drc linear map of vector spaces}%
\Index%32%D%drc linear span in vector space
   {\drc linear span in vector space}%
   {linear span, vector space}%
\Index%37%D%drc linear starT- representation of group
   {\drc linear \sT representation of group}%
   {drc linear starT representation of group}%
\Index%20%D%Drc linearly dependent vector fields
   {\Drc linearly dependent vector fields}%
   {linearly dependent vector fields}%
\Index%21%D%drc linearly dependent vectors
   {\drc linearly dependent vectors}%
   {linearly dependent, vector space}%
\Index%33%D%Drc linearly independent vector fields
   {\Drc linearly independent vector fields}%
   {linearly independent vector fields}%
\Index%34%D%drc linearly independent vectors
   {\drc linearly independent vectors}%
   {linearly independent, vector space}%
\Index%449%D%drc system of linear equations
   {\drc system of linear equations}%
   {system of linear equations}%
\Index%27%D%drc vector
   {\drc vector}%
   {drc vector}%
\Index%30%D%drc vector function
   {\drc vector function}%
   {drc vector function}%
\Index%28%D%drc vector space
   {\drc vector space}%
   {drc vector space}%
\Index%42%D%Dstar component of coordinates of vector \Vector r
   {\Ds component of coordinates of vector $\Vector r$}%
   {Dstar component of coordinates of vector, D vector space}%
\Index%41%D%Dstar- vector bundle
   {$\mathcal D\star$\Hyph vector bundle}%
   {Dstar vector bundle}%
\Index%39%D%Dstar- vector field
   {$\mathcal D\star$\Hyph vector field}%
   {Dstar vector field}%
\Index%40%D%Dstar- vector space
   {$D\star$\hyph  vector space}%
   {Dstar vector space}%
\Index%43%D%Dstar-linear composition of vector fields
   {$\mathcal D\star$\hyph linear composition of vector fields}%
   {linear composition of vector fields}%
\Index%45%D%Dstar-product of vector field over scalar
   {$\mathcal D\star$\hyph product of vector field over scalar}%
   {Dstar product of vector field over scalar, vector space}%
\Index%44%D%Dstar-product of vector over scalar
   {$D\star$\hyph product of vector over scalar}%
   {Dstar product of vector over scalar, vector space}%
\Index%159%D%dual space of drc vector space
   {dual space of \drc vector space}%
   {dual space of drc vector space}%
\Index%382%D%duality principle for biring
   {duality principle for biring}%
   {duality principle for biring}%
\Index%383%D%duality principle for biring of matrices
   {duality principle for biring of matrices}%
   {duality principle for biring of matrices}%
\SetIndexSpace%
\Index%517%E%effective representation of \mathfrak{F}- algebra A
   {effective representation of $\mathfrak{F}$\Hyph algebra $A$}%
   {effective representation of algebra}%
\Index%521%E%effective representation of division ring
   {effective representation of division ring}%
   {effective representation of division ring}%
\Index%519%E%effective representation of fibered \mathfrak{F}- algebra
   {effective representation of fibered $\mathfrak{F}$\Hyph algebra}%
   {effective representation of fibered F-algebra}%
\Index%518%E%effective representation of group
   {effective representation of group}%
   {effective representation of group}%
\Index%324%E%effective tower of Tstar- representations
   {effective tower of \Ts representations}%
   {effective tower of representations}%
\Index%514%E%effective Tstar- representation of fibered division ring
   {effective \Ts representation of fibered division ring}%
   {effective representation of fibered division ring}%
\Index%520%E%effective Tstar- representation of fibered group
   {effective \Ts representation of fibered group}%
   {effective representation of fibered group}%
\Index%516%E%effective Tstar- representation of group
   {effective \Ts representation of group}%
   {effective Tstar representation of group}%
\Index%557%E%endomorphism of representation of F algebra
   {endomorphism of representation of $\mathfrak{F}$\Hyph algebra}%
   {endomorphism of representation}%
\Index%327%E%endomorphism of representation regular on generating set X
   {endomorphism of representation regular on generating set $X$}%
   {endomorphism of representation, regular on set}%
\Index%328%E%endomorphism of representation singular on generating set X
   {endomorphism of representation singular on generating set $X$}%
   {endomorphism of representation, singular on set}%
\Index%559%E%endomorphism of tower of representations
   {endomorphism of tower of representations}%
   {endomorphism of tower of representations}%
\Index%564%E%endomorphism of tower of representations regular on tuple of generating sets
   {endomorphism of tower of representations regular on tuple of generating sets}%
   {endomorphism of representation, regular on tuple}%
\Index%565%E%endomorphism of tower of representations singular on tuple of generating sets
   {endomorphism of tower of representations singular on tuple of generating sets}%
   {endomorphism of representation, singular on tuple}%
\Index%50%E%enhanced Lie group
   {enhanced Lie group}%
   {enhanced Lie group}%
\Index%51%E%essential parameters
   {essential parameters}%
   {essential parameters}%
\Index%539%E%Euclidean metric on division ring
   {Euclidean metric on division ring}%
   {Euclidean metric on division ring}%
\Index%547%E%Euclidean scalar product in D- vector space
   {Euclidean scalar product in $D$\Hyph vector space}%
   {Euclidean scalar product, vector space}%
\Index%542%E%Euclidean scalar product on division ring
   {Euclidean scalar product on division ring}%
   {Euclidean scalar product on division ring}%
\Index%437%E%extended matrix of drc linear equations
   {extended matrix of \drc linear equations}%
   {extended matrix, system of drc linear equations}%
\Index%438%E%extended matrix of rcd linear equations
   {extended matrix of \rcd linear equations}%
   {extended matrix, system of rcd linear equations}%
\Index%386%E%extension of correspondence
   {extension of correspondence}%
   {extension of correspondence}%
\Index%515%E%extreme line
   {extreme line}%
   {extreme line}%
\SetIndexSpace%
\Index%420%F%fibered \mathfrak{F}- algebra
   {fibered $\mathfrak{F}$\Hyph algebra}%
   {fibered F-algebra}%
\Index%421%F%fibered \mathfrak{F}- subalgebra
   {fibered $\mathfrak{F}$\Hyph subalgebra}%
   {fibered F-subalgebra}%
\Index%434%F%fibered coordinate Drc isomorphism
   {fibered coordinate \Drc isomorphism}%
   {fibered coordinate drc isomorphism}%
\Index%432%F%fibered correspondence from \mathcal{A} to \mathcal{B}
   {fibered correspondence from $\mathcal{A}$ to $\mathcal{B}$}%
   {fibered correspondence from A to B}%
\Index%430%F%fibered correspondence in \mathcal{A}
   {fibered correspondence in $\mathcal{A}$}%
   {fibered correspondence in A}%
\Index%431%F%fibered correspondence of homomorphism
   {fibered correspondence of homomorphism}%
   {fibered correspondence of homomorphism}%
\Index%427%F%fibered equivalence
   {fibered equivalence}%
   {fibered equivalence}%
\Index%422%F%fibered group
   {fibered group}%
   {fibered group}%
\Index%436%F%fibered identification morphism
   {fibered identification morphism}%
   {fibered identification morphism}%
\Index%424%F%fibered little group
   {fibered little group}%
   {fibered little group}%
\Index%435%F%fibered morphism from bundle \mathcal{A} into \mathcal{B}
   {fibered morphism from bundle $\mathcal{A}$ into $\mathcal{B}$}%
   {fibered morphism from A into B}%
\Index%433%F%fibered natural morphism
   {fibered natural morphism}%
   {fibered natural morphism}%
\Index%426%F%fibered ordering
   {fibered ordering}%
   {fibered ordering}%
\Index%425%F%fibered preordering
   {fibered preordering}%
   {fibered preordering}%
\Index%428%F%fibered ring
   {fibered ring}%
   {fibered ring}%
\Index%423%F%fibered stability group
   {fibered stability group}%
   {fibered stability group}%
\Index%429%F%fibered subset
   {fibered subset}%
   {fibered subset}%
\Index%479%F%field-strength tensor
   {field-strength tensor}%
   {field-strength tensor}%
\Index%497%F%filter \mathfrak{F} converges to A
   {filter $\mathfrak{F}$ converges to $A$}%
   {filter converges}%
\Index%337%F%first Newton law
   {first Newton law}%
   {First Newton law}%
\Index%442%F%free Tstar- representation of fibered group
   {free \Ts representation of fibered group}%
   {free representation of fibered group}%
\Index%441%F%free Tstar- representation of group
   {free \Ts representation of group}%
   {free representation of group}%
\Index%339%F%Frenet transport
   {Frenet transport}%
   {Frenet transport}%
\Index%506%F%function continuous with respect to set of arguments
   {function continuous with respect to set of arguments}%
   {function continuous with respect to set of arguments}%
\Index%503%F%function is projective over field R and continuous in direction over field R
   {function is projective over field $R$ and continuous in direction over field $R$}%
   {projective function is continuous in direction over field R, division ring}%
\Index%502%F%function of \gi n D- valued variables
   {function of $\gi n$ $D$\Hyph valued variables}%
   {function of n D valued variables}%
\Index%501%F%function of D- vector space \Vector{V} to D|Hyph vector space \Vector W differentiable in Gateaux sense
   {function of $D$\Hyph vector space $\Vector{V}$ to $D$|Hyph vector space $\Vector W$ differentiable in the G\^ateaux sense}%
   {function differentiable in Gateaux sense, D vector space}%
\Index%505%F%function of division ring differentiable in Gateaux sense
   {function of division ring differentiable in the G\^ateaux sense}%
   {function differentiable in Gateaux sense, division ring}%
\Index%504%F%function of division ring Dstar differentiable in the Fr\'echet sense
   {function of division ring \Ds differentiable in the Fr\'echet sense}%
   {function Dstar differentiable in Frechet sense, division ring}%
\Index%588%F%functional matrix
   {functional matrix}%
   {functional matrix}%
\Index%500%F%fundamental sequence in valued division ring
   {fundamental sequence in valued division ring}%
   {fundamental sequence, valued division ring}%
\SetIndexSpace%
\Index%52%G%G- reference frame
   {$G$\Hyph reference frame}%
   {G reference frame}%
\Index%53%G%G-basis\ of vector space
   {\Gbasis\ of vector space}%
   {G-basis}%
\Index%54%G%G-coords\ of basis
   {\Gcoords\ of basis}%
   {G-coordinates}%
\Index%55%G%G-space
   {\Gspace}%
   {GSpace}%
\Index%13%G%Gateaux crd derivative of map \Vector f of D-vector space \Vector V to D-vector space \Vector W
   {the G\^ateaux \crd derivative of map $\Vector f$ of $D$\hyph vector space $\Vector V$ to $D$\hyph vector space $\Vector W$}%
   {Gateaux crd derivative of map, D vector space}%
\Index%398%G%Gateaux derivative of map \Vector f of normed D- vector space \Vector{V} to normed D- vector space \Vector{W}
   {the G\^ateaux derivative of map $\Vector f$ of normed $D$\Hyph vector space $\Vector{V}$ to normed $D$\Hyph vector space $\Vector{W}$}%
   {Gateaux derivative of map, D vector space}%
\Index%397%G%Gateaux derivative of map f
   {the G\^ateaux derivative of map $f$}%
   {Gateaux derivative of map, division ring}%
\Index%399%G%Gateaux derivative of order n of map \Vector f
   {the G\^ateaux derivative of order $n$ of map $\Vector f$}%
   {Gateaux derivative of Order n, D vector space}%
\Index%400%G%Gateaux derivative of order n of map f of division ring
   {the G\^ateaux derivative of order $n$ of map $f$ of division ring}%
   {Gateaux derivative of Order n, division ring}%
\Index%396%G%Gateaux derivative of second order of map of division ring
   {the G\^ateaux derivative of second order of map of division ring}%
   {Gateaux derivative of Second Order, division ring}%
\Index%157%G%Gateaux differential of map \Vector f of normed D- vector space \Vector{V} to normed D- vector space \Vector{W}
   {the G\^ateaux differential of map $\Vector f$ of normed $D$\Hyph vector space $\Vector{V}$ to normed $D$\Hyph vector space $\Vector{W}$}%
   {Gateaux differential of map, D vector space}%
\Index%158%G%Gateaux differential of map f
   {the G\^ateaux differential of map $f$}%
   {Gateaux differential of map, division ring}%
\Index%155%G%Gateaux differential of second order of map of division ring
   {the G\^ateaux differential of second order of map of division ring}%
   {Gateaux differential of Second Order, division ring}%
\Index%156%G%Gateaux differential of second order of mapping \Vector f
   {the G\^ateaux differential of second order of mapping $\Vector f$}%
   {Gateaux differential of Second Order, D vector space}%
\Index%38%G%Gateaux drc derivative of map \Vector f of D- vector space \Vector V to D- vector space \Vector W
   {the G\^ateaux \drc derivative of map $\Vector f$ of $D$\Hyph vector space $\Vector V$ to $D$\Hyph vector space $\Vector W$}%
   {Gateaux drc derivative of map, D vector space}%
\Index%46%G%Gateaux Dstar derivative of map f of division ring D
   {the G\^ateaux \Ds derivative of map $f$ of division ring $D$}%
   {Gateaux Dstar derivative of map, division ring}%
\Index%453%G%Gateaux mixed partial derivative of map f^j with respect to variables v^i, v^j
   {the G\^ateaux mixed partial derivative of map $f^j$ with respect to variables $v^i$, $v^j$}%
   {Gateaux partial derivative of Second Order, D vector space}%
\Index%512%G%Gateaux partial derivative of map f^j with respect to variable v^i
   {the G\^ateaux partial derivative of map $f^j$ with respect to variable $v^i$}%
   {Gateaux partial derivative, D vector space}%
\Index%510%G%Gateaux partial drc derivative of map f^b with respect to variable x^a
   {the G\^ateaux partial \drc derivative of map $f^b$ with respect to variable $x^a$}%
   {Gateaux partial crd derivative of map with respect to variable, D vector space}%
\Index%511%G%Gateaux partial drc derivative of map f^b with respect to variable x^a
   {the G\^ateaux partial \drc derivative of map $f^b$ with respect to variable $x^a$}%
   {Gateaux partial drc derivative of map with respect to variable, D vector space}%
\Index%76%G%Gateaux starD derivative of map f of division ring D
   {the G\^ateaux \sD derivative of map $f$ of division ring $D$}%
   {Gateaux starD derivative of map, division ring}%
\Index%332%G%generating set of representation
   {generating set of representation}%
   {generating set of representation}%
\Index%331%G%generating set of subrepresentation
   {generating set of subrepresentation}%
   {generating set of subrepresentation}%
\Index%289%G%generator of additive map
   {generator of additive map}%
   {generator of additive map, division ring}%
\Index%136%G%geometrical object defined in drc vector space
   {geometrical object defined in \drc vector space}%
   {geometrical object, drc vector space}%
\Index%133%G%geometrical object in coordinate representation
   {geometrical object in coordinate representation}%
   {geometrical object, coordinate vector space}%
\Index%134%G%geometrical object in coordinate representation defined in drc vector space
   {geometrical object in coordinate representation defined in \drc vector space}%
   {geometrical object, coordinate drc vector space}%
\Index%132%G%geometrical object in vector space
   {geometrical object in vector space}%
   {geometrical object, vector space}%
\Index%135%G%geometrical object of type A
   {geometrical object of type $A$ in vector space}%
   {geometrical object of type A, vector space}%
\Index%146%G%group algebra
   {group algebra}%
   {group algebra}%
\SetIndexSpace%
\Index%292%H%Hadamard inverse of matrix
   {Hadamard inverse of matrix}%
   {Hadamard inverse of matrix}%
\Index%546%H%hermitian conjugated vector
   {hermitian conjugated vector}%
   {hermitian conjugated vector}%
\Index%541%H%hermitian conjugation in division ring
   {hermitian conjugation in division ring}%
   {hermitian conjugation, division ring}%
\Index%544%H%hermitian metric on division ring
   {hermitian metric on division ring}%
   {hermitian metric on division ring}%
\Index%350%H%hermitian scalar product in D- vector space
   {hermitian scalar product in $D$\Hyph vector space}%
   {hermitian scalar product, vector space}%
\Index%545%H%hermitian scalar product on division ring
   {hermitian scalar product on division ring}%
   {hermitian scalar product on division ring}%
\Index%138%H%holonomic coordinates of connection
   {holonomic coordinates of connection}%
   {holonomic coordinates of connection}%
\Index%137%H%holonomic coordinates of vector
   {holonomic coordinates of vector}%
   {vector holonomic coordinates}%
\Index%296%H%homogeneous bundle of fibered group
   {homogeneous bundle of fibered group}%
   {homogeneous bundle of fibered group}%
\Index%294%H%homogeneous map of degree k over field F
   {homogeneous map of degree $k$ over field $F$}%
   {homogeneous map of degree over field, D vector space}%
\Index%295%H%homogeneous space of group
   {homogeneous space of group}%
   {homogeneous space of group}%
\Index%139%H%homomorphism of fibered \mathfrak{F}- algebras
   {homomorphism of fibered $\mathfrak{F}$\Hyph algebras}%
   {homomorphism of fibered F-algebras}%
\Index%140%H%homomorphism of fibered groups
   {homomorphism of fibered groups}%
   {homomorphism of fibered groups}%
\SetIndexSpace%
\Index%121%I%infinitesimal generator
   {infinitesimal generator}%
   {infinitesimal generator}%
\Index%174%I%infinitesimal generators of group Lie
   {infinitesimal generators of group Lie}%
   {infinitesimal generators of group Lie}%
\Index%384%I%invariance principle
   {invariance principle}%
   {invariance principle}%
\Index%385%I%invariance principle in vector space
   {invariance principle in vector space}%
   {invariance principle, vector space}%
\Index%291%I%inverse fibered correspondence
   {inverse fibered correspondence}%
   {inverse fibered correspondence}%
\Index%290%I%inverse reduced fibered correspondence
   {inverse reduced fibered correspondence}%
   {inverse reduced fibered correspondence}%
\Index%173%I%isomorphism of fibered \mathfrak{F}- algebras
   {isomorphism of fibered $\mathfrak{F}$\Hyph algebras}%
   {isomorphism of fibered F-algebras}%
\Index%266%I%isomorphism of repesentations of F algebra
   {isomorphism of repesentations of $\mathfrak{F}$\Hyph algebra}%
   {isomorphism of repesentations of F algebra}%
\SetIndexSpace%
\Index%592%J%Jacobian complete system of differential equations
   {Jacobian complete system of differential equations}%
   {Jacobian complete system of differential equations}%
\Index%593%J%Jacobian complete system of drc differential equations
   {Jacobian complete system of \drv differential equations}%
   {Jacobian complete system of drc differential equations}%
\Index%527%J%the Jacobi Gateaux matrix of map of D- vector space
   {the Jacobi\Hyph G\^ateaux matrix of map of $D$\Hyph vector space}%
   {Jacobi Gateaux matrix of map, D vector space}%
\SetIndexSpace%
\Index%522%K%kernel of additive map of D- vector space
   {kernel of additive map of $D$\Hyph vector space}%
   {kernel of additive map, D vector space}%
\Index%523%K%kernel of additive map of division ring
   {kernel of additive map of division ring}%
   {kernel of additive map, division ring}%
\Index%526%K%kernel of inefficiency of representation of fibered group
   {kernel of inefficiency of representation of fibered group}%
   {kernel of inefficiency of representation of fibered group}%
\Index%525%K%kernel of inefficiency of representation of group
   {kernel of inefficiency of representation of group}%
   {kernel of inefficiency of representation of group}%
\Index%524%K%kernel of inefficiency of Tstar- representation of group G
   {kernel of inefficiency of \Ts representation of group $G$}%
   {kernel of inefficiency of Tstar representation of group}%
\Index%493%K%Killing equation
   {Killing equation}%
   {Killing equation}%
\Index%494%K%Killing equation of second type
   {Killing equation of second type}%
   {Killing equation second type}%
\Index%123%K%Killing vector of second type
   {Killing vector of second type}%
   {Killing vector second type}%
\Index%445%K%Kronecker symbol
   {Kronecker symbol}%
   {Kronecker symbol}%
\SetIndexSpace%
\Index%301%L%left defined Lie algebra of Lie group
   {left defined Lie algebra of Lie group}%
   {left defined Lie algebra}%
\Index%229%L%left invariant vector field
   {left invariant vector field}%
   {left invariant vector}%
\Index%237%L%left module
   {left module}%
   {left module}%
\Index%240%L%left shift on fibered group
   {left shift on fibered group}%
   {Tstar shift, fibered group}%
\Index%238%L%left shift on group
   {left shift on group}%
   {left shift}%
\Index%239%L%left shift on group
   {left shift on group}%
   {left shift, group}%
\Index%236%L%left structural constant of Lie algebra
   {left structural constant of Lie algebra}%
   {left structural constant of Lie algebra}%
\Index%228%L%left vector space
   {left vector space}%
   {left vector space}%
\Index%231%L%left-side contravariant representation of group
   {left-side contravariant representation of group}%
   {left-side contravariant representation of group}%
\Index%230%L%left-side covariant representation of group
   {left-side covariant representation of group}%
   {left-side covariant representation of group}%
\Index%232%L%left-side representation of \mathfrak{F}- algebra A in H algebra M
   {left-side representation of $\mathfrak{F}$\Hyph algebra $A$ in $\mathfrak{H}$\Hyph algebra $M$}%
   {left-side representation of algebra}%
\Index%233%L%left-side representation of fibered \mathfrak{F}- algebra
   {left-side representation of fibered $\mathfrak{F}$\Hyph algebra}%
   {left-side representation of fibered F-algebra}%
\Index%234%L%left-side transformation
   {left-side transformation}%
   {left-side transformation}%
\Index%235%L%left-side transformation on bundle
   {left-side transformation on bundle}%
   {left-side transformation of bundle}%
\Index%104%L%Lie algebra of Lie group
   {Lie algebra of Lie group}%
   {algebra Lie group Lie}%
\Index%402%L%Lie derivative
   {Lie derivative}%
   {Lie derivative}%
\Index%404%L%Lie derivative of connection
   {Lie derivative of connection}%
   {Lie derivative of connection}%
\Index%403%L%Lie derivative of metric
   {Lie derivative of metric}%
   {Lie derivative of metric}%
\Index%118%L%Lie group basic operators
   {Lie group basic operators}%
   {Lie group basic operators}%
\Index%246%L%lift of correspondence
   {lift of correspondence}%
   {lift of correspondence}%
\Index%245%L%lift of map
   {lift of map}%
   {lift of map}%
\Index%365%L%limit of correspondence with respect to the filter
   {limit of correspondence with respect to the filter}%
   {limit of correspondence with respect to the filter}%
\Index%366%L%limit of filter
   {limit of filter}%
   {limit of filter}%
\Index%364%L%limit of sequence in valued division ring
   {limit of sequence in valued division ring}%
   {limit of sequence, valued division ring}%
\Index%367%L%limit set of filter
   {limit set of filter}%
   {limit set of filter}%
\Index%241%L%linear map of D- vector spaces
   {linear map of $D$\Hyph vector spaces}%
   {linear map, vector space}%
\Index%242%L%linear map of D- vector spaces over field F
   {linear map of $D$\Hyph vector spaces over field $F$}%
   {linear map over field, vector space}%
\Index%243%L%linear map of division ring
   {linear map of division ring}%
   {linear map, division ring}%
\Index%244%L%linear representation of group
   {linear representation of group}%
   {linear representation of group}%
\Index%250%L%little group
   {little group}%
   {little group}%
\Index%247%L%local reference frame
   {local reference frame}%
   {local reference frame}%
\Index%248%L%locally compact at point p space
   {locally compact at point $p$ space}%
   {locally compact at point space}%
\Index%249%L%locally compact space
   {locally compact space}%
   {locally compact space}%
\Index%568%L%loop of automorphisms of representation
   {loop of automorphisms of representation}%
   {loop of automorphisms of representation}%
\Index%373%L%Lorentz transformation
   {Lorentz transformation}%
   {Lorentz transformation}%
\SetIndexSpace%
\Index%56%M%m- dimensional parallelepiped
   {$m$\Hyph dimensional parallelepiped}%
   {m dimensional parallelepiped}%
\Index%57%M%m- vector
   {$m$\Hyph vector}%
   {m-vector}%
\Index%556%M%manifold with D- affine connections
   {manifold with $D$\Hyph affine connections}%
   {manifold with D- affine connections}%
\Index%310%M%map of D- vector space, multiplicative over field
   {map of $D$\Hyph vector space, multiplicative over field}%
   {map multiplicative over field, D vector space}%
\Index%311%M%map of D- vector space, projective over field
   {map of $D$\Hyph vector space, projective over field}%
   {map projective over field, D vector space}%
\Index%315%M%map of division ring, projective over commutative ring
   {map of division ring, projective over commutative ring}%
   {map projective over commutative ring, ring}%
\Index%313%M%map of ring D, linear over commutative ring F
   {map of ring $D$, linear over commutative ring $F$}%
   {linear map over commutative ring, ring}%
\Index%314%M%map of ring, multiplicative over commutative ring
   {map of ring, multiplicative over commutative ring}%
   {map multiplicative over commutative ring, ring}%
\Index%312%M%map of rings R_1, ..., R_n polylinear over commutative ring P
   {map polylinear over commutative ring}%
   {map polylinear over commutative ring, ring}%
\Index%316%M%map of type G on manifold
   {map of type $G$ on manifold}%
   {map of type G on manifold}%
\Index%572%M%map polylinear over finite dimensional algebras
   {map polylinear over finite dimensional algebras}%
   {map polylinear over finite dimensional algebras}%
\Index%407%M%mapping space
   {mapping space}%
   {mapping space}%
\Index%584%M%matrix of Astar linear map
   {matrix of $A\star$\Hyph linear map}%
   {matrix of Astar linear map}%
\Index%252%M%matrix of bilinear function
   {matrix of bilinear function}%
   {matrix of bilinear function}%
\Index%251%M%matrix of drc linear map
   {matrix of \drc linear map}%
   {matrix of drc linear map}%
\Index%595%M%matrix of endomorphisms of H algebra
   {matrix of endomorphisms of $\mathfrak{H}$\Hyph algebra}%
   {matrix of endomorphisms of H algebra}%
\Index%253%M%matrix of fibered Drc linear map
   {matrix of fibered \Drc linear map}%
   {matrix of fibered drc linear map}%
\Index%591%M%matrix of linear mappings
   {matrix of linear mappings}%
   {matrix of linear mappings}%
\Index%552%M%matrix of quadratic map
   {matrix of quadratic map}%
   {matrix of quadratic map, division ring}%
\Index%254%M%metric-affine manifold
   {metric-affine manifold}%
   {metric-affine manifold}%
\Index%265%M%morphism from tower of Tstar- representations into tower of Tstar- representations
   {morphism from tower of \Ts representations into tower of \Ts representations}%
   {morphism from tower of representations into tower of representations}%
\Index%269%M%morphism of fibered Tstar- representations from \mathcal{F} into \mathcal{G}
   {morphism of fibered \Ts representations from $\mathcal{F}$ into $\mathcal{G}$}%
   {morphism of fibered representations from f into g}%
\Index%268%M%morphism of representations from f into g
   {morphism of representations from $f$ into $g$}%
   {morphism of representations from f into g}%
\Index%267%M%morphism of representations of \mathfrak{F}- algebra in \mathfrak{H}- algebra
   {morphism of representations of $\mathfrak{F}$\Hyph algebra in $\mathfrak{H}$\Hyph algebra}%
   {morphism of representations of F algebra in H algebra}%
\Index%264%M%morphism of Tstar- representations of fibered \mathfrak{F}- algebra
   {morphism of \Ts representations of fibered $\mathfrak{F}$\Hyph algebra}%
   {morphism of representations of fibered F algebra}%
\Index%147%M%movement on basis manifold
   {movement on basis manifold}%
   {movement transformation}%
\SetIndexSpace%
\Index%58%N%n- ary fibered relation
   {$n$\Hyph ary fibered relation}%
   {fibered relation}%
\Index%279%N%nonmetricity
   {nonmetricity}%
   {nonmetricity}%
\Index%272%N%nonsingular bilinear function
   {nonsingular bilinear function}%
   {nonsingular bilinear function}%
\Index%274%N%nonsingular transformation
   {nonsingular transformation}%
   {nonsingular transformation}%
\Index%286%N%norm of map \Vector A of normed D-vector space
   {norm of map $\Vector A$ of normed $D$\hyph vector space}%
   {norm of map, D vector space}%
\Index%285%N%norm of mapping of division ring
   {norm of mapping of division ring}%
   {norm of map, division ring}%
\Index%283%N%norm on D- vector space
   {norm on $D$\Hyph vector space}%
   {norm on D vector space}%
\Index%287%N%normed D- vector space
   {normed $D$\Hyph vector space}%
   {normed D vector space}%
\SetIndexSpace%
\Index%300%O%operation on bundle
   {operation on bundle}%
   {operation on bundle}%
\Index%408%O%opposite fibered preordering
   {opposite fibered preordering}%
   {opposite fibered preordering}%
\Index%306%O%orbit of representation of fibered group
   {orbit of representation of fibered group}%
   {orbit of representation of fibered group}%
\Index%305%O%orbit of representation of group
   {orbit of representation of group}%
   {orbit of representation of group}%
\Index%304%O%orbit of Tstar- representation of group
   {orbit of \Ts representation of group}%
   {orbit of Tstar  representation of group}%
\Index%307%O%orthonornal basis
   {orthonornal basis}%
   {Orthonornal Basis}%
\Index%308%O%orthonornal basis
   {orthonornal basis}%
   {Orthonornal Basis, division ring}%
\SetIndexSpace%
\Index%317%P%parallelogram
   {parallelogram}%
   {parallelogram}%
\Index%513%P%partial additive map of variable v^i
   {partial additive map of variable $v^i$}%
   {partial additive map of variable}%
\Index%333%P%passive representation
   {passive representation}%
   {passive representation}%
\Index%570%P%passive transformation of the basis manifold of representation
   {passive transformation of the basis manifold of representation}%
   {passive transformation of basis, representation}%
\Index%335%P%passive transformation on basis manifold
   {passive transformation on basis manifold}%
   {passive transformation}%
\Index%336%P%passive transformation on the set of drc bases
   {passive transformation on the set of \drc bases}%
   {passive transformation, vector space}%
\Index%321%P%passive Tstar- representation
   {passive \Ts representation}%
   {passive Tstar representation}%
\Index%416%P%pfaffian derivative
   {pfaffian derivative}%
   {pfaffian derivative}%
\Index%343%P%polyadditive map of rings
   {polyadditive map of rings}%
   {polyadditive map of rings}%
\Index%342%P%polyadditive mapping of (n)-D-vector spaces
   {polyadditive mapping of $(n)$\hyph $D$\hyph vector spaces}%
   {polyadditive map of D vector spaces}%
\Index%346%P%polylinear map of rings
   {polylinear map of rings}%
   {polylinear map of rings}%
\Index%345%P%polylinear skew symmetric map
   {polylinear skew symmetric map}%
   {polylinear map skew symmetric, division ring}%
\Index%347%P%polylinear symmetric map
   {polylinear symmetric map}%
   {polylinear map symmetric, division ring}%
\Index%344%P%polyvector
   {polyvector}%
   {polyvector}%
\Index%352%P%potential energy
   {potential energy}%
   {potential energy}%
\Index%388%P%product of geometrical object and constant
   {product of geometrical object and constant}%
   {product of geometrical object and constant}%
\Index%389%P%product of geometrical object and constant in vector space
   {product of geometrical object and constant in vector space}%
   {product of geometrical object and constant, vector space}%
\Index%390%P%product of groups
   {product of groups}%
   {product of groups}%
\Index%391%P%product of morphisms of representations of \mathfrak{F}- algebra
   {product of morphisms of representations of $\mathfrak{F}$\Hyph algebra}%
   {product of morphisms of representations of F algebra}%
\Index%561%P%product of morphisms of tower of representations
   {product of morphisms of tower of representations}%
   {product of morphisms of tower of representations}%
\Index%392%P%product of morphisms of Tstar- representations of fibered \mathfrak{F}- algebra
   {product of morphisms of \Ts representations of fibered $\mathfrak{F}$\Hyph algebra}%
   {product of morphisms of representations of fibered F algebra}%
\Index%395%P%product of objects in category
   {product of objects in category}%
   {product of objects in category}%
\Index%387%P%projection of bundle \mathcal{E} along fiber E
   {projection of bundle $\mathcal{E}$ along fiber $E$}%
   {projection of bundle along fiber}%
\Index%540%P%pseudo-Euclidean metric on division ring
   {pseudo\Hyph Euclidean metric on division ring}%
   {pseudo-Euclidean metric on division ring}%
\Index%548%P%pseudo-Euclidean scalar product in D- vector space
   {pseudo\Hyph Euclidean scalar product in $D$\Hyph vector space}%
   {pseudo-Euclidean scalar product, vector space}%
\Index%543%P%pseudo-Euclidean scalar product on division ring
   {pseudo-Euclidean scalar product on division ring}%
   {pseudo-Euclidean scalar product on division ring}%
\SetIndexSpace%
\Index%555%Q%quadratic form in division ring
   {quadratic form in division ring}%
   {quadratic form, division ring}%
\Index%551%Q%quadratic map of division ring
   {quadratic map of division ring}%
   {Quadratic Map of Division Ring}%
\Index%180%Q%quasi affine transformation on basis manifold
   {quasi affine transformation on basis manifold}%
   {quasi affine transformation}%
\Index%181%Q%quasi affine transformation on basis manifold
   {quasi affine transformation on basis manifold}%
   {quasi affine drc transformation}%
\Index%182%Q%quasi movement on basis manifold
   {quasi movement on basis manifold}%
   {quasi movement, division ring}%
\Index%183%Q%quasi movement on basis manifold
   {quasi movement on basis manifold}%
   {quasi movement}%
\Index%103%Q%quaternion algebra E over the field F
   {quaternion algebra $E$ over the field $F$}%
   {quaternion algebra over the field}%
\Index%496%Q%quotient bundle
   {quotient bundle}%
   {quotient bundle}%
\SetIndexSpace%
\Index%59%R%r row of matrix
   {\sups row of matrix}%
   {r row}%
\Index%29%R%r rows drc vector space
   {\sups rows \drc vector space}%
   {sups rows drc vector space}%
\Index%72%R%r-row of matrix
   {$r$\hyph row of matrix}%
   {r-row}%
\Index%554%R%rank of quadratic map of division ring
   {rank of quadratic map of division ring}%
   {rank of quadratic map, division ring}%
\Index%65%R%RC inverse element of biring
   {\RC inverse element of biring}%
   {rc-inverse element}%
\Index%61%R%RC major minor
   {\RC major minor}%
   {RC-major minor}%
\Index%63%R%RC matrix group
   {\RC matrix group}%
   {rc-matrix group}%
\Index%273%R%RC nonsingular drc system of linear equations
   {\RC nonsingular \drc system of linear equations}%
   {nonsingular system of linear equations}%
\Index%64%R%RC nonsingular matrix
   {\RC nonsingular matrix}%
   {RC nonsingular matrix}%
\Index%68%R%RC power
   {\RC power}%
   {rc power}%
\Index%590%R%rc product of functional matrices
   {$\RCcirc$\Hyph product of functional matrices}%
   {rc product of functional matrices}%
\Index%66%R%RC product of matrices
   {\RC product of matrices}%
   {rc-product of matrices}%
\Index%62%R%RC quasideterminant
   {\RC quasideterminant}%
   {RC-quasideterminant}%
\Index%67%R%RC rank of matrix
   {\RC rank of matrix}%
   {rc-rank of matrix}%
\Index%60%R%RC singular matrix
   {\RC singular matrix}%
   {RC singular matrix}%
\Index%598%R%RCcirc nonsingular matrix of endomorphisms
   {$\RCcirc$\Hyph nonsingular matrix of endomorphisms}%
   {RCcirc nonsingular matrix of endomorphisms}%
\Index%601%R%RCcirc nonsingular matrix of linear mappings
   {$\RCcirc$\Hyph nonsingular matrix of linear mappings}%
   {RCcirc nonsingular matrix of linear mappings}%
\Index%602%R%RCcirc nonsingular system of linear equations
   {$\RCcirc$\Hyph nonsingular system of linear equations}%
   {RCcirc nonsingular system of linear equations}%
\Index%597%R%RCcirc quasideterminant
   {$\RCcirc$\Hyph quasideterminant}%
   {RCcirc-quasideterminant definition}%
\Index%599%R%RCcirc singular matrix of endomorphisms
   {$\RCcirc$\Hyph singular matrix of endomorphisms}%
   {RCcirc singular matrix of endomorphisms}%
\Index%600%R%RCcirc singular matrix of linear mappings
   {$\RCcirc$\Hyph singular matrix of linear mappings}%
   {RCcirc singular matrix of linear mappings}%
\Index%49%R%rcd basis dual to drc basis of vector space
   {\rcd basis dual to \drc basis of vector space}%
   {basis dual to basis, drc vector space}%
\Index%70%R%rcd representation of group
   {\rcd representation of group}%
   {rcd linear representation of group}%
\Index%69%R%rcd vector space
   {\rcd vector space}%
   {rcd vector space}%
\Index%377%R%reduced Cartesian product of bundles
   {reduced Cartesian product of bundles}%
   {reduced Cartesian product of bundles}%
\Index%378%R%reduced Cartesian product of total spaces
   {reduced Cartesian product of total spaces}%
   {reduced Cartesian product of total spaces}%
\Index%379%R%reduced fibered correspondence from \mathcal{A} to \mathcal{B}
   {reduced fibered correspondence from $\mathcal{A}$ to $\mathcal{B}$}%
   {reduced fibered correspondence from A to B}%
\Index%380%R%reduced fibered correspondence in \mathcal{A}
   {reduced fibered correspondence in $\mathcal{A}$}%
   {reduced fibered correspondence in A}%
\Index%381%R%reducible biring
   {reducible biring}%
   {reducible biring}%
\Index%450%R%reference frame in event space
   {reference frame in event space}%
   {reference frame in event space}%
\Index%263%R%reference frame manifold
   {reference frame manifold}%
   {reference frame manifold}%
\Index%439%R%reflexive 2- ary fibered relation
   {reflexive $2$\Hyph ary fibered relation}%
   {reflexive 2 ary fibered relation}%
\Index%330%R%regular endomorphism of representation
   {regular endomorphism of representation}%
   {regular endomorphism of representation}%
\Index%566%R%regular endomorphism of tower of representations
   {regular endomorphism of tower of representations}%
   {regular endomorphism of tower of representations}%
\Index%553%R%regular quadratic map in division ring
   {regular quadratic map in division ring}%
   {regular quadratic map, division ring}%
\Index%370%R%representation of \mathfrak{F}- algebra A in category B
   {representation of $\mathfrak{F}$\Hyph algebra $A$ in category $\mathcal B$}%
   {representation of F algebra in category}%
\Index%371%R%representation of \mathfrak{F}- algebra A in H algebra M
   {representation of $\mathfrak{F}$\Hyph algebra $A$ in$\mathfrak{H}$\Hyph algebra $M$}%
   {representation of algebra}%
\Index%528%R%representation of F- algebra in representation
   {representation of $\mathcal F$\Hyph algebra in representation}%
   {representation of F algebra in representation}%
\Index%529%R%representation of F- algebra in tower of representations
   {representation of $\mathcal F$\Hyph algebra in tower of representations}%
   {representation of F algebra in tower of representations}%
\Index%372%R%representation of group
   {representation of group}%
   {representation of group}%
\Index%368%R%representative of geometrical object in drc vector space
   {representative of geometrical object in \drc vector space}%
   {representative of geometrical object, drc vector space}%
\Index%369%R%representative of geometrical object in vector space
   {representative of geometrical object in vector space}%
   {representative of geometrical object, vector space}%
\Index%474%R%restriction of correspondence \Phi to set C
   {restriction of correspondence $\Phi$ to set $C$}%
   {restriction of correspondence}%
\Index%302%R%right defined Lie algebra of Lie group
   {right defined Lie algebra of Lie group}%
   {right defined Lie algebra}%
\Index%354%R%right invariant vector field
   {right invariant vector field}%
   {right invariant vector}%
\Index%361%R%right module
   {right module}%
   {right module}%
\Index%362%R%right shift on group
   {right shift on group}%
   {right shift}%
\Index%363%R%right shift on group
   {right shift on group}%
   {right shift, group}%
\Index%360%R%right structural constant of Lie algebra
   {right structural constant of Lie algebra}%
   {right structural constant of Lie algebra}%
\Index%353%R%right vector space
   {right vector space}%
   {right vector space}%
\Index%356%R%right-side contravariant representation of group
   {right-side contravariant representation of group}%
   {right-side contravariant representation of group}%
\Index%355%R%right-side covariant representation of group
   {right-side covariant representation of group}%
   {right-side covariant representation of group}%
\Index%357%R%right-side representation of \mathfrak{F}- algebra A in H algebra M
   {right-side representation of $\mathfrak{F}$\Hyph algebra $A$ in $\mathfrak{H}$\Hyph algebra $M$}%
   {right-side representation of algebra}%
\Index%358%R%right-side representation of fibered \mathfrak{F}- algebra
   {right-side representation of fibered $\mathfrak{F}$\Hyph algebra}%
   {right-side representation of fibered F-algebra}%
\Index%359%R%right-side transformation
   {right-side transformation}%
   {right-side transformation}%
\Index%188%R%ring has characteristic 0
   {ring has characteristic $0$}%
   {ring has characteristic 0}%
\Index%189%R%ring has characteristic p
   {ring has characteristic $p$}%
   {ring has characteristic p}%
\Index%127%R%row vector
   {row vector}%
   {row vector}%
\Index%71%R%Rstar- module
   {$R\star$\Hyph module}%
   {Rstar-module}%
\SetIndexSpace%
\Index%451%S%scalar potential
   {scalar potential}%
   {scalar potential}%
\Index%129%S%second Newton law
   {second Newton law}%
   {Second Newton law}%
\Index%562%S%set of coordinates of representation
   {set of coordinates of representation}%
   {coordinate set of representation}%
\Index%563%S%set of tuples of coordinates of tower of representations
   {set of tuples of coordinates of tower of representations}%
   {coordinate set of tower of representations}%
\Index%405%S%simple polyvector
   {simple polyvector}%
   {simple polyvector}%
\Index%297%S%single transitive representation of \mathfrak{F}- algebra A
   {single transitive representation of $\mathfrak{F}$\Hyph algebra $A$}%
   {single transitive representation of algebra}%
\Index%299%S%single transitive representation of fibered \mathfrak{F}- algebra
   {single transitive representation of fibered $\mathfrak{F}$\Hyph algebra}%
   {single transitive representation of fibered F-algebra}%
\Index%298%S%single transitive representation of group
   {single transitive representation of group}%
   {single transitive representation of group}%
\Index%131%S%singular linear map of division ring
   {singular linear map of division ring}%
   {singular additive map, division ring}%
\Index%130%S%singular linear mapp of D- vector space
   {singular linear mapp of $D$\Hyph vector space}%
   {singular additive map, D vector space}%
\Index%226%S%skew product of vectors
   {skew product of vectors}%
   {skew product of vectors}%
\Index%406%S%space of orbits of Tstar- representation
   {space of orbits of \Ts representation}%
   {space of orbits of Ts representation}%
\Index%452%S%speed of deviation
   {speed of deviation}%
   {speed of deviation}%
\Index%74%S%SRCstar,TRCstar- linear map of vector bundles
   {$(\mathcal S\RCstar,\mathcal T\RCstar)$\Hyph linear map of vector bundles}%
   {src trc linear map of vector bundles}%
\Index%73%S%SRCstar,TRCstar- linear map of vector spaces
   {$(S\RCstar,T\RCstar)$\Hyph linear map of vector spaces}%
   {src trc linear map of vector spaces}%
\Index%112%S%Sstar, star T-bimodule
   {($S\star$, $\star T$)\hyph bimodule}%
   {(Sstar,starT)-bimodule}%
\Index%145%S%stability group
   {stability group}%
   {stability group}%
\Index%457%S%stable set of representation
   {stable set of representation}%
   {stable set of representation}%
\Index%458%S%standard component of Gateaux differential of map f
   {standard component of the G\^ateaux differential of map $f$}%
   {standard component of Gateaux differential, division ring}%
\Index%460%S%standard component of linear map of division ring
   {standard component of linear map of division ring}%
   {standard component of linear map, division ring}%
\Index%462%S%standard component of polylinear map f of division ring
   {standard component of polylinear map $f$ of division ring}%
   {standard component of polylinear map, division ring}%
\Index%459%S%standard component of quadratic map f over field F
   {standard component of quadratic map $f$ over field $F$}%
   {standard component of quadratic map, division ring}%
\Index%463%S%standard component of tensor
   {standard component of tensor}%
   {standard component of tensor, division ring}%
\Index%574%S%standard component of tensor in tensor product of algebras
   {standard component of tensor in tensor product of algebras}%
   {standard component of tensor, algebra}%
\Index%461%S%standard component over field F of bilitnear map f
   {standard component over field $F$ of bilitnear map $f$}%
   {standard component of bilinear map, division ring}%
\Index%472%S%standard coordinates of basis
   {standard coordinates of basis}%
   {standard coordinates of basis}%
\Index%471%S%standard coordinates of drc basis
   {standard coordinates of \drc basis}%
   {standard coordinates of drc basis}%
\Index%456%S%standard F- component of additive map f
   {standard $F$\Hyph component of additive map $f$}%
   {standard component of additive map, division ring}%
\Index%464%S%standard F- representation of additive map of division ring
   {standard $F$\Hyph representation of additive map of division ring}%
   {additive map, standard representation, division ring}%
\Index%465%S%standard representation of Gateaux differential of map of division ring over field F
   {standard representation of the G\^ateaux differential of map of division ring over field $F$}%
   {Gateaux differential, standard representation, division ring}%
\Index%467%S%standard representation of linear map of division ring
   {standard representation of linear map of division ring}%
   {linear map, standard representation, division ring}%
\Index%468%S%standard representation of matrix
   {standard representation of matrix}%
   {Standard representation}%
\Index%470%S%standard representation of polylinear map of division ring
   {standard representation of polylinear map of division ring}%
   {polylinear map, standard representation, division ring}%
\Index%466%S%standard representation of quadratic map of division ring over field F
   {standard representation of quadratic map of division ring over field $F$}%
   {quadratic map, standard representation, division ring}%
\Index%469%S%standard representation over field F of bilinear map of division ring
   {standard representation over field $F$ of bilinear map of division ring}%
   {bilinear map, standard representation, division ring}%
\Index%586%S%star A- product of Astar linear map over scalar
   {$\star A$\Hyph product of $A\star$\Hyph linear mapping over scalar}%
   {starA product of Astar linear map over scalar}%
\Index%77%S%star D- product of drc linear map A over scalar
   {$\star D$\Hyph product of \drc linear map $A$ over scalar}%
   {starD product of drc linear map over scalar}%
\Index%78%S%star D-vector space
   {$\star D$\hyph vector space}%
   {starD-vector space}%
\Index%79%S%star R-module
   {$\star R$\hyph module}%
   {starR-module}%
\Index%75%S%starD component of coordinates of vector \Vector r
   {\sD component of coordinates of vector $\Vector r$}%
   {starD component of coordinates of vector, D vector space}%
\Index%82%S%starT- representation of \mathfrak{F}- algebra A in H algebra M
   {\sT representation of $\mathfrak{F}$\Hyph algebra $A$ in $\mathfrak{H}$\Hyph algebra $M$}%
   {starT representation of algebra}%
\Index%83%S%starT- representation of fibered \mathfrak{F}- algebra
   {\sT representation of fibered $\mathfrak{F}$\Hyph algebra}%
   {starT representation of fibered F-algebra}%
\Index%84%S%starT- representation of fibered group
   {\sT representation of fibered group}%
   {starT representation of fibered group}%
\Index%80%S%starT- shift
   {\sT shift}%
   {starT shift}%
\Index%81%S%starT- shift on fibered group
   {\sT shift on fibered group}%
   {starT shift, fibered group}%
\Index%85%S%starT- transformation
   {\sT transformation}%
   {starT transformation}%
\Index%86%S%starT- transformation on bundle
   {\sT transformation on bundle}%
   {starT transformation of bundle}%
\Index%473%S%structural constants of division ring D over field F
   {structural constants of division ring $D$ over field $F$}%
   {structural constants of division ring over field}%
\Index%340%S%subbundle
   {subbundle}%
   {subbundle}%
\Index%341%S%subbundle of Dstar-vector space
   {subbundle of $\mathcal D\star$\hyph vector space}%
   {subbundle of Dstar vector bundle}%
\Index%322%S%subrepresentation generated by set X
   {subrepresentation generated by set $X$}%
   {subrepresentation generated by set}%
\Index%334%S%subrepresentation of representation
   {subrepresentation of representation}%
   {subrepresentation of representation}%
\Index%475%S%sum of drc linear maps
   {sum of \drc linear maps}%
   {sum of drc linear maps, drc vector spaces}%
\Index%476%S%sum of geometrical objects
   {sum of geometrical objects}%
   {sum of geometrical objects}%
\Index%477%S%sum of geometrical objects in vector space
   {sum of geometrical objects in vector space}%
   {sum of geometrical objects, vector space}%
\Index%446%S%symmetric 2- ary fibered relation
   {symmetric $2$\Hyph ary fibered relation}%
   {symmetric 2 ary fibered relation}%
\Index%447%S%symmetric bilinear map of D- vector space to division ring
   {symmetric bilinear map of $D$\Hyph vector space to division ring}%
   {symmetric bilinear map, vector space to division ring}%
\Index%143%S%symmetry group
   {symmetry group}%
   {symmetry group}%
\Index%144%S%symmetry group
   {symmetry group}%
   {SymmetryGroup}%
\Index%448%S%synchronization of reference frame
   {synchronization of reference frame}%
   {synchronization of reference frame}%
\SetIndexSpace%
\Index%499%T%Taylor polynomial
   {Taylor polynomial}%
   {Taylor polynomial, division ring}%
\Index%440%T%Taylor series
   {Taylor series}%
   {Taylor series, division ring}%
\Index%573%T%tensor product of algebras
   {tensor product of algebras}%
   {tensor product of algebras}%
\Index%480%T%tensor product of D- vector spaces
   {tensor product of $D$\Hyph vector spaces}%
   {tensor product of D vector spaces}%
\Index%484%T%tensor product of division rings
   {tensor product of division rings}%
   {tensor product of division rings}%
\Index%481%T%tensor product of Dstar vector spaces
   {tensor product of \Ds vector spaces}%
   {tensor product of Dstar vector spaces}%
\Index%483%T%tensor product of representations
   {tensor product of representations}%
   {tensor product of representations}%
\Index%482%T%tensor product of rings over commutative ring
   {tensor product of rings over commutative ring}%
   {tensor product of rings}%
\Index%47%T%the Fr\'echet Dstar derivative of map f of division ring D at point x
   {the Fr\'echet \Ds derivative of map $f$ of division ring $D$ at point $x$}%
   {Frechet Dstar derivative of map, division ring}%
\Index%486%T%topological D- vector space
   {topological $D$\Hyph vector space}%
   {topological D vector space}%
\Index%488%T%topological division ring
   {topological division ring}%
   {topological division ring}%
\Index%487%T%topological drc vector space
   {topological \drc vector space}%
   {topological drc vector space}%
\Index%498%T%torsion form
   {torsion form}%
   {torsion form}%
\Index%478%T%torsion tensor
   {torsion tensor}%
   {torsion tensor}%
\Index%120%T%tower of bundles
   {tower of bundles}%
   {tower of bundles}%
\Index%119%T%tower of representations of \mathfrak{F}- algebras
   {tower of representations of $\Vector{\mathfrak{F}}$\Hyph algebras}%
   {tower of representations of algebras}%
\Index%531%T%tower of subrepresentations
   {tower of subrepresentations}%
   {tower of subrepresentations}%
\Index%532%T%tower of subrepresentations of tower of representations generated by tuple of sets
   {tower of subrepresentations of tower of representations $\Vector f$ generated by tuple of sets $\VX X$}%
   {subrepresentation generated by tuple of sets}%
\Index%376%T%transformation coordinated with equivalence
   {transformation coordinated with equivalence}%
   {transformation coordinated with equivalence}%
\Index%374%T%transformation of H algebra
   {transformation of $\mathfrak{H}$\Hyph algebra}%
   {transformation of H algebra}%
\Index%375%T%transformation on bundle
   {transformation on bundle}%
   {transformation of bundle}%
\Index%489%T%transitive 2- ary fibered relation
   {transitive $2$\Hyph ary fibered relation}%
   {transitive 2 ary fibered relation}%
\Index%490%T%transitive representation of \mathfrak{F}- algebra A
   {transitive representation of $\mathfrak{F}$\Hyph algebra $A$}%
   {transitive representation of algebra}%
\Index%492%T%transitive representation of fibered \mathfrak{F}- algebra
   {transitive representation of fibered $\mathfrak{F}$\Hyph algebra}%
   {transitive representation of fibered F-algebra}%
\Index%491%T%transitive representation of group
   {transitive representation of group}%
   {transitive representation of group}%
\Index%88%T%Tstar- linear composition of  vectors
   {\Ts linear composition of  vectors}%
   {linear composition of  vectors}%
\Index%87%T%Tstar- matrices vector space
   {\Ts matrices vector space}%
   {matrices vector space}%
\Index%90%T%Tstar- representation of \mathfrak{F}- algebra A in H algebra M
   {\Ts representation of $\mathfrak{F}$\Hyph algebra $A$ in $\mathfrak{H}$\Hyph algebra $M$}%
   {Tstar representation of algebra}%
\Index%91%T%Tstar- representation of fibered \mathfrak{F}- algebra
   {\Ts representation of fibered $\mathfrak{F}$\Hyph algebra}%
   {Tstar representation of fibered F-algebra}%
\Index%89%T%Tstar- shift
   {\Ts shift}%
   {Tstar shift}%
\Index%92%T%Tstar- transformation
   {\Ts transformation}%
   {Tstar transformation}%
\Index%93%T%Tstar- transformation on bundle
   {\Ts transformation on bundle}%
   {Tstar transformation of bundle}%
\Index%535%T%tuple of coordinates of element relative to tuple of sets
   {tuple of coordinates of element $\Vector a$ relative to tuple of sets $\VX X$}%
   {coordinates of element, tower of representations}%
\Index%534%T%tuple of generating sets of tower of representations
   {tuple of generating sets of tower of representations}%
   {tuple of generating sets of tower of representations}%
\Index%533%T%tuple of generating sets of tower subrepresentations
   {tuple of generating sets of tower subrepresentations}%
   {tuple of generating sets of tower subrepresentations}%
\Index%530%T%tuple of stable sets of tower of representation
   {tuple of stable sets of tower of representation}%
   {tuple of stable sets of tower of representations}%
\Index%587%T%twin representations of associative algebra
   {twin representations of associative algebra}%
   {twin representations of associative algebra}%
\Index%320%T%twin representations of division ring
   {twin representations of division ring}%
   {twin representations of division ring}%
\Index%319%T%twin representations of fibered group
   {twin representations of fibered group}%
   {twin representations of fibered group}%
\Index%318%T%twin representations of group
   {twin representations of group}%
   {twin representations of group}%
\SetIndexSpace%
\Index%160%U%unit sphere in division ring
   {unit sphere in division ring}%
   {unit sphere in division ring}%
\Index%172%U%unitarity law for  Dstar- vector fields
   {unitarity law for  $\mathcal D\star$\Hyph vector fields}%
   {unitarity law, Dstar vector fields}%
\Index%578%U%unitarity law for Astar- vector space
   {unitarity law for $A\star$\Hyph vector space}%
   {unitarity law, Astar vector space}%
\Index%171%U%unitarity law for Dstar- vector space
   {unitarity law for $D\star$\Hyph vector space}%
   {unitarity law, Dstar vector space}%
\SetIndexSpace%
\Index%288%V%valued division ring
   {valued division ring}%
   {valued division ring}%
\Index%125%V%vector bundle
   {vector bundle}%
   {vector bundle}%
\Index%126%V%vector potential
   {vector potential}%
   {vector potential}%
\Index%485%V%vector space type
   {vector space type}%
   {vector space type}%

\CloseIndex

%auto-ignore
\def\indexname{Special Symbols and Notations}
\OpenIndex

\SetIndexSpace%A%0
\Symb%A/0
   {$(^a_b)$\hyph\CR quasideterminant}%
   {a b CR quasideterminant definition}%
\Symb%A/0
   {$(^a_b)$\hyph \RC quasideterminant}%
   {a b RC-quasideterminant definition}%
\Symb%A/0
   {$(^a_b)$\hyph $\RCcirc$\Hyph quasideterminant}%
   {a b RCcirc-quasideterminant definition}%
\Symb%A/0
   {minor}%
   {A from b a}%
\Symb%A/0
   {minor}%
   {A from columns T}%
\Symb%A/0
   {minor}%
   {A from rows S}%
\Symb%A/0
   {minor}%
   {A without column a}%
\Symb%A/0
   {minor}%
   {A without columns T}%
\Symb%A/0
   {minor}%
   {A without row b}%
\Symb%A/0
   {minor}%
   {A without rows S}%
\Symb%A/0
   {affine space}%
   {affine space, division ring}%
\Symb%A/0
   {affine space}%
   {An}%
\Symb%A/0
   {\subs row ($c$\hyph row) of matrix}%
   {c row}%
\Symb%A/0
   {component of linear map $\Vector{A}$ of $D$\Hyph vector space}%
   {component of linear map, D vector space}%
\Symb%A/0
   {component $p$ of polyadditive map $\Vector A$}%
   {component of polyadditive map, D vector space}%
\Symb%A/0
   {\CR power of element $A$ of biring}%
   {cr power}%
\Symb%A/0
   {\CR inverse element of biring}%
   {cr-inverse element}%
\Symb%A/0
   {\CR product of matrices}%
   {cr-product of matrices}%
\Symb%A/0
   {\dcr vector}%
   {dcr vector}%
\Symb%A/0
   {derivative of left shift}%
   {derivative of left shift}%
\Symb%A/0
   {derivative of left shift}%
   {derivative of left shift, 1-Parameter Group}%
\Symb%A/0
   {}%
   {derivative of right shift}%
\Symb%A/0
   {}%
   {derivative of right shift}%
\Symb%A/0
   {derivative of right shift}%
   {derivative of right shift, 1-Parameter Group}%
\Symb%A/0
   {derivative of left shift}%
   {derivative of Tstar shift}%
\Symb%A/0
   {\drc vector}%
   {drc vector}%
\Symb%A/0
   {hermitian conjugation in division ring}%
   {hermitian conjugation, division ring}%
\Symb%A/0
   {loop of automorphisms of representation $f$}%
   {loop of automorphisms of representation}%
\Symb%A/0
   {norm of map $\Vector A$ of normed $D$\hyph vector space}%
   {norm of map, D vector space}%
\Symb%A/0
   {derivative}%
   {overline nabla_l, definition 2}%
\Symb%A/0
   {partial additive map of variable $v^i$}%
   {partial additive map of variable}%
\Symb%A/0
   {\sups row ($r$\hyph row) of matrix}%
   {r row}%
\Symb%A/0
   {\RC power of element $A$ of biring}%
   {rc power}%
\Symb%A/0
   {\RC inverse element of biring}%
   {rc-inverse element}%
\Symb%A/0
   {\RC product of matrices}%
   {rc-product of matrices}%
\Symb%A/0
   {\RC quasideterminant}%
   {RC-quasideterminant definition}%
\Symb%A/0
   {$\RCcirc$\Hyph quasideterminant}%
   {RCcirc-quasideterminant definition}%
\Symb%A/0
   {set of additive maps of $D$\Hyph vector space $\Vector{V}$ into $D$\Hyph vector space $\Vector{W}$}%
   {set additive maps, D vector space}%
\Symb%A/0
   {set of additive maps of ring $R_1$ into ring $R_2$}%
   {set additive maps, ring}%
\Symb%A/0
   {set of polyadditive mappings}%
   {set polyadditive maps, D vector space}%
\Symb%A/0
   {set of polyadditive maps of rings $R_1$, ..., $R_n$ into module $S$}%
   {set polyadditive maps, ring}%
\Symb%A/0
   {set of polylinear maps of rings $R_1$, ..., $R_n$ into module $S$}%
   {set polylinear maps, ring}%
\Symb%A/0
   {skew product of vectors $\Vector a_1$, ..., $\Vector a_m$}%
   {skew product of vectors}%
\Symb%A/0
   {standard component of tensor in tensor product of algebras}%
   {standard component of tensor, algebra}%
\Symb%A/0
   {right shift}%
   {starT shift}%
\Symb%A/0
   {\sT shift}%
   {starT shift, fibered group}%
\Symb%A/0
   {tensor product of algebras}%
   {tensor product of algebras}%
\Symb%A/0
   {left shift}%
   {Tstar shift}%
\Symb%A/0
   {\Ts shift}%
   {Tstar shift, fibered group}%
\Symb%A/0
   {anholonomic coordinates of vector}%
   {vector anholonomic coordinates}%
\Symb%A/0
   {holonomic coordinates of vector}%
   {vector holonomic coordinates}%

\SetIndexSpace%B%0
\Symb%B/0
   {basis manifold of \drc vector space $\mathcal{V}$}%
   {basis manifold of drc vector space}%
\Symb%B/0
   {basis manifold of vector space}%
   {basis manifold of vector space}%
\Symb%B/0
   {basis manifold of representation $f$}%
   {basis manifold representation F algebra}%
\Symb%B/0
   {basis manifold of affine space}%
   {Basis Manifold, Affine Space}%
\Symb%B/0
   {basis manifold of central affine space}%
   {BCAn}%
\Symb%B/0
   {basis manifold of Euclid space}%
   {BEn}%
\Symb%B/0
   {Cartesian power of bundle}%
   {Cartesian power of bundle}%
\Symb%B/0
   {Cartesian power $A$ of set $B$}%
   {Cartesian power of set}%
\Symb%B/0
   {\drc basis manifold of affine space}%
   {drc Basis Manifold, Affine Space, division ring}%
\Symb%B/0
   {basis manifold of central affine space}%
   {FCAn}%
\Symb%B/0
   {basis manifold of Euclid space}%
   {FEn}%
\Symb%B/0
   {lattice of subrepresentations of representation $f$}%
   {lattice of subrepresentations}%
\Symb%B/0
   {lattice of towers of subrepresentations of tower of representations $\Vector f$}%
   {lattice of subrepresentations, tower of representations}%
\Symb%B/0
   {product of objects $B_1$, ..., $B_n$ in category $\mathcal A$}%
   {product of objects in category, 1 n}%
\Symb%B/0
   {structural constants of division ring $D$ over field $F$}%
   {structural constants of division ring over field}%

\SetIndexSpace%C%0
\Symb%C/0
   {central affine space}%
   {CAn}%
\Symb%C/0
   {central affine space}%
   {central affine space}%
\Symb%C/0
   {$\CRcirc$\Hyph product of functional matrices}%
   {cr product of functional matrices}%
\Symb%C/0
   {left structural constant of Lie algebra}%
   {left structural constant of Lie algebra}%
\Symb%C/0
   {right structural constant of Lie algebra}%
   {right structural constant of Lie algebra}%

\SetIndexSpace%D%0
\Symb%D/0
   {basis vector of \sT representation}%
   {basis vector of starT representation}%
\Symb%D/0
   {basis vector of \sT representation}%
   {basis vector of starT representation, coordinates}%
\Symb%D/0
   {basis vector of \Ts representation}%
   {basis vector of Tstar representation}%
\Symb%D/0
   {basis vector of \Ts representation}%
   {basis vector of Tstar representation, coordinates}%
\Symb%D/0
   {\subs rows \dcr vector space}%
   {c rows dcr vector space}%
\Symb%D/0
   {component of the G\^ateaux derivative of map $\Vector f(\Vector x)$}%
   {component of Gateaux derivative of map, D vector space}%
\Symb%D/0
   {component of the G\^ateaux derivative of map $\Vector f(\Vector x)$}%
   {component of Gateaux derivative of map, D vector space, short}%
\Symb%D/0
   {component of the G\^ateaux differential of map $f(x)$}%
   {component of Gateaux derivative of map, division ring}%
\Symb%D/0
   {component of the G\^ateaux derivative of second order of map $\Vector f(\Vector x)$}%
   {component of Gateaux derivative of Second Order, D vector space}%
\Symb%D/0
   {component of the G\^ateaux derivative of second order of map $f(x)$ of division ring}%
   {component of Gateaux derivative of Second Order, division ring}%
\Symb%D/0
   {coordinate \Drc vector bundle}%
   {coordinate drc vector bundle}%
\Symb%D/0
   {coordinate \drc vector space}%
   {coordinate drc vector space}%
\Symb%D/0
   {coordinate reference frame}%
   {coordinate reference frame, extensive definition}%
\Symb%D/0
   {diagonal in bundle $\mathcal{A}$}%
   {diagonal in bundle, 1}%
\Symb%D/0
   {direct product of division rings $D_1$, ..., $D_n$}%
   {direct product of division rings, 1 n}%
\Symb%D/0
   {the Fr\'echet \Ds derivative of map $f$ of division ring}%
   {Frechet Dstar derivative of map, division ring}%
\Symb%D/0
   {the G\^ateaux \crd derivative of map $\Vector f$ of $D$\hyph vector space $\Vector V$ to $D$\hyph vector space $\Vector W$}%
   {Gateaux crd derivative of map, D vector space}%
\Symb%D/0
   {the G\^ateaux derivative of map $\Vector f$ of normed $D$\Hyph vector space $\Vector{V}$ to normed $D$\Hyph vector space $\Vector{W}$}%
   {Gateaux derivative of map, D vector space}%
\Symb%D/0
   {the G\^ateaux derivative of map $f$}%
   {Gateaux derivative of map, division ring}%
\Symb%D/0
   {the G\^ateaux derivative of map $f$}%
   {Gateaux derivative of map, fraction, division ring}%
\Symb%D/0
   {the G\^ateaux derivative of order $n$ of map $\Vector f$}%
   {Gateaux derivative of Order n, D vector space}%
\Symb%D/0
   {the G\^ateaux derivative of order $n$ of map $f$ of division ring}%
   {Gateaux derivative of Order n, division ring}%
\Symb%D/0
   {the G\^ateaux derivative of order $n$ of map $f$ of division ring}%
   {Gateaux derivative of Order n, fraction, division ring}%
\Symb%D/0
   {the G\^ateaux derivative of second order of map $\Vector f$}%
   {Gateaux derivative of Second Order, D vector space}%
\Symb%D/0
   {the G\^ateaux derivative of second order of map $f$ of division ring}%
   {Gateaux derivative of Second Order, division ring}%
\Symb%D/0
   {the G\^ateaux derivative of second order of map $f$ of division ring}%
   {Gateaux derivative of Second Order, fraction, division ring}%
\Symb%D/0
   {the G\^ateaux differential of map $\Vector f$ of normed $D$\Hyph vector space $\Vector{V}$ to normed $D$\Hyph vector space $\Vector{W}$}%
   {Gateaux differential of map, D vector space}%
\Symb%D/0
   {the G\^ateaux differential of map $f$}%
   {Gateaux differential of map, division ring}%
\Symb%D/0
   {the G\^ateaux differential of second order of mapping $\Vector f$}%
   {Gateaux differential of Second Order, D vector space}%
\Symb%D/0
   {the G\^ateaux differential of second order of mapping $f$ of division ring}%
   {Gateaux differential of Second Order, division ring}%
\Symb%D/0
   {the G\^ateaux \drc derivative of map $\Vector f$ of $D$\Hyph vector space $\Vector V$ to $D$\Hyph vector space $\Vector W$}%
   {Gateaux drc derivative of map, D vector space}%
\Symb%D/0
   {the G\^ateaux \Ds derivative of map $f$ of division ring $D$}%
   {Gateaux Dstar derivative of map, division ring}%
\Symb%D/0
   {the G\^ateaux Jacobian of map of $D$\Hyph vector space}%
   {Gateaux Jacobian of map, D vector space}%
\Symb%D/0
   {the G\^ateaux partial \drc derivative of map $f^b$ with respect to variable $v^a$}%
   {Gateaux partial crd derivative of map, 1, D vector space}%
\Symb%D/0
   {the G\^ateaux partial \drc derivative of map $f^b$ with respect to variable $v^a$}%
   {Gateaux partial crd derivative of map, 2, D vector space}%
\Symb%D/0
   {the G\^ateaux partial \drc derivative of map $f^b$ with respect to variable $v^a$}%
   {Gateaux partial crd derivative of map, 3, D vector space}%
\Symb%D/0
   {the G\^ateaux mixed partial derivative of map $f^j$ with respect to variables $v^i$, $v^j$}%
   {Gateaux partial derivative of Second Order, D vector space}%
\Symb%D/0
   {the G\^ateaux partial derivative of map $f^j$ with respect to variable $v^i$}%
   {Gateaux partial derivative, D vector space}%
\Symb%D/0
   {the G\^ateaux partial \drc derivative of map $f^b$ with respect to variable $v^a$}%
   {Gateaux partial drc derivative of map, 1, D vector space}%
\Symb%D/0
   {the G\^ateaux partial \drc derivative of map $f^b$ with respect to variable $v^a$}%
   {Gateaux partial drc derivative of map, 2, D vector space}%
\Symb%D/0
   {the G\^ateaux partial \drc derivative of map $f^b$ with respect to variable $v^a$}%
   {Gateaux partial drc derivative of map, 3, D vector space}%
\Symb%D/0
   {the G\^ateaux \sD derivative of map $f$ of division ring $D$}%
   {Gateaux starD derivative of map, division ring}%
\Symb%D/0
   {matrices vector space}%
   {matrices vector space}%
\Symb%D/0
   {Cartan derivative}%
   {overbrace D}%
\Symb%D/0
   {derivative}%
   {overline D}%
\Symb%D/0
   {derivative $e_{(k)}$}%
   {partial(k)}%
\Symb%D/0
   {\sups rows \drc vector space}%
   {r rows drc vector space}%
\Symb%D/0
   {speed of deviation}%
   {speed of deviation}%
\Symb%D/0
   {standard component of the G\^ateaux differential of map $f$}%
   {standard component of Gateaux differential, division ring}%
\Symb%D/0
   {tensor product of division rings}%
   {tensor product of division rings}%
\Symb%D/0
   {tensor product of rings}%
   {tensor product of rings}%
\Symb%D/0
   {vector space type}%
   {vector space type}%

\SetIndexSpace%E%0
\Symb%E/0
   {affine basis}%
   {Affine Basis}%
\Symb%E/0
   {basis of vector space}%
   {Basis e}%
\Symb%E/0
   {basis in vector space $\mathcal{V}$}%
   {basis in V}%
\Symb%E/0
   {basis of vector space}%
   {basis, vector space}%
\Symb%E/0
   {basis of $(n)$\hyph vector space}%
   {basis,n vector space}%
\Symb%E/0
   {Cartesian power of total spaces}%
   {Cartesian power of total spaces}%
\Symb%E/0
   {Cartesian product of total spaces}%
   {Cartesian product of total spaces, definition 1}%
\Symb%E/0
   {central affine basis}%
   {Central Affine Basis}%
\Symb%E/0
   {affine basis}%
   {drc affine basis, division ring}%
\Symb%E/0
   {basis for \Drc vector bundle}%
   {drc basis, vector bundle}%
\Symb%E/0
   {form of reference frame}%
   {dual forms, reference frame}%
\Symb%E/0
   {Euclid space}%
   {Euclid space}%
\Symb%E/0
   {Euclid space}%
   {Euclid space, division ring}%
\Symb%E/0
   {identical transformation of bundle}%
   {identical transformation of bundle}%
\Symb%E/0
   {orthonornal basis}%
   {Orthonornal Basis}%
\Symb%E/0
   {pseudo Euclid space}%
   {pseudo Euclid space}%
\Symb%E/0
   {pseudo Euclid space}%
   {pseudo Euclid space, division ring}%
\Symb%E/0
   {quaternion algebra over the field $F$}%
   {quaternion algebra over the field}%
\Symb%E/0
   {quaternion division algebra over the field}%
   {quaternion division algebra over the fieldL}%
\Symb%E/0
   {reduced Cartesian product of total spaces}%
   {reduced Cartesian product of total spaces, definition 1}%
\Symb%E/0
   {set of nonsingular \sT transformations of bundle $\mathcal{E}$}%
   {set of starT nonsingular transformations of bundle}%
\Symb%E/0
   {set of nonsingular \Ts transformations of bundle $\mathcal{E}$}%
   {set of Tstar nonsingular transformations of bundle}%
\Symb%E/0
   {standard coordinates of basis}%
   {standard coordinates of basis}%
\Symb%E/0
   {standard coordinates of reference frame}%
   {standard coordinates of reference frame}%
\Symb%E/0
   {vector field of reference frame}%
   {vector field of reference frame}%
\Symb%E/0
   {vector of basis}%
   {vector of basis}%

\SetIndexSpace%F%0
\Symb%F/0
   {active representation of loop $\mathfrak A(f)$ in basis manifold $\mathcal B(f)$}%
   {active representation in basis manifold}%
\Symb%F/0
   {coordinates of basis in \subs rows \dcr vector space}%
   {basis coordinates, c rows dcr vector space}%
\Symb%F/0
   {coordinates of basis in \sups rows \drc vector space}%
   {basis coordinates, r rows drc vector space}%
\Symb%F/0
   {basis for \subs rows \dcr vector space}%
   {basis, c rows dcr vector space}%
\Symb%F/0
   {basis for \sups rows \drc vector space}%
   {basis, r rows drc vector space}%
\Symb%F/0
   {central affine basis}%
   {Central Affine Basis, division ring}%
\Symb%F/0
   {component of linear map $f$ of division ring}%
   {component of linear map, division ring}%
\Symb%F/0
   {component of polylinear map of division ring}%
   {component of polylinear map, division ring}%
\Symb%F/0
   {fibered morphism from bundle $\mathcal{A}$ into $\mathcal{B}$}%
   {fibered morphism from A into B}%
\Symb%F/0
   {filter $\mathfrak{F}$ converges to set $A$}%
   {filter converges}%
\Symb%F/0
   {homomorphism of fibered $\mathfrak{F}$\Hyph algebras}%
   {homomorphism of fibered F-algebras}%
\Symb%F/0
   {inverse fibered correspondence}%
   {inverse fibered correspondence, 1}%
\Symb%F/0
   {inverse reduced fibered correspondence}%
   {inverse reduced fibered correspondence, 1}%
\Symb%F/0
   {map to Cartesian product}%
   {map to Cartesian product}%
\Symb%F/0
   {norm of map $f$ of division ring}%
   {norm of map, division ring}%
\Symb%F/0
   {representation orbit of group $G$}%
   {orbit of Tstar representation of group}%
\Symb%F/0
   {orthonornal basis}%
   {Orthonornal Basis, division ring}%
\Symb%F/0
   {reference frame}%
   {reference frame}%
\Symb%F/0
   {reference frame, extensive definition}%
   {reference frame, extensive definition}%
\Symb%F/0
   {standard $F$\Hyph component of additive map $f$}%
   {standard component of additive map, division ring}%
\Symb%F/0
   {standard component of biadditive map $f$ over field $F$}%
   {standard component of biadditive map, division ring}%
\Symb%F/0
   {standard component of linear map $f$ of division ring}%
   {standard component of linear map, division ring}%
\Symb%F/0
   {standard component of polylinear map $f$ of division ring}%
   {standard component of polylinear map, division ring}%
\Symb%F/0
   {standard component of quadratic map $f$ over field $F$}%
   {standard component of quadratic map, division ring}%
\Symb%F/0
   {standard component of tensor}%
   {standard component of tensor, division ring}%

\SetIndexSpace%G%0
\Symb%G/0
   {affine transformation group}%
   {affine transformation group}%
\Symb%G/0
   {\CR matrix group}%
   {cr-matrix group}%
\Symb%G/0
   {affine transformation group}%
   {drc affine transformation group}%
\Symb%G/0
   {fibered little group of section $h$}%
   {fibered little group}%
\Symb%G/0
   {fibered stability group of section $h$}%
   {fibered stability group}%
\Symb%G/0
   {algebra Lie of group Lie}%
   {g}%
\Symb%G/0
   {left defined algebra Lie of group Lie}%
   {gl}%
\Symb%G/0
   {right defined algebra Lie of group Lie}%
   {gr}%
\Symb%G/0
   {group of homomorphisms of vector space $\mathcal{V}$}%
   {GV}%
\Symb%G/0
   {little group of $x$}%
   {little group}%
\Symb%G/0
   {orbit of effective covariant \sT representation of fibered group}%
   {orbit of effective covariant starT representation of fibered group}%
\Symb%G/0
   {orbit of effective covariant \sT representation of group}%
   {orbit of effective covariant starT representation of group}%
\Symb%G/0
   {orbit of effective covariant \Ts representation of fibered group}%
   {orbit of effective covariant Tstar representation of fibered group}%
\Symb%G/0
   {orbit of effective covariant \Ts representation of group}%
   {orbit of effective covariant Tstar representation of group}%
\Symb%G/0
   {product of groups $G_1$, ..., $G_n$}%
   {product of groups, 1 n}%
\Symb%G/0
   {\RC matrix group}%
   {rc-matrix group}%
\Symb%G/0
   {stability group of $x$}%
   {stability group}%

\SetIndexSpace%H%0
\Symb%H/0
   {Hadamard inverse of matrix}%
   {Hadamard inverse of matrix}%
\Symb%H/0
   {quaternion algebra over real field}%
   {quaternion algebra over real field}%

\SetIndexSpace%I%0
\Symb%I/0
   {infinitesimal generator of representation}%
   {infinitesimal generator of representation}%
\Symb%I/0
   {Lie group infinitesimal generators}%
   {Lie group infinitesimal generators}%

\SetIndexSpace%J%0
\Symb%J/0
   {closure operator of representation $f$}%
   {closure operator, representation}%
\Symb%J/0
   {closure operator of tower of representations $\Vector f$}%
   {closure operator, tower of representations}%
\Symb%J/0
   {tower of subrepresentations of tower of representations $\Vector f$ generated by tuple of sets $\VX X$}%
   {subrepresentation generated by tuple of sets}%

\SetIndexSpace%K%0
\Symb%K/0
   {kernel of additive map of $D$\Hyph vector space}%
   {kernel of additive map, D vector space}%
\Symb%K/0
   {kernel of additive map of division ring}%
   {kernel of additive map, division ring}%

\SetIndexSpace%L%0
\Symb%L/0
   {left shift}%
   {left shift}%
\Symb%L/0
   {Lie derivative of connection}%
   {Lie derivative of connection}%
\Symb%L/0
   {Lie derivative of metric}%
   {Lie derivative of metric}%
\Symb%L/0
   {limit of correspondence $\Phi$ with respect to the filter $\mathfrak{F}$}%
   {limit of correspondence with respect to the filter}%
\Symb%L/0
   {limit of sequence in valued division ring}%
   {limit of sequence, valued division ring}%
\Symb%L/0
   {passive transformation}%
   {passive transformation}%
\Symb%L/0
   {\rcd vector space of \drc linear maps of \drc vector space $\Vector{V}$ into \drc vector space $\Vector{W}$}%
   {set drc linear maps, drc vector space}%
\Symb%L/0
   {set of linear maps of $D$\Hyph vector space $\Vector{V}$ into $D$\Hyph vector space $\Vector{W}$}%
   {set linear maps, D vector space}%
\Symb%L/0
   {set of left-side nonsingular transformations of set $M$}%
   {set of left-side nonsingular transformations}%
\Symb%L/0
   {set of polylinear maps of algebras $A_1$, ..., $A_n$ into algebra $A$}%
   {set polylinear maps, finite dimensional algebra}%
\Symb%L/0
   {\drc vector space of \rcd linear maps of \rcd vector space $\Vector{V}$ into \rcd vector space $\Vector{W}$}%
   {set rcd linear maps, rcd vector space}%
\Symb%L/0
   {set of \sT representations of division ring $S$ in additive group of division ring $R$}%
   {set sT representations, division ring}%
\Symb%L/0
   {set of \Ts representations of division ring $S$ in additive group of division ring $R$}%
   {set Ts representations, division ring}%

\SetIndexSpace%M%0
\Symb%M/0
   {set of \sT transformations of set $M$}%
   {set of starT transformations}%
\Symb%M/0
   {set of transformations of set $M$}%
   {set of transformations}%
\Symb%M/0
   {set of \Ts transformations of set $M$}%
   {set of Tstar transformations}%
\Symb%M/0
   {space of orbits of effective \sT covariant representation of the group}%
   {space of orbits of effective sT representation}%
\Symb%M/0
   {space of orbits of effective \Ts covariant representation of the group}%
   {space of orbits of effective Ts representation}%
\Symb%M/0
   {space of orbits of \Ts representation $f$ of group $G$ in set $M$}%
   {space of orbits of Ts representation}%

\SetIndexSpace%O%0
\Symb%O/0
   {geometrical object in coordinate representation defined in \drc vector space}%
   {geometrical object, coordinate drc vector space}%
\Symb%O/0
   {geometrical object in coordinate representation}%
   {geometrical object, coordinate vector space}%
\Symb%O/0
   {geometrical object defined in \drc vector space}%
   {geometrical object, drc vector space}%
\Symb%O/0
   {orbit of representation of fibered group $\mathcal{G}$}%
   {orbit of representation of fibered group}%
\Symb%O/0
   {orbit of representation of the group $G$}%
   {orbit of representation of group}%

\SetIndexSpace%P%0
\Symb%P/0
   {bundle}%
   {bundle}%
\Symb%P/0
   {bundle of level $2$}%
   {bundle of level 2}%
\Symb%P/0
   {bundle of level $n$}%
   {bundle of level n}%
\Symb%P/0
   {}%
   {Cartesian power of bundle}%
\Symb%P/0
   {Cartesian product of bundles}%
   {Cartesian product of bundles, definition 1}%
\Symb%P/0
   {reduced Cartesian product of bundles}%
   {reduced Cartesian product of bundles, definition 1}%
\Symb%P/0
   {set of nonsingular \sT transformations of bundle $\bundle{}pE{}$}%
   {set of starT nonsingular transformations of bundle, projection}%
\Symb%P/0
   {set of nonsingular \Ts transformations of bundle $\bundle{}pE{}$}%
   {set of Tstar nonsingular transformations of bundle, projection}%

\SetIndexSpace%R%0
\Symb%R/0
   {active transformation}%
   {active transformation}%
\Symb%R/0
   {Cartan curvature}%
   {Cartan curvature}%
\Symb%R/0
   {\CR rank of matrix}%
   {cr-rank of matrix}%
\Symb%R/0
   {diagonal in bundle  $\bundle{}pA{}$}%
   {diagonal in bundle, 2}%
\Symb%R/0
   {diagonal in bundle $\mathcal{A}$}%
   {diagonal in reduced bundle, 2}%
\Symb%R/0
   {\Ds component of coordinates of vector $\Vector r$}%
   {Dstar component of coordinates of vector, D vector space}%
\Symb%R/0
   {curvature}%
   {GLn curvature_overline}%
\Symb%R/0
   {$\RCcirc$\Hyph product of functional matrices}%
   {rc product of functional matrices}%
\Symb%R/0
   {\RC rank of matrix}%
   {rc-rank of matrix}%
\Symb%R/0
   {right shift}%
   {right shift}%
\Symb%R/0
   {set of right-side nonsingular transformations of set $M$}%
   {set of right-side nonsingular transformations}%
\Symb%R/0
   {\sD component of coordinates of vector $\Vector r$}%
   {starD component of coordinates of vector, D vector space}%

\SetIndexSpace%S%0
\Symb%S/0
   {composition of fibered correspondences}%
   {composition of fibered correspondences}%
\Symb%S/0
   {inverse fibered correspondence}%
   {inverse fibered correspondence, 2}%
\Symb%S/0
   {inverse reduced fibered correspondence}%
   {inverse reduced fibered correspondence, 2}%
\Symb%S/0
   {linear span in vector space}%
   {linear span, vector space}%

\SetIndexSpace%T%0
\Symb%T/0
   {category of \Ts representations of $\mathfrak{F}$\Hyph algebra $A$}%
   {category of Tstar representations of F algebra}%
\Symb%T/0
   {category of \Ts representations of $\mathfrak{F}$\Hyph algebra from category $\mathcal A$}%
   {category of Tstar representations of F algebra from category}%
\Symb%T/0
   {tangent plane to group $G$}%
   {TaG}%

\SetIndexSpace%V%0
\Symb%V/0
   {coordinate vector space}%
   {coordinate vector space}%
\Symb%V/0
   {coordinates in vector space}%
   {coordinates in vector space}%
\Symb%V/0
   {direct product of $D_i\RCstar$\hyph vector spaces $\Vector V_1$, ..., $\Vector V_n$}%
   {direct product, drc vector space, 1 n}%
\Symb%V/0
   {dual space of \drc vector space $\Vector V$}%
   {dual space of drc vector space}%
\Symb%V/0
   {hermitian conjugated vector}%
   {hermitian conjugated vector}%
\Symb%V/0
   {\dcr vector space}%
   {left CR vector space}%
\Symb%V/0
   {\drc vector space}%
   {left RC vector space}%
\Symb%V/0
   {\crd vector space}%
   {right CR vector space}%
\Symb%V/0
   {\rcd vector space}%
   {right RC vector space}%
\Symb%V/0
   {tensor product of $D$\Hyph vector spaces}%
   {tensor product of D vector spaces}%
\Symb%V/0
   {tensor product of \Ds vector spaces}%
   {tensor product of Dstar vector spaces}%
\Symb%V/0
   {vector space}%
   {V}%

\SetIndexSpace%W%0
\Symb%W/0
   {set of coordinates of representation $J_f(X)$}%
   {coordinate set of representation}%
\Symb%W/0
   {set of tuples of coordinates of tower of representations $\Vector J(\Vector f,\VX X)$}%
   {coordinate set of tower of representations}%
\Symb%W/0
   {coordinates of basis $X'$ relative to basis $X$ of representation}%
   {coordinates of basis relative to basis, representation}%
\Symb%W/0
   {coordinates of element $m$ relative to set $X$}%
   {coordinates of element relative to set, representation}%
\Symb%W/0
   {tuple of coordinates of element $\Vector a*$ relative to tuple of sets $\VX X$}%
   {coordinates of element, tower of representations}%
\Symb%W/0
   {coordinates of passive map of basis manifold of representation}%
   {coordinates of passive map of basis manifold of representation}%
\Symb%W/0
   {geometrical object in vector space}%
   {geometrical object, vector space}%

\SetIndexSpace%X%0
\Symb%X/0
   {local basis of affine space}%
   {local basis of affine space}%
\Symb%X/0
   {anholonomic coordinate}%
   {x(k)}%

\SetIndexSpace%Z%0
\Symb%Z/0
   {center of ring $D$}%
   {center of ring}%

\SetIndexSpace%Delta%1
\Symb%Delta/1
   {deviation of trajectories}%
   {deviation of trajectories}%
\Symb%Delta/1
   {identical transformation}%
   {identical transformation}%
\Symb%Delta/1
   {image of vector $\Vector e_k\in\Basis e$ under isomorphism to coordinate vector space}%
   {image of vector e_k, coordinate vector space}%
\Symb%Delta/1
   {Kronecker symbol}%
   {Kronecker symbol}%

\SetIndexSpace%Gamma%1
\Symb%Gamma/1
   {anholonomic coordinates of connection}%
   {anholonomic coordinates of connection}%
\Symb%Gamma/1
   {Cartan symbol}%
   {Cartan symbol}%
\Symb%Gamma/1
   {connection}%
   {conection overline}%
\Symb%Gamma/1
   {connection coefficients in $D$\Hyph affine space}%
   {connection coefficients, D affine space}%
\Symb%Gamma/1
   {connection in $D$\Hyph affine manifold}%
   {connection, affine manifold}%
\Symb%Gamma/1
   {$D$\Hyph affine connection coefficients on manifold}%
   {D affine connection coefficients, manifold}%
\Symb%Gamma/1
   {holonomic coordinates of connection}%
   {holonomic coordinates of connection}%
\Symb%Gamma/1
   {Cartan connection}%
   {overbrace Gamma i kl}%
\Symb%Gamma/1
   {set of sections of bundle}%
   {set of sections of bundle}%

\SetIndexSpace%Lambda%1
\Symb%Lambda/1
   {inverse operator to operator $\psi_l$}%
   {inverse operator to operator psi l}%
\Symb%Lambda/1
   {inverse operator to operator $\psi_r$}%
   {inverse operator to operator psi r}%

\SetIndexSpace%Omega%1
\Symb%Omega/1
   {anholonomity object}%
   {anholonomity object}%

\SetIndexSpace%Psi%1
\Symb%Psi/1
   {}%
   {Lie Basic Operator L}%
\Symb%Psi/1
   {}%
   {Lie Basic Operator L}%
\Symb%Psi/1
   {basic operator of group Lie}%
   {Lie Basic Operator L, 1-Parameter Group}%
\Symb%Psi/1
   {}%
   {Lie Basic Operator R}%
\Symb%Psi/1
   {}%
   {Lie Basic Operator R}%
\Symb%Psi/1
   {basic operator of group Lie}%
   {Lie Basic Operator R, 1-Parameter Group}%

\SetIndexSpace%Fi%2
\Symb%Fi/2
   {Lie group composition law}%
   {Lie group composition law}%

\SetIndexSpace%Nabla%2
\Symb%Nabla/2
   {Cartan derivative}%
   {overbrace nabla_l}%
\Symb%Nabla/2
   {derivative}%
   {overline nabla_l, definition 1}%

\SetIndexSpace%Phi%2
\Symb%Phi/2
   {restriction of correspondence $\Phi$ to set $C$}%
   {restriction of correspondence}%

\SetIndexSpace%Pi%2
\Symb%Pi/2
   {Cartesian product of bundles}%
   {Cartesian product of bundles, definition 2}%
\Symb%Pi/2
   {Cartesian product of total spaces}%
   {Cartesian product of total spaces, definition 2}%
\Symb%Pi/2
   {direct product of division rings $D_i$, $i\in I$}%
   {direct product of division rings}%
\Symb%Pi/2
   {direct product of division rings $D_1$, ..., $D_n$}%
   {direct product of division rings, i 1 n}%
\Symb%Pi/2
   {direct product of $D_i\RCstar$\hyph vector spaces $\Vector V_i$, $i\in I$}%
   {direct product, drc vector space}%
\Symb%Pi/2
   {direct product of $D_i\RCstar$\hyph vector spaces}%
   {direct product, drc vector space, i 1 n}%
\Symb%Pi/2
   {product of groups $G_i$, $i\in I$}%
   {product of groups}%
\Symb%Pi/2
   {product of groups $G_1$, ..., $G_n$}%
   {product of groups, i 1 n}%
\Symb%Pi/2
   {product of objects $\{B_i,i\in I\}$ in category $\mathcal A$}%
   {product of objects in category}%
\Symb%Pi/2
   {product of objects $B_1$, ..., $B_n$ in category $\mathcal A$}%
   {product of objects in category, i 1 n}%
\Symb%Pi/2
   {reduced Cartesian product of bundles}%
   {reduced Cartesian product of bundles, definition 2}%
\Symb%Pi/2
   {reduced Cartesian product of total spaces}%
   {reduced Cartesian product of total spaces, definition 2}%

\SetIndexSpace%S%2
\Symb%S/2
   {fibered subset}%
   {fibered subset}%
\Symb%S/2
   {subbundle}%
   {subbundle}%

\CloseIndex

\end{document}